# EULER'S GRAPH WORLD
# - MORE CONJECTURES ON GRACEFULNESS BOUNDARIES-III

**Suryaprakash Nagoji Rao***

Surusha5152@hotmail.com
Jul-Sep 2014

**ABSTRACT**. Euler graphs with only one (two) type(s) of cycles under (mod 4) operation were studied in Part-I(II). Here we consider the class of Euler graphs with only three types of cycles under (mod 4). This gives rise to four cases viz., graphs having cycle types (0,1,2), (0,1,3), (0,2,3), (1,2,3). Some constructions of Euler graphs of each class are given. We prove the existence of degree two node in part cases and conjecture the existence of degree two node in every graph of order p>5. As a special case, we conjecture that regularity is nonexistent in all four cases for graphs of order p>5. That is, regular Euler graphs of order p>5 with exactly three types of cycles don't exist. In other words, a regular Euler graph of order p>5 with three types of cycles has a cycle of fourth type. An example of a nonplanar Euler graph in the first three cases is given. Graphs of this type satisfying Rosa-Golomb criterion are nongraceful. In other cases necessary conditions are given and the graphs are conjectured graceful leading to better understanding of gracefulness boundaries.

## 1. INTRODUCTION

The word *graph* will mean a finite, undirected graph without loops and multiple edges. Unless otherwise stated a graph is connected. For terminology and notation not defined here we refer to Harary (1972), Mayeda (1972), Buckly and Harary (1988). A graph G is called a *labeled graph* when each node u is assigned a label $\varphi(u)$ and each edge uv is assigned the label $\varphi(uv)=|\varphi(u)-\varphi(v)|$. In this case $\varphi$ is called a *labeling* of G. Define $N(\varphi)=\{n\in\{0,1,...,q_0\}: \varphi(u)=n$, for some $u\in V\}$, $E(\varphi)=\{e\in\{1,2,...,q_0\}: |\varphi(u)-\varphi(v)|=e$, for some edge $uv\in E(G)\}$. Elements of $N(\varphi)$, $(E(\varphi))$ are called *node (edge) labels* of G with respect to $\varphi$. A (p,q)-graph G is *gracefully labeled* if there is a labeling $\varphi$ of G such that $N(\varphi)\subseteq\{0,1,...,q\}$ and $E(\varphi)=\{1,2,...,q\}$. Such a labeling is called a *graceful labeling* of G. A *graceful graph* can be gracefully labeled, otherwise it is a *nongraceful graph* (see Rosa (1967,1991), Golomb (1972), Sheppard (1976) and Guy (1977) for chronology 1969-1977 and Bermond (1978). For a recent survey on graph labelings and varied applications of labeled graphs refer Gallion (2013).

A *(p,q)-graph* has p nodes and q edges; p and q are called *order* and *size* of the graph respectively. A *cycle graph* consists of a single cycle. $C_n$ or *n-cycle* denotes a cycle graph of length n≥3. *Euler trail* is a trail in a graph which visits every edge exactly once. Similarly, Euler circuit is Euler trail which starts and ends on the same node. *Euler graphs* admit Euler circuits and are characterized by the simple criterion that all nodes are of even degree. That is, evenness of node degrees characterizes Euler graphs. *Pendant free graph* G has node degrees two or more. Euler graphs are pendant free graphs. Regular graphs of degree >1 are pendant free graphs. *Core graph* is obtained from G by deleting all pendant nodes recursively till no pendant nodes exist. Core graph of a graph is pendant free. Core graph of a tree is empty graph. Core graph of a unicyclic graph is cycle graph. Core graph of pendant free graph is the graph itself.

By *planting a graph G onto a graph H* we mean identifying a node of G with a node of H. For a pendant free graph G of order p, a *graphforest* GF(G) is constructed from G by planting any number of trees at each node of G. GF(G) is trivial when it is of order p, ie., GF(G)=G and the trees planted are $K_1$s. We assume GF(G) is nontrivial and is of order larger than p. The smallest nontrivial graphforest has pendant node planted at a node. The choice of planting a tree at a node is random. That is, trees of any order and in any number may be planted at each node. The requirement of G a pendant free graph is for convenience. If G is not pendant free then start with core graph of G. Note that for any graph G the simplest graph structure containing G of higher order is graphforest GF(G). A graphforest is *graphtree* when one nontrivial tree is planted. Graphforest of Euler graph G is called *Eulerforest* EF(G). *Cycleforest* CF($C_n$) and *Treeforest* TF(T) are similarly defined. The class of cycleforests is precisely the class of unicyclic graphs. Note that graphforest of a forest is forest and graphforest of a tree is a tree.

# EULER GRAPHS

Several properties of Euler graphs are known including characterizations (e.g., Zsolt (2010)).

Gracefulness exhibits affinity with Euler graphs as in:

**Theorem** G (Rosa (1967), Golomb (1972)). A necessary condition for Euler $(p,q)$-graph to be graceful is that $[(q+1)/2]$ is even.

Euler graph of this type is called here a *'Rosa-Golomb graph'*. This implies that Euler graphs with $q \equiv 1$ or $2 \pmod 4$ are nongraceful. Most of the known nongraceful graphs contain a subgraph isomorphic to Rosa-Golomb graph. The necessary condition for Euler graphs leads to understanding the intrinsic cycle structural properties of graceful Euler graphs. We consider four subclasses of Euler graphs with only three types of cycles. We denote by $\varepsilon_{ijk}$ the class of Euler graphs with only cycles $C_n$, $n \equiv i,j$ or $k \pmod 4$, $i,j,k = 0,1,2,3$; with at least one cycle of each of the three types. For example, $\varepsilon_{012}$ is the class of Euler graphs having only cycles of the type $C_n$, $n \equiv 0,1\&2 \pmod 4$.

# ODD CHORDS

An *odd chord* of a cycle in a graph is an edge which joins two nodes of a cycle and divides it into two parts such that one part is an odd cycle. A chord which is not odd is called even, in other words, an *even chord* of a cycle divides it into two even parts. A cycle without chords is a *chordless cycle*. A bipartite graph has no odd chords. An edge added within a bipartition of a bipartite graph results in odd chords.

**Observation 1.** A chord of a cycle $C_n$, $n \equiv 2 \pmod 4$ divides it into two even cycles of different parity. That is, if one is $\equiv 2 \pmod 4$ then the other is $\equiv 0 \pmod 4$. Whereas, a chord of $n \equiv 0 \pmod 4$ cycle divides into two even cycles of same parity. That is, if one is $\equiv 0 \pmod 4$ ($\equiv 2 \pmod 4$)) then the other also is $\equiv 0 \pmod 4$ ($\equiv 2 \pmod 4$)).

**Observation 2.** A chord of a cycle $C_n$, $n \equiv 1 \pmod 4$ ($\equiv 3 \pmod 4$)) divides it into an odd cycle and an even cycle. That is, if odd cycle is $\equiv 1 \pmod 4$ ($\equiv 3 \pmod 4$)) then the other is even and is $\equiv 2 \pmod 4$ ($\equiv 0 \pmod 4$)).

# CYCLE INTERSECTION

**Theorem 1.** If two cycles intersect in $m>0$ multiple paths $P_t$, $t>0$ of different lengths in common then removal of the common paths results in $m+1$ disjoint cycles. Each such cycle is made up of two parts one from each of the two cycles.

**Proof.** We shall prove this for two intersections, with each intersection having at least one edge as shown in Fig.1. The general case follows similarly. The cycles $C_1$ and $C_2$ have two intersections $I_1$ and $I_2$ with at least one edge in common. Removal of the common paths except for the end nodes from the graph results in three cycles $NC_i$, $i=1,2,3$. The cycles clearly are made up of two parts one of $C_1$ and the other $C_2$.

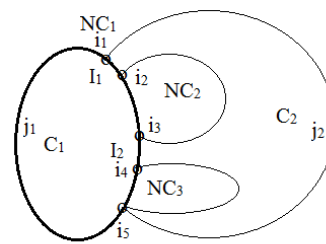

Fig.1 Cycles with two (multiple) intersections.

Hereafter, by *two intersecting cycles* means they intersect only once. That is, removal of the path of intersection results in only one cycle consisting of two parts, one part each from the two cycles.

## GRACEFULNESS BOUNDARIES THROUGH EULER GRAPHS

### Euler Graphs with Only Three Types of Cycles

Four cases arise according as Euler graph has only three of the four cycle types: (0,1,2), (0,1,3), (0,2,3), (1,2,3). A graph with each block being a i-cycle graph, i=0,1,2,3 is Euler. Such a graph with three distinct types of cycles i,j,k from {0,1,2,3} belongs to $\varepsilon_{ijk}$.

For a given cycle decomposition of a Euler graph G, classify the cycles $C_n$ of G into four classes $\mathcal{C}_i$ - the class consisting of cycles $C_n$, n≡i(mod 4), i=0,1,2,3. Denote by $\xi_i=|\mathcal{C}_i|$, number of cycles in $\mathcal{C}_i$ (see Rao (2014), [10]).

Problem. What is the minimum order and size of Euler graph in each of these four cases?

### Case 0,1&2. $\varepsilon_{012}$: Euler graphs with only three types of cycles $C_n$, n≡0,1&2(mod 4)

**Observation 3.** We claim that the order of graphs in $\varepsilon_{012}$ satisfies p≥8. Such a graph has a 6-cycle and so p≥6. Any such 6-cycle is chordless by Observation 1. The graph has cycles of type (0,1,2). Any 5-cycle in it also is chordless and has a maximum of three edges in common with that of 6-cycle. This makes the graph of order at least 7 with two nodes of degree three each. For the graph to be Euler each of these nodes is adjacent to another node making the graph of order at least 8. The graph in Fig.2a is an example of Euler (8,10)-graph with only cycles of type (0,1,2). By the argument above this is the smallest Euler graph from $\varepsilon_{012}$.  Note that this graph is a nongraceful graph.

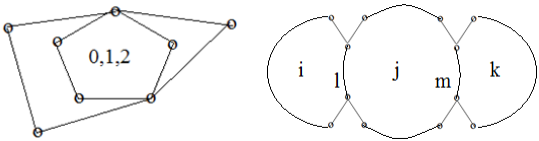

Fig.2 An example of Euler graph from $\varepsilon_{012}$ and three intersecting cycles.

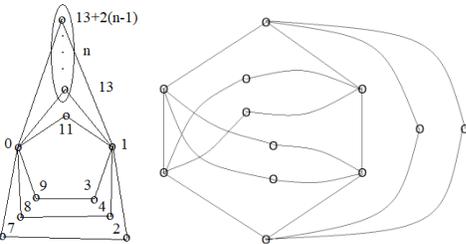

Fig.3 A graceful Euler graph and a nonplanar Euler graph from $\varepsilon_{012}$.

Three intersecting cycles i,j,k are shown in Fig.2b with each intersection of length >0. A graceful Euler graph from $\varepsilon_{012}$ is shown in Fig.3a. Graphs from $\varepsilon_{012}$ need not be planar. A nonplanar graph from $\varepsilon_{012}$ is given in Fig.3b.

| Table-1 Case 0,1,2 | | | |
|---|---|---|---|
| Cycle Type | | Combined Cycle type | |
| Cycle-1 | Cycle-2 | Int Even | Int Odd |
| 0 | 0 | 0 | 2 |
| 0 | 1 | 1 | 3 |
| 0 | 2 | 2 | 0 |
| 1 | 1 | 2 | 0 |
| 1 | 2 | 3 | 1 |
| 2 | 2 | 0 | 2 |

The Combined Cycle (CC) rule for even or odd intersection of two cycles is tabulated in Table-1. For example, two cycles of type 0 intersecting in odd path result in a combined cycle type 2. Two cycles of type 1 and 2 intersecting in even path results in a CC type 3 and so is open which in this case 0,1,2 is fourth type shown in yellow. The intersection of cycle types 0 and 1 with odd intersection also yields a CC of type 3 and so open. In other cases the CC is closed, that is it is of type 0,1 or 2.

**Theorem 2.** The CC is of type 3 whenever two cycles of type 0 and 1 (1 and 2) of a graph G from $\varepsilon_{012}$ intersect in a path of length l>0 is odd (even).

**Theorem 3.** Size of a graph from $\varepsilon_{012}$ satisfies $q \equiv \xi_1 + 2\xi_2 \pmod 4$.

Proof follows from Equation (1) (Rao 2014 [10]). Further, if $\xi_1 + 2\xi_2 \equiv 1 \text{ or } 2 \pmod 4$ then the graphs are nongraceful; else $\xi_1 + 2\xi_2 \equiv 0 \text{ or } 3 \pmod 4$ and the graphs are candidates for gracefulness. Graphs in $\varepsilon_{012}$ satisfy that $\xi_0 > 0$, $\xi_1 > 0$ and $\xi_2 > 0$.

**Theorem 4.** Two necessary conditions follow: If $\xi_1 + 2\xi_2 \equiv 0 \pmod 4$ then $\xi_1$ is even. If $\xi_1 + 2\xi_2 \equiv 3 \pmod 4$ then $\xi_1 - 3$ is even or $\xi_1$ is odd.

**Conjecture 1.** Graphs from $\varepsilon_{012}$ satisfying $\xi_1 + 2\xi_2 \equiv 0 \text{ or } 3 \pmod 4$ are graceful.

**Three Intersecting Cycles**

Every Euler graph G from $\varepsilon_{ijk}$, i,j,k=0,1,2,3 has two intersecting cycles $C_i$, $C_j$, i≠j. Let u-v path of length >0 be the intersection. Then there is at least one node u' adjacent to u. Eight cases arise as in Fig.4 (a) to (h). Case (a) is same as Case (b) with u'=v. Further, Case (c) and Case (d) are same. If there is a node w≠u' without a w-u or w-v path then there is a node w' as shown in Figs.4(e),(f),(g),(h). We consider the five cases Figs.4(a), (b), (c), (e), (f) for further analysis. Note here that Case (e) is same as Case (c) with w'=u. Cases (c) and (d) are same. Case (h) includes Cases (a) and (b). The three cycles are denoted by $C_1, C_2, C_3$ as shown in Fig.4 Throughout the paper the label i, i=0,1,2,3 with big font size in the figures for the cycles indicates a cycle $C_i$, i≡0,1,2 or 3 (mod 4). The intersection path length is denoted by l,m,n. The intersection path denoted by d or e with small font size denotes odd or even parity of the path.

An Euler graph which is not a cycle graph has nodes of even degree with at least two nodes of degree ≥4. By this argument, G has a node u'' adjacent or a path to u'. We consider all possible cases for u''. This introduces new cycle in the resulting graph. If this cycle is any one of types i,j,k then we proceed further to the next level else the case stops. That is, the resulting graph has a cycle not of the type i,j,k and so the graph has four types of cycles, a contradiction. E.g., if G is from $\varepsilon_{012}$ then a cycle of the type ≡3 (mod 4) exists by adding a new node. In the following figures such a cycle is shown either by thick path or blue or black color. In some cases, the second level analysis is also shown. Continue cases at first level lead to second level. The new node and the new path containing it is selected so as to minimize the number of second level cases. In general, the subcases keep increasing from second level depending on the continue cases in the first level. We have two cases where the second level subcases all stop.

The path adding a new node creates new cycles. An existing cycle may get divided into two cycles by this path. The possible types for these cycles come from the Table-1. More new cycles get created in combination with these two cycles and existing cycles. Where possible, parity of intersection paths may be determined. There may be intersection paths with new cycles not satisfying parity of intersection path making the case stop. Such paths are shown in orange color.

*Euler's Graph World - More Conjectures on Gracefulness Boundaries-III*

MGB-Cijk-3C-20Jan16

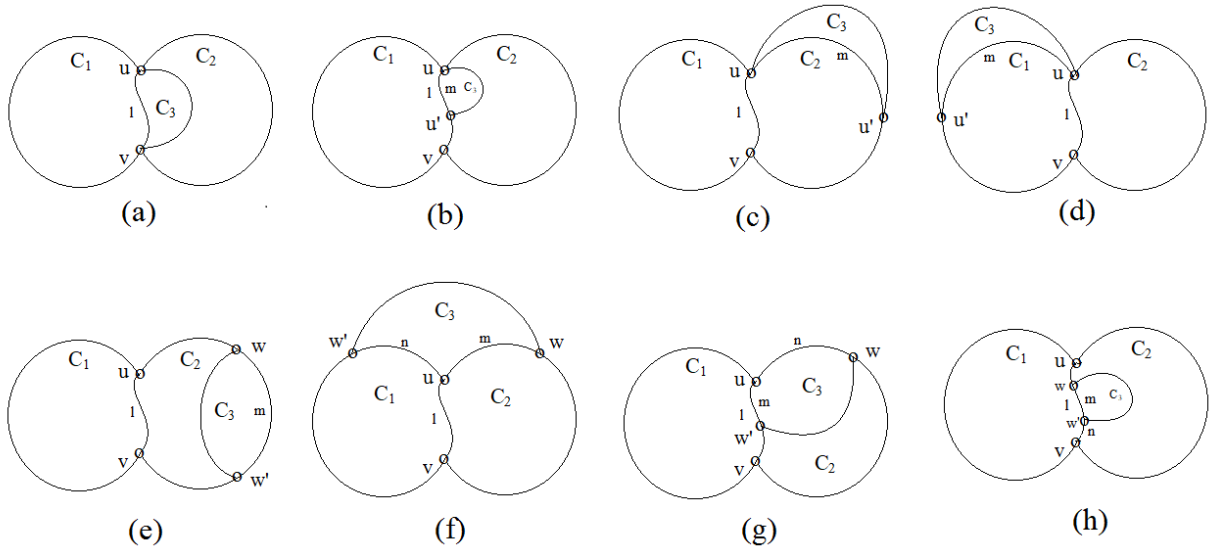

Fig.4 Three intersecting cycle configurations.

**Observation 4**. If three cycles A,B and C in a Euler graph intersect $A \cap B$, $B \cap C$, $A \cap C$, $A \cap B \cap C$ are nonempty. Then $A \cap B \cap C$ is a node x and each of $A \cap B$, $B \cap C$, $A \cap C$ is a path containing x. Without loss to generality the cycles $C_1, C_2, C_3$ may be chosen such that the u-x, v-x and w-x paths are minimum.

In general, we conjecture that regular graphs are sparse in $\varepsilon_{012}$. We prove this by verifying the existence of a node of degree two in some subcases in the cases (a), (b), (c), (e) and (f). It is hoped that the cases at some level stop completely with no continue case.

**Conjecture 2**. There exists a degree two node in every graph G of order p>7 from $\varepsilon_{012}$.

If this conjecture is true then as a corollary it follows that:

**Conjecture 3**. Regular Euler graphs G of order >7 in $\varepsilon_{012}$ are nonexistent.

or

A regular Euler graph from $\varepsilon_{012}$ of order p>7 with three types of cycles also contains fourth type.

**Theorem 5**. There exists degree two node in a graph G of order p>7 from $\varepsilon_{012}$ in the subcases which stop.

Part Proof. We consider part proof in the cases (a),(b),(c),(e) and (f). The proof technique used is to show the existence of fourth cycle type or existence of two cycles which intersect with parity leading to fourth cycle type as in Table-1. The cases (a), (b), (c) and (d) consider the u-u' path from the node u, in the first case u'=v. For proof, it suffices to consider the first three cases.

**Case (a)** No case is possible.

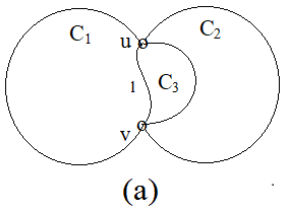

**Case (b)** u'≠v

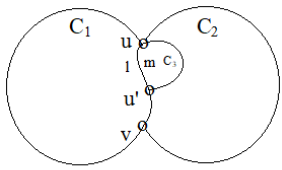

(b)

Three intersecting cycles may be as shown in Case (b) above with one cycle having partial intersection with the common path. Here, u'≠v. In this case, the possibilities leading to closed cycles are listed in the Table-2 for each of the two cases. Apart from the three cycles, depending on the intersection there are three combined cycles as shown. Since each cycle may be one of the four types under (mod 4), these combined cycles (mod 4) are shown in the last three columns as in the Table-2. Eliminating the cases containing all fourth type of cycle from the last three columns resulted in the two cases as shown in Table-2. Note that the cases C1, C2 and C2, C1 are same in this case. For each combination an example of Euler graph is shown in the Fig.5.

Table-2 Possible cases in Case (b).

| i | j | k | l | m | i+j-2l | i+k-2m | j+k-2m |
|---|---|---|---|---|--------|--------|--------|
| 0 | 1 | 2 | 0 | 1 | 1 | 0 | 1 |
| 1 | 2 | 0 | 1 | 0 | 1 | 1 | 2 |

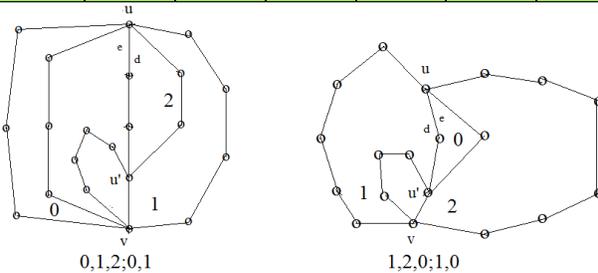

Fig.5 Euler graph examples for each of the two cases in Table-2.

Two representations in drawing the figures are considered in this case as in Figs.6 and 7. In the following figures, e,d stand for even, odd respectively. For example, here e is the even parity of u-v path and d is the odd parity of u-u' path. The node u' is of degree 3 and so there is a node adjacent to u'. In general, we consider possible paths at u' say u'-w paths. There may be other possibilities for u' such as the node v. The purpose of this choice is to minimize possible cases and the cycle division by the chosen path leads to maximum stop cases.   A case is '*stop case*' if it contains fourth type of cycle otherwise it continues to next level.

The Fig.6 considers possible such paths. In the first figure of Fig.6, the u'-w path divides the cycle 1 into two cycles. These cycles may be 0,1,2 type. Possible cycles are 0 and 1 with u'-w path parity e; 1 and 2 with odd parity. On the right, possible cycles with parity is listed for each cycle 0,1,2 from the rule Table-1. Some of these cases stop as in second figure of Fig.6, the case 0,2;d. It is shown as x against the case. The cases with x stop, other cases continue to next level.

MGB-C3-case2-01201-19Oct16

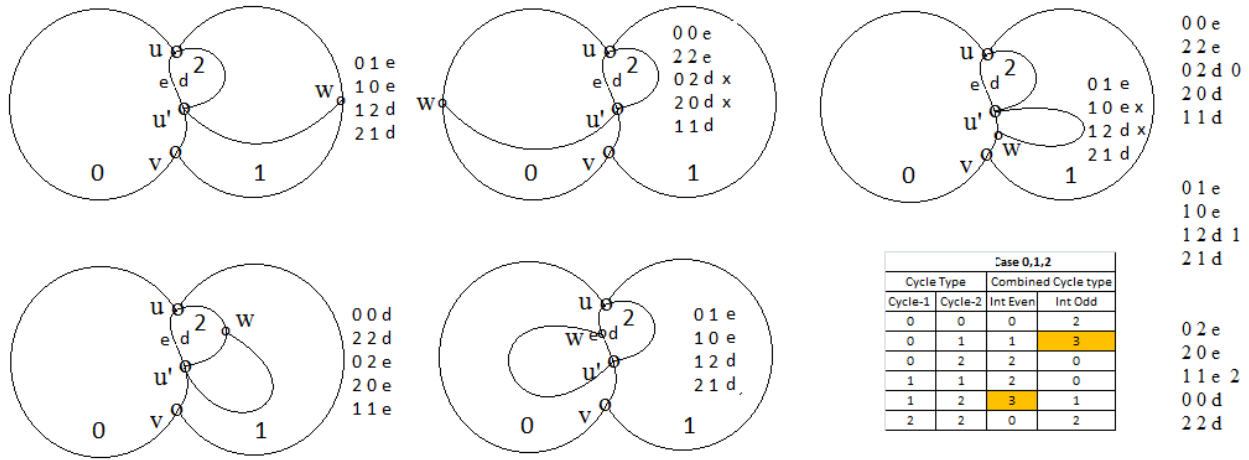

Fig.6 The five first level possibilities for a new node w and u'-w path with cycle division in the case 0,1,2;e,d.

Using second representation, there are two possibilities for i,j,k with $C_i$ and $C_j$ intersecting in u-v path and $C_k$ intersecting $C_i$, $C_j$ with w-v path which is a part of u-v path in common. $C_0$, $C_1$, $C_2$ cycles are as shown in Fig.4b. Consider the first case l (e-even) and m (d-odd) parity which is equivalent to the case 0,1,2;e,d. Five subcases arise for u'-x path as shown in Fig.7.

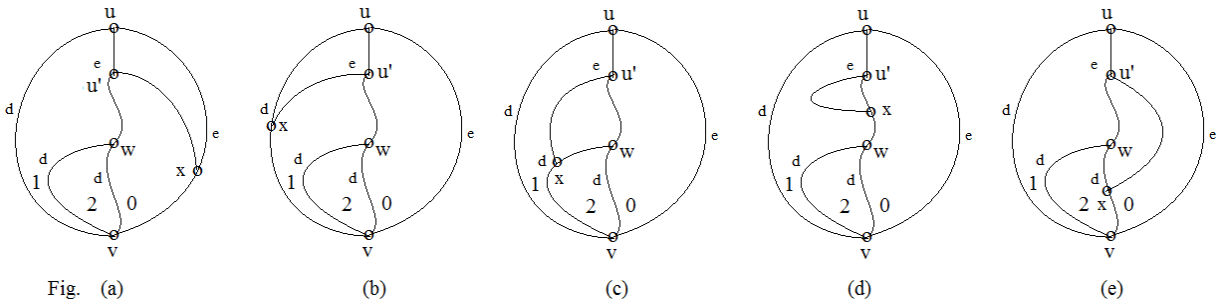

Fig.7 The five first level possibilities for the u'-x path.

Consider the node u' such that uu' is an edge existence comes from the uu'-v path. If u' is of degree 2 we are done else there is a node x so that u'-x path exists. This gives rise to five cases as shown in Figs.6 and 7 for the first level. In the subcase Fig.7a, the path u'-x divided the cycle $C_0$ into two cycles uu'-x-u; u'-w-v-x-u'. From Table-1 we can find possibilities for these cycles with parity. Fig.8 shows 20 first level possibilities for u'-x path with cycle division and stop cases for the case 1,0,2;e,d.

Note that, the choice for the cycle for u'-x path may not be unique. Some subcases stop as one of the new cycles being the fourth type. In this case, it is cycle type 3. It is shown in blue or black or bold paths. Paths with orange color indicate a cycle of type 3 containing end nodes of the path and so a stop case. For example, Fig.8b in the second row has orange path as odd intersection of cycle types 1 and 0 resulting in cycle type 3. At times the possible new paths are shown in blue color or bold paths. The u'-x path may divide one or more cycles. For example, Fig.7d may be considered dividing cycle 1 or cycle 0. Accordingly, the possibilities occur from Table-1. We look for the u'-x path so that no cases continue or minimum cases continue.

Similar analysis is accomplished throughout the paper. Analysis for the first or second level is presented. Some cases stop with the presence of a cycle of fourth type. The continue cases lead to more cases as the level increases.

MGB-C102CI-01-01e-11Nov15

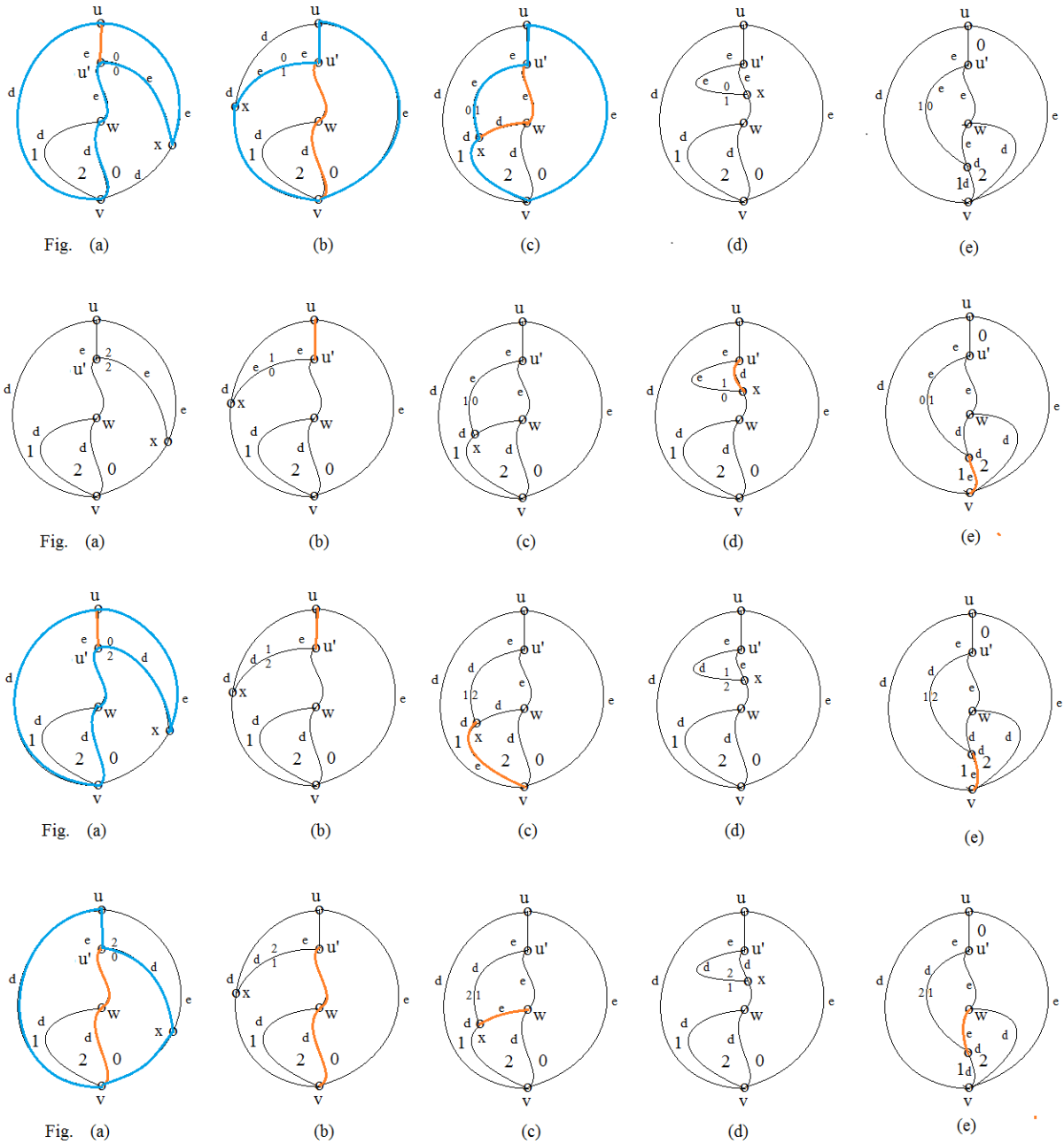

Fig. 8 The 20 first level possibilities for u'-x path with cycle division and stop cases for the case 1,0,2;e,d.

Fig.8a and Fig.8c in second row; Fig.8d in rows 3 and 4 continue for the second level. Second levels are shown in the first two cases. Considering Fig.8a, u'-x is even type and divides the cycle 0, uu'-w-x-v-u into cycle types 2,2;e. This case continues further to next level as all cycle types are 0,1,2. The path u'-x is even implies there is a node y so that u'y is an edge. If y is of degree two then we are done. If not, Fig.9 shows the 31 possibilities of a y-y' path to take the case 1,0,2;e,d;2,2,e further for second level for verification whether a case stops or continues with cycles of the type 0,1,2.

MGB-C102CI-01-22e-06Dec15

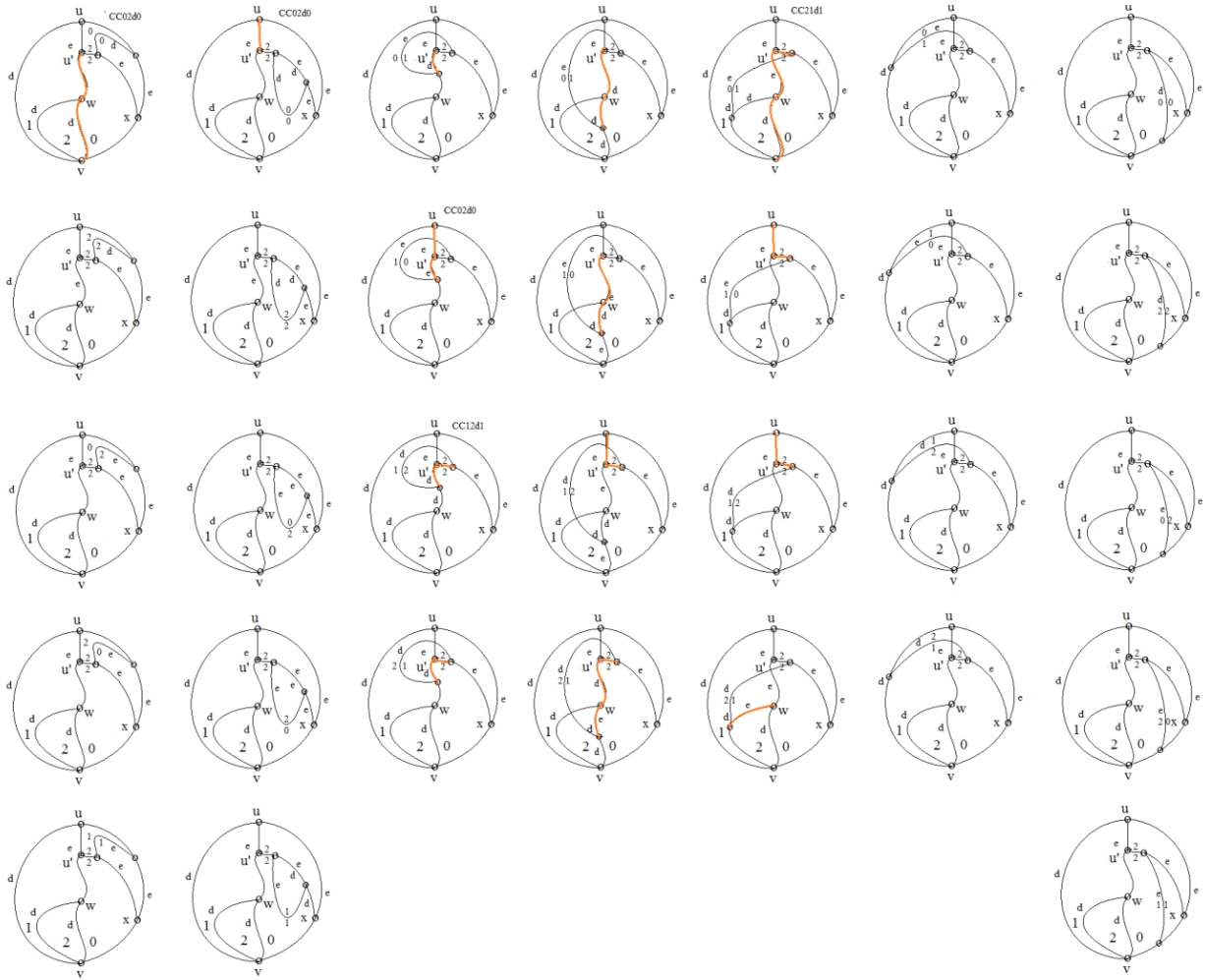

Fig.9 The 31 second level possibilities for u'-x path with 2,2,e division with cycle division and stop cases.

Fig.8c in second row u'-x is even type and divides the cycle uu'-w-x-v-u into cycle types 1,0;e. This case continues further to next level as all cycle types are 0,1,2. The path u'-x is even implies there is a node y so that u'y is an edge. If y is of degree two then we are done. If not, Fig.10 shows the 32 possibilities of a y-y' path to take the case further for second level for verification whether a case stops or ends in cycles of the type 0,1,2 to continue to third level.

MGB-C102CI-01-10e-11Nov15

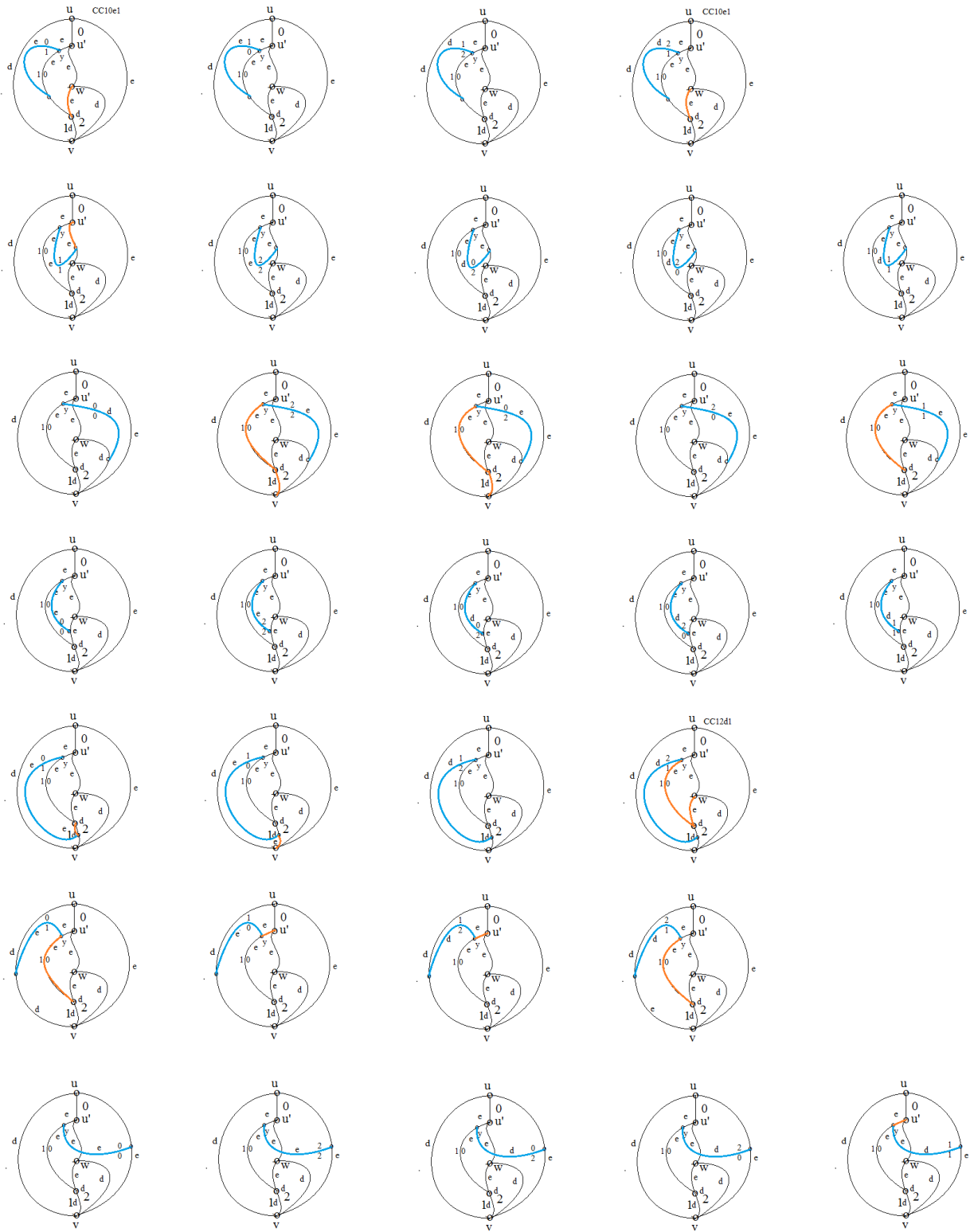

Fig.10 The 32 second level possibilities at the node y with division of cycles and stop cases in case 1,0,2;e,d;1,0,e.

MGB-C120CI-10-06Dec15.jpg

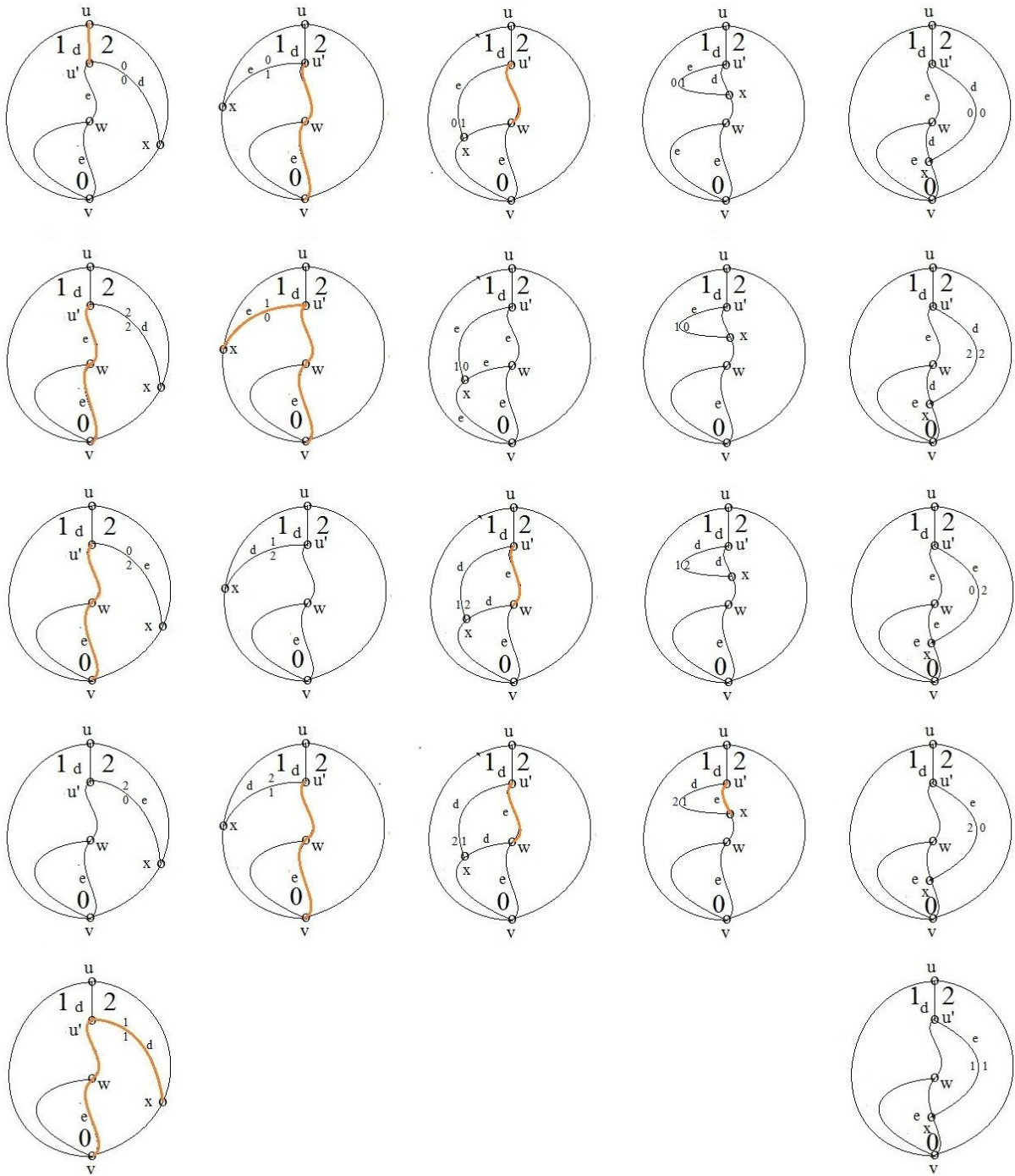

Fig. 11 The first level 22 possibilities for u'-x path at u' with cycle division and stop cases in the case 1,2,0;d,e.

**Case (c)**

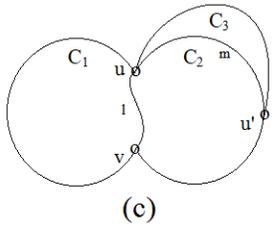

(c)

Table-3 Possible cases for the Case (c).

| i | j | k | l | m | i+j-2l | j+k-2m | i+j+k-2l-2m |
|---|---|---|---|---|--------|--------|-------------|
| 0 | 1 | 2 | 0 | 1 | 1 | 1 | 1 |
| 0 | 2 | 1 | 0 | 1 | 2 | 1 | 1 |
| 1 | 0 | 2 | 0 | 1 | 1 | 0 | 1 |

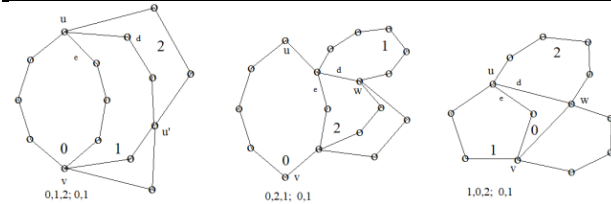

Fig. 12 Euler graph form $\varepsilon_{012}$ for each case in the Table-3.

MGB-C3-case1-01201-19Oct16

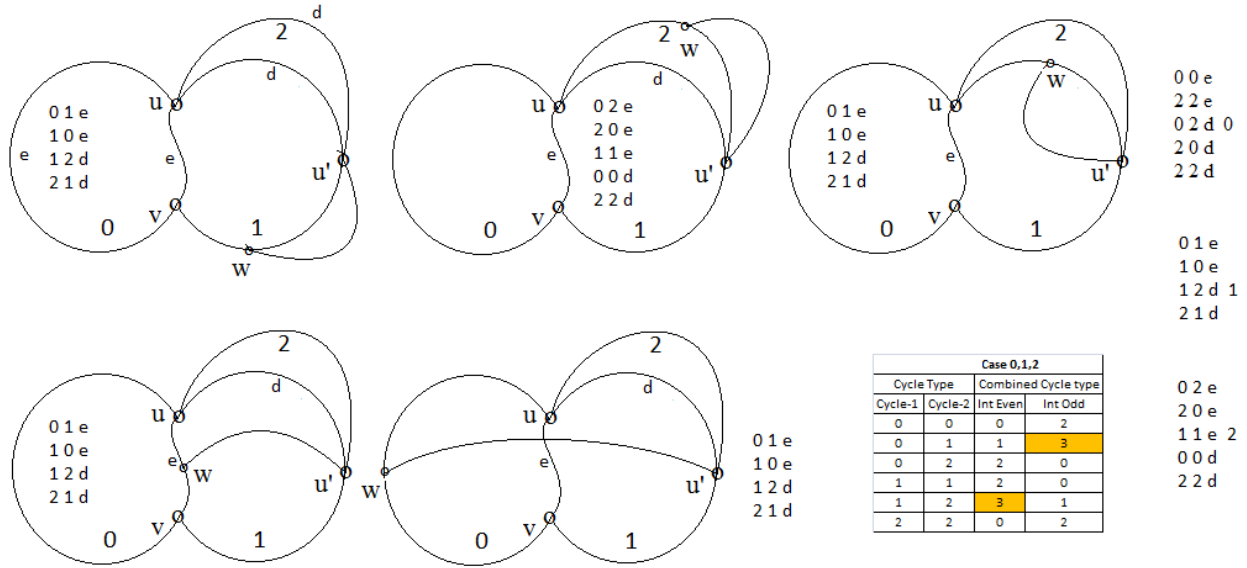

Fig. 13 The five first level possibilities for u'-w path at u' for the case 0,1,2;e,d.

MGB-C3-case1-02101-19Oct16

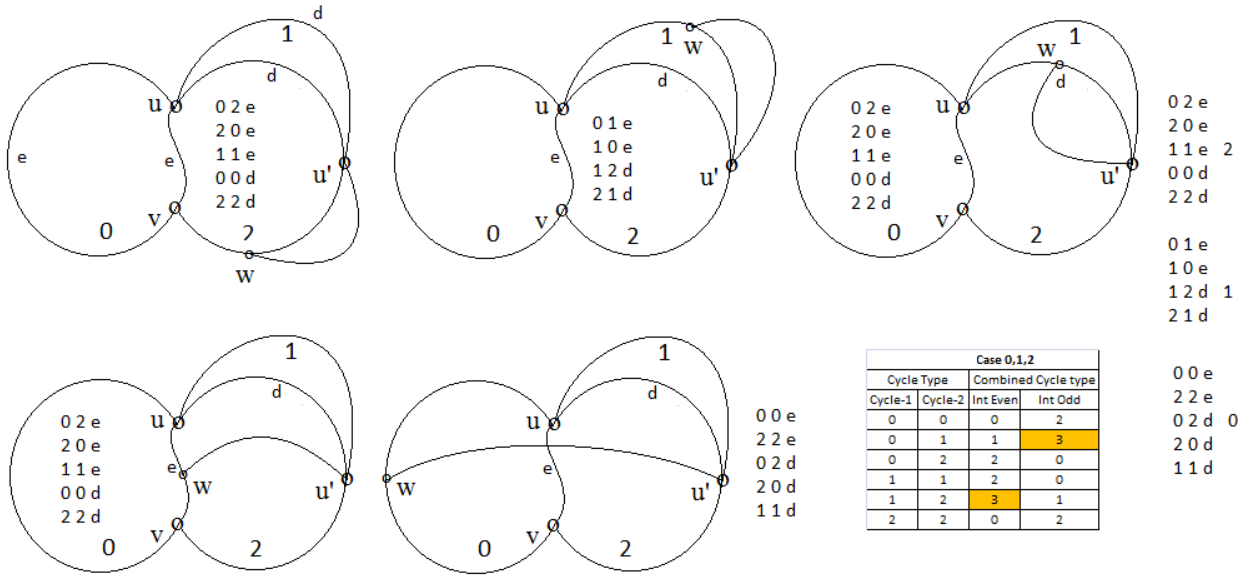

Fig. 14 The five first level possibilities for u'-w path at u' for the case 0,2,1;e,d.

MGB-C3-case1-10201-19Oct16

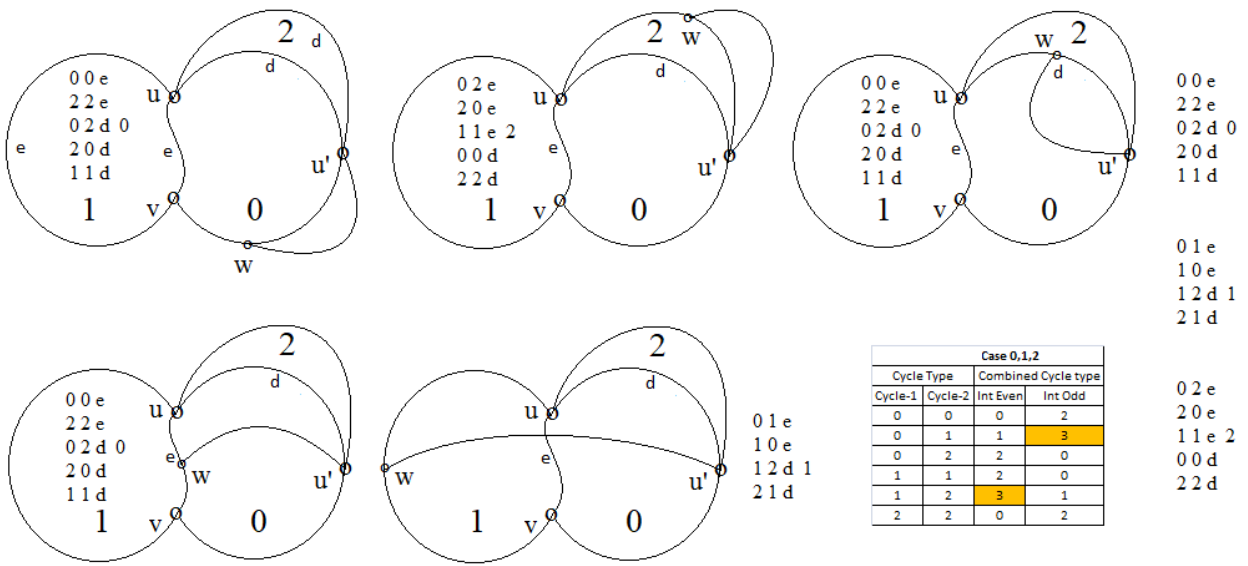

Fig. 15 The five first level possibilities for the u'-w path for the case 1,2,0;e,d

**Case (d)**  Same as Case c.

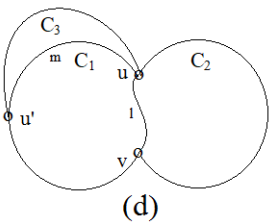

(d)

**Case (e)**

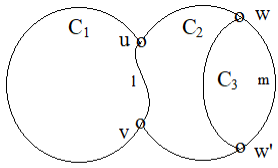

(e)

Table-4 Possible cases for the case (e).

| i | j | k | l | m | i+j-2l | j+k-2m | i+j+k-2l-2m |
|---|---|---|---|---|--------|--------|-------------|
| 0 | 1 | 2 | 0 | 1 | 1 | 1 | 1 |
| 0 | 2 | 1 | 0 | 1 | 2 | 1 | 1 |
| 1 | 0 | 2 | 0 | 1 | 1 | 0 | 1 |

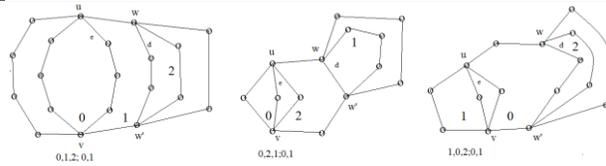

Fig.16a Euler graph form $\varepsilon_{012}$ for each case in the Table-4.

MGB-Cijk-3CC-012-ed-14Jun20

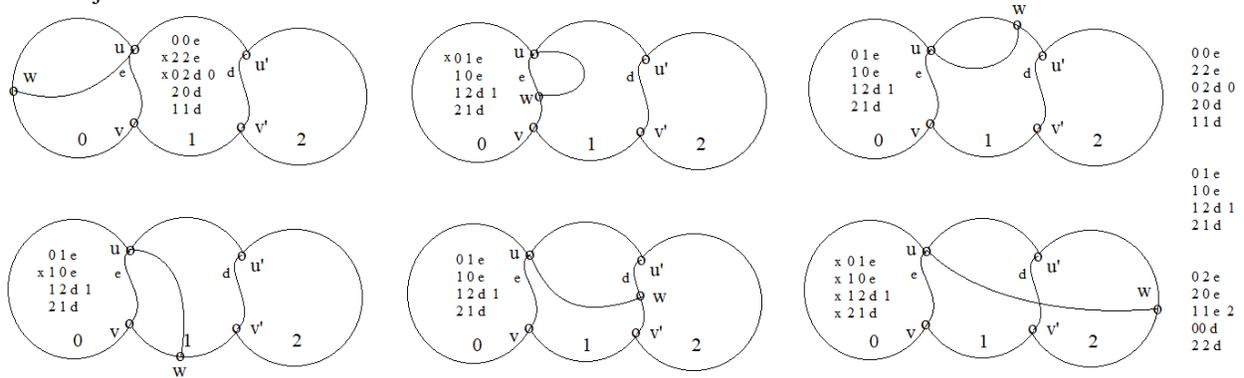

Fig. 16b The six first level possibilities for u-w path with cycle division and stop cases in the case 0,1,2;e,d.

**Case (f)**   Graphs with three intersecting cycles, any two intersecting in different paths

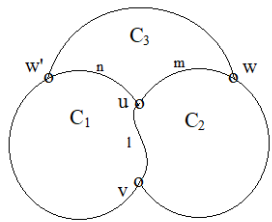

(f)

Table-5 Possible cases in the case (f).

| i | j | k | l | m | n | i+j-2l | j+k-2m | i+k-2n | i+j+k-2l-2(m+n) |
|---|---|---|---|---|---|--------|--------|--------|-----------------|
| 0 | 1 | 2 | 0 | 1 | 0 | 1 | 1 | 2 | 1 |
| 0 | 2 | 1 | 0 | 1 | 0 | 2 | 1 | 1 | 1 |

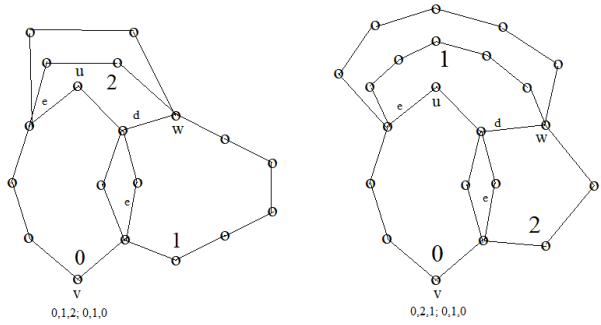

Fig.17 Euler graph form $\varepsilon_{012}$ for each case in the Table-5.

The v-x, u-x and w-x path lengths (mod 4) are in 1,2; 3,1; 3,3 cells. The cells (i,j) = 2,1; 3,2; 2,3; are i+j-l(v-x)-l(u-x) (mod 4).

MGB-C012-10-eed-16Jul15

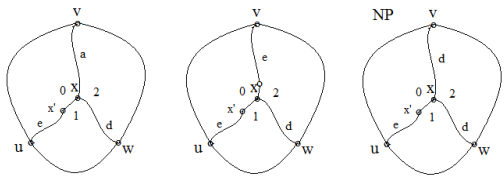

Fig.18 General case aed followed by the cases eed and ded. The ded case is not possible.

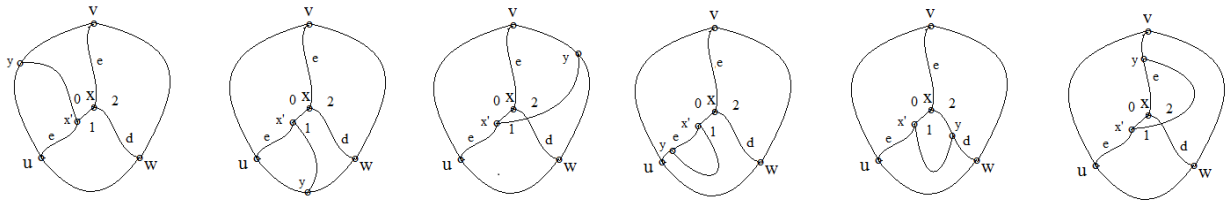

Fig.19 The six first level possibilities for x'-y path in the eed case.

Consider a graph G from $\varepsilon_{012}$. By definition G has at least one cycle of type $T_0$, $T_1$ and $T_2$. G has two intersecting cycles $C_i$ and $C_j$ i≠j; i,j=0,1,2. Two cases arise according as G has three distinct intersecting cycles $C_0$, $C_1$ and $C_2$ or not.

Firstly, let the cycles $C_0$, $C_1$, $C_2$ intersect in the node x and the common u-x, v-x, w-x paths as in Fig.18. From Table-1 u-x path is even, w-x path is odd and v-x path may be even or odd as in Fig.18a. Table-1 shows types of cycles including combined cycle resulting from u-v, v-w, w-u paths. $C_0$ and $C_2$ cycles may have even or odd intersection. Further, the u-x and w-x paths are even and odd respectively follows leading to two cases as in Fig.18b,c. The case Fig.18c is impossible follows as the combined cycle v-u-w-v is of type $T_3$.

Consider the case Fig.19a eed. Since u-x path is even there is a node x' adjacent to x. If x' is of degree 2 then we are done. Otherwise, there is x'-y path. There are six cases for y as in Fig.19a-f. In all these cases we shall show stop cases by the existence of a cycle type $T_3$ or existence of two cycles with intersection not as in Table-1.

**Euler's Graph World - More Conjectures on Gracefulness Boundaries-III**

MGB-C012-eed-16Jul15

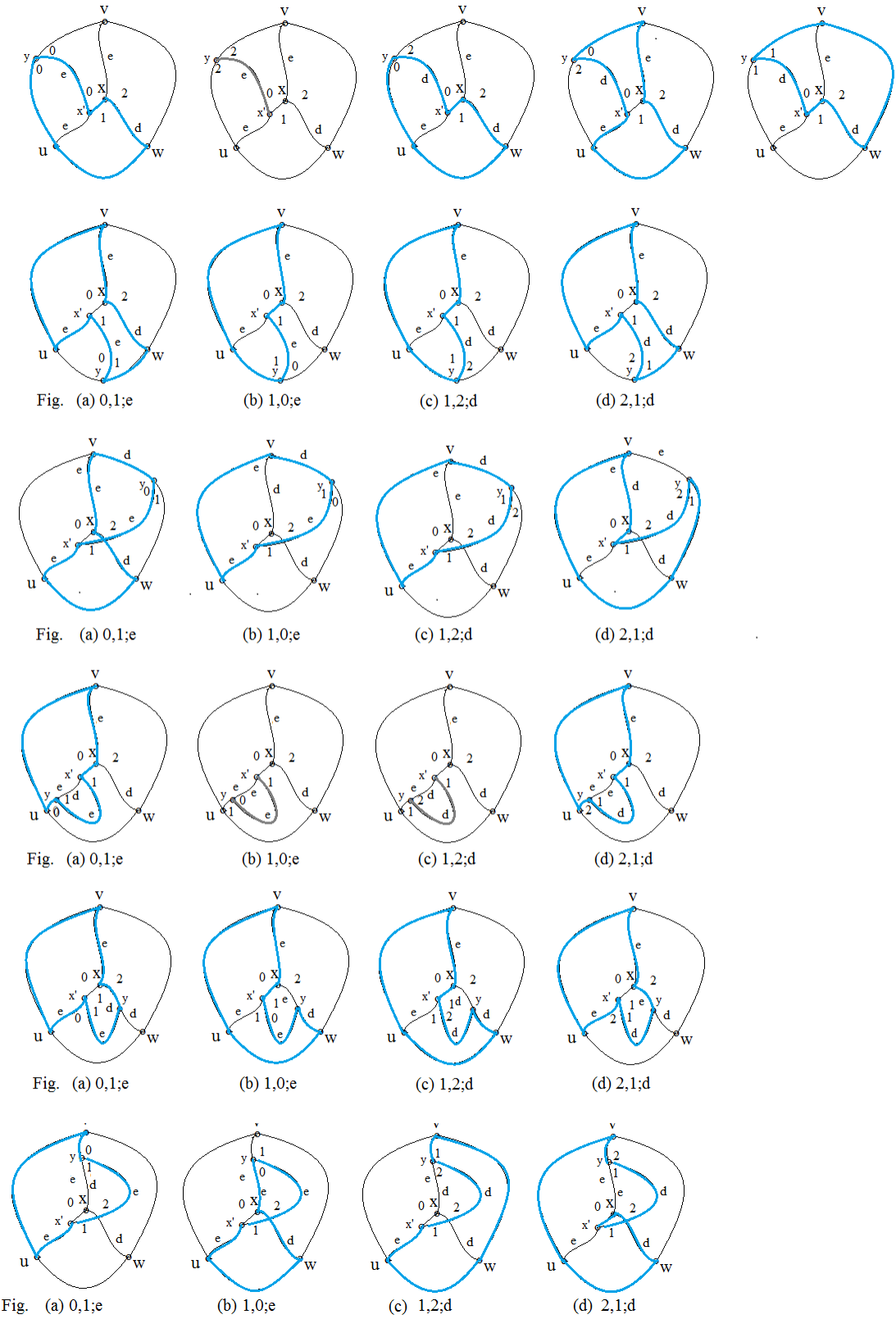

Fig. 20 The 25 first level possibilities for x'-y path with dividing cycles and stop cases.

In case-a, $C_0$ gets divided into two cycles u-x'-y-u and v-x,x'-y-v. From Table-1 the possibilities are 0,0;e, 0,2;d, 2,0;d, 1,1;d, 2,2;e. Here, the triple i,j;p stands for intersecting cycles of type $T_i$, $T_j$ with path of intersection of length l>0 with parity p. First four cases are impossible as two cycles of type $T_0$ and $T_1$ exist with odd intersection.

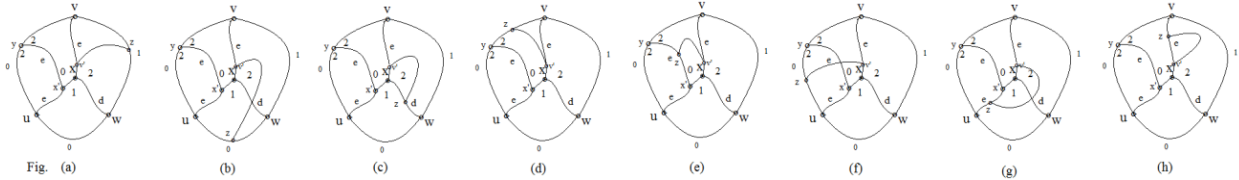

Fig. 21 The eight second level possibilities for v'-z path from x-v path in 2,2;e case.

The case-b 2,2,e leads to eight subcases as in Fig.21 (a) to (h).
The cases in Fig.21a,b,c,g,h are impossible as the v'-z path divides type $T_1$ cycle into 0,1;e; 1,0;e; 1,2;d; or 2,1;d.

The cases d,e lead to division of cycle type $T_2$ into 0,0;e; 0,2;d; 2,0,d; 1,1,d; 2,2,e. Impossibilities of these cases follow from the intersections.

The case f leads to division of cycle type $T_0$ into: 0,0;e; 0,2;d; 2,0;d; 1,1;d; 2,2;e.

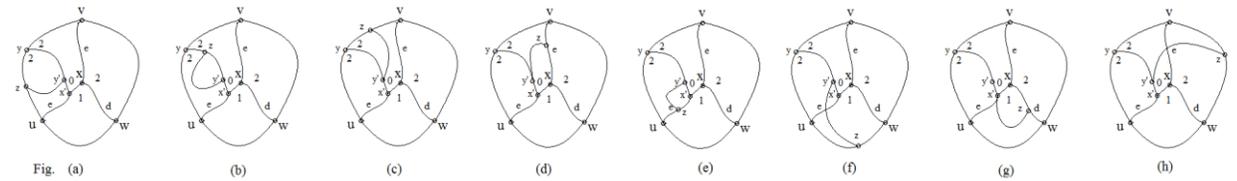

Fig.22 An alternative eight possibilities for y'-z path from x'-y path in 2,2;e case.

Fig.22 Case-e long x'-y path divides $C_2$ into two cycles: x'-y(short)-x'(long), x',x,w-u-y-x' with possibilities: 0,1;e, ruled out as $T_1$ and $C_1$ intersect in even x'-y long path with CC of type $T_1$. So short x'-y path is odd. But then, this $T_1$ cycle intersects with $C_3$ in odd v-x path from Table-1, implying v-x path is both even and odd, a contradiction.
1,0;e,
1,2;d,
2,1;d is ruled out on the similar lines that of 0,1;e.

In each of the other four Fig.22 cases-c,d,f,g there is a cycle type $T_1$ or combined cycle of $C_1$ and $C_2$ with even intersection resulting in cycle type $T_1$. So the x'-y path divides this cycle of type $T_1$ into two cycles. From Table-1, the possibilities are 0,1;e; 1,0;e; 1,2;d; and 2,1,d. In each of these four cases cycle type $T_3$ exists and so impossible.

In Fig.23 the x'y'-y path is even and so there is an edge x'y'. Consider a path at y' and its possibilities are shown in Fig.23. Also shown are a cycle it divides and the possibilities for the cycles from Table-1. Some cases stop as shown in blue or orange color.

MGB-C012-10-eed-16Jul15

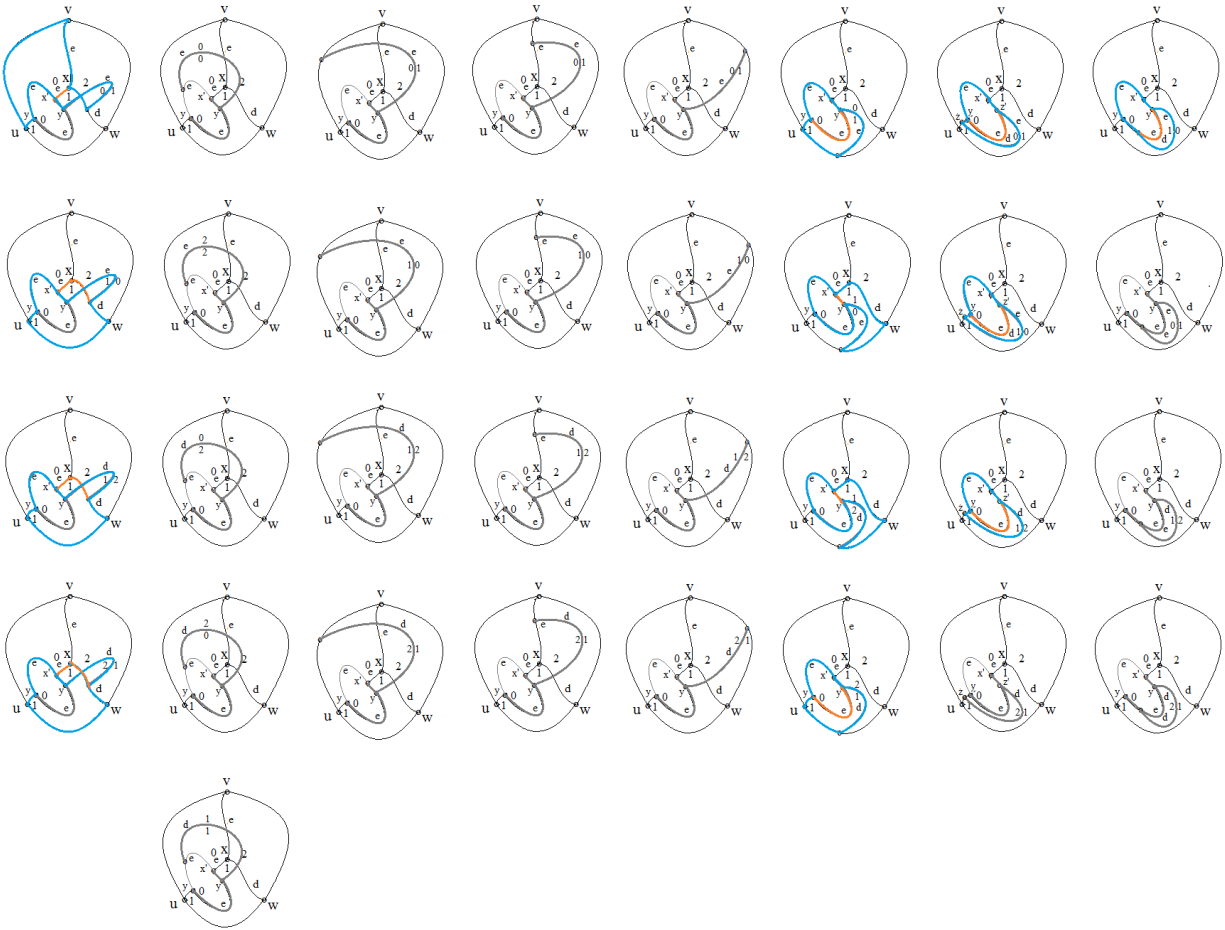

Fig. 23 The 33 second level possibilities for y-y' path in the case 012-10 for eed case.

MGB-C012-12-eed-16Jul15

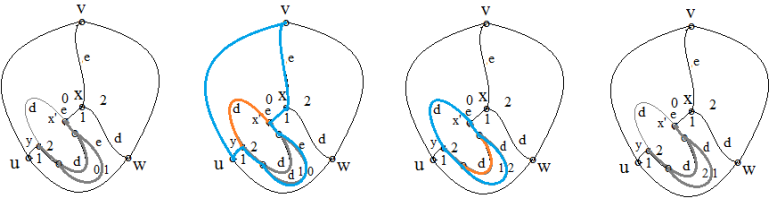

Fig. 24 Case 012-12-eed

Fig.25 corresponds to the continue case 2,2;e. Let y' be a node adjacent to x'. Eight possibilities arise for the y'-z path as shown in Fig.25. It shows cases which stop. In this case all cases stop.

MGB-C012-22-eed-COMP-16Jul15

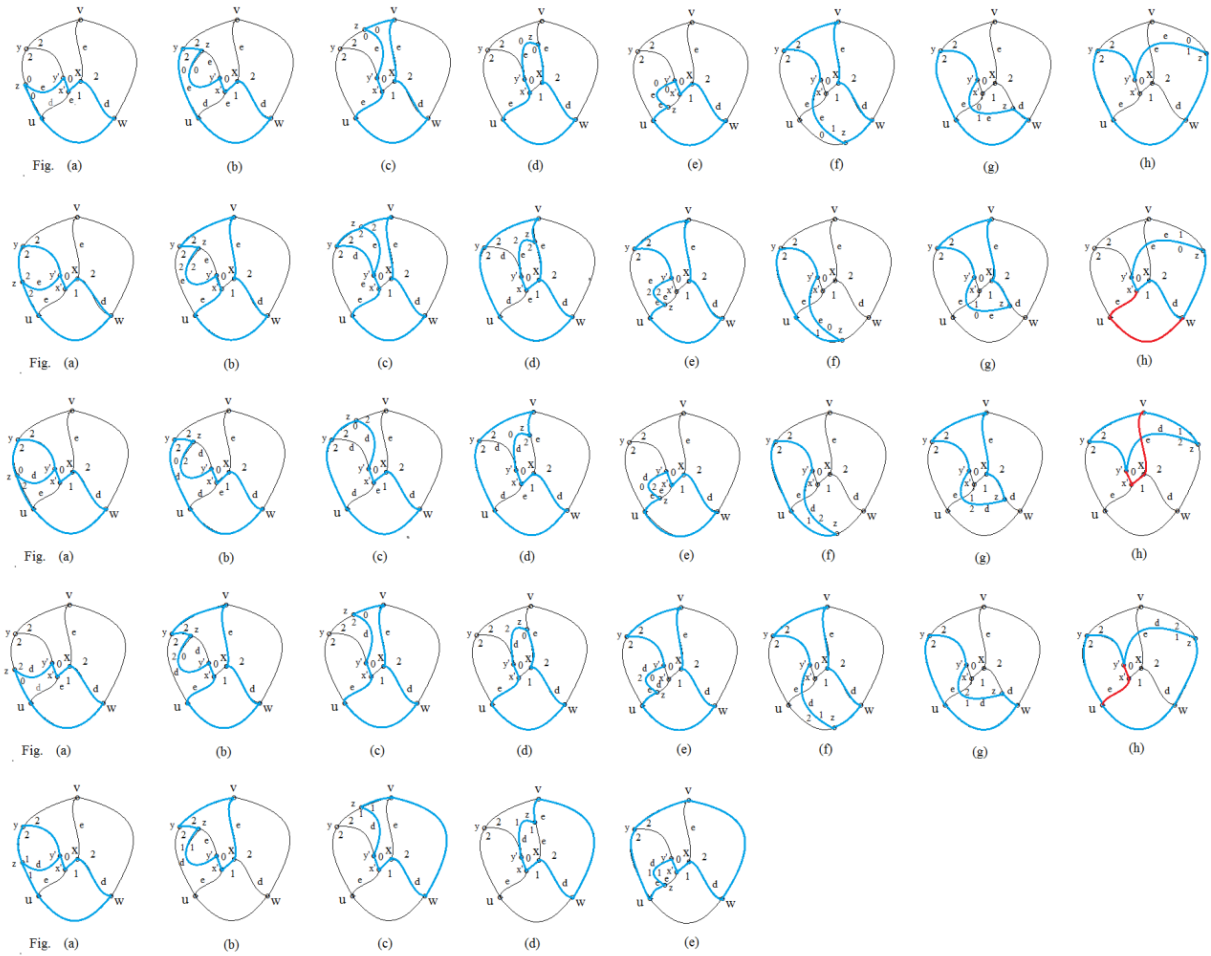

Fig. 25 The 32 second level cases in 0,1,2-22-eed with cycle division. Complete – All subcases stop

Only first level considered for the continue case 0,1,2;eed case.
MGB-C012-eed-16Jul15

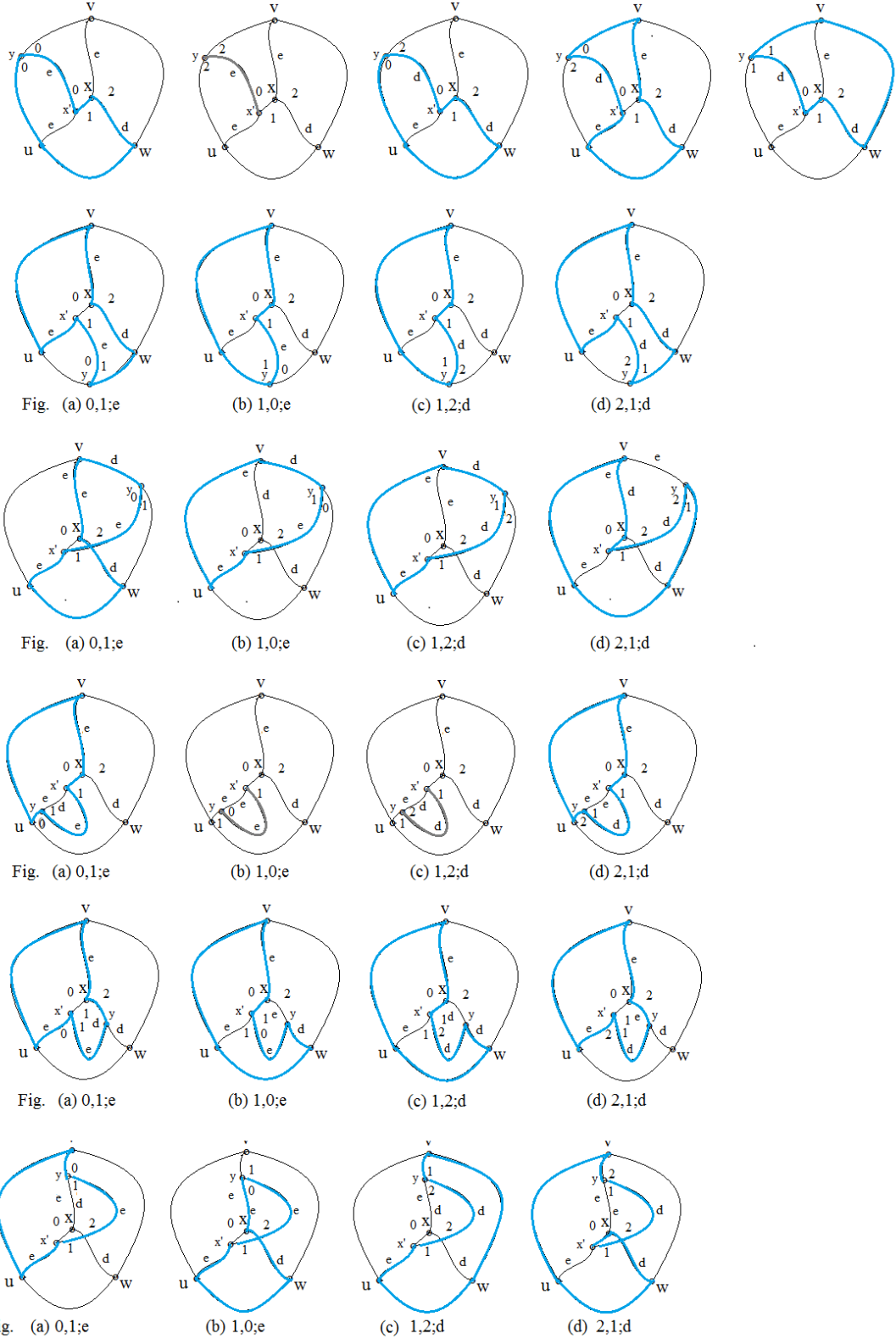

Fig.26 The 25 first level possibilities for x'-y path in the case 0,1,2-eed and stop cases.

Cases to continue to second level:

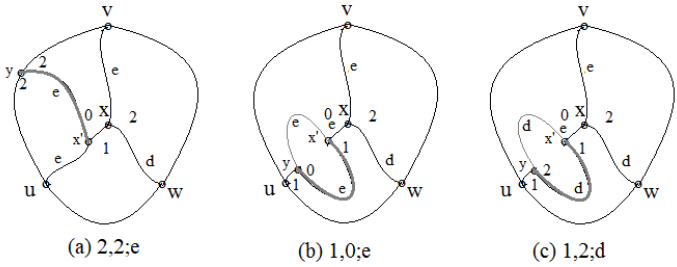

| Table-2 Case 0,1,2 | | | |
|---|---|---|---|
| Cycle Type | | Combined Cycle type | |
| Cycle-1 | Cycle-2 | Int Even | Int Odd |
| 0 | 0 | 0 | 2 |
| 0 | 1 | 1 | 3 |
| 0 | 2 | 2 | 0 |
| 1 | 1 | 2 | 0 |
| 1 | 2 | 3 | 1 |
| 2 | 2 | 0 | 2 |

(a) 2,2;e   (b) 1,0;e   (c) 1,2;d

Fig.27 No solution cases

**Case (g)** Not studied.

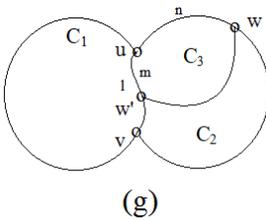

(g)

**Case (h)** Not studied here.

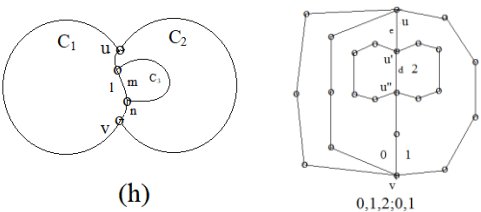

(h)

Fig.28 An example of Euler graph.

**Corollary 5.1.** A regular Euler graph in $\varepsilon_{012}$ is nonexistent in the subcases which stop.

**Case 0,1&3. $\varepsilon_{013}$: Euler graphs with only three types of cycles $C_n$, n≡0,1&3(mod 4)**

**Observation 5.** Order of a graph in $\varepsilon_{013}$ satisfies p≥6. Such a graph $\varepsilon_{013}$ has cycle types (0,1,3). So order is at least 5. Consider a 5-cycle. A 3 & 4-cycles may be obtained by joining two nodes of the 5-cycle at distance 2. This results in a graph with all three types of cycles. But the graph is not Euler. This is achieved by adding a new node adjacent to the same two nodes. The resulting Euler (6,8)-graph is shown in Fig.29a with a graceful numbering and belongs to $\varepsilon_{013}$. A second example is Pythagoras diagram, Fig.29b consisting of three squares on the edges of a triangle, is Euler (9,12)-graph and belongs to $\varepsilon_{013}$. A graceful Euler graph from $\varepsilon_{013}$ is shown in Fig.29c.

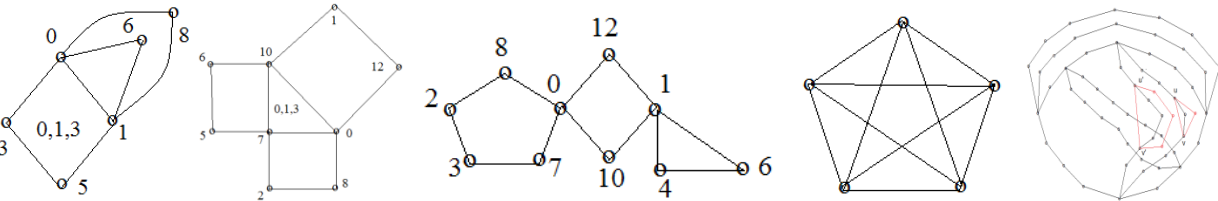

Fig.29 Euler graphs from $\varepsilon_{013}$.

**Observation 6.** Complete graph $K_5$, see Fig.29d, has cycle types 0,1,3. So $K_5$ is an example of a regular Euler graph of degree 4 from $\varepsilon_{013}$. This is the only regular Euler graph known to the author with only three types of cycles 0,1,3. The existence of higher order regular graphs from $\varepsilon_{013}$ remains open. An example of a nonplanar graph from $\varepsilon_{013}$ is given in Fig.29e. The combined cycle (CC) rules in this case are given in the Table-7.

| Table-7    Case 0,1,3 | | | |
|---|---|---|---|
| Cycle Type | | Combined Cycle type | |
| Cycle-1 | Cycle-2 | Int Even | Int Odd |
| 0 | 0 | 0 | 2 |
| 0 | 1 | 1 | 3 |
| 0 | 3 | 3 | 1 |
| 1 | 1 | 2 | 0 |
| 1 | 3 | 0 | 2 |
| 3 | 3 | 2 | 0 |

**Theorem 6.** The CC is of type 2 whenever two cycles of type 0,0 or 1,3 (1,1 or 3,3) of a graph G from $\varepsilon_{013}$ intersect in a path of length l>0 is odd (even).

**Theorem 7.** Size of a graph from $\varepsilon_{013}$ satisfies $q \equiv \xi_1 + 3\xi_3 \pmod{4}$.

Proof follows from Equation (1) (Rao 2014 [10]). Further, if $\xi_1 + 3\xi_3 \equiv 1 \text{ or } 2 \pmod{4}$ then the graphs are nongraceful; else $\xi_1 + 3\xi_3 \equiv 0 \text{ or } 3 \pmod{4}$ and the graphs are candidates for gracefulness. Graphs in $\varepsilon_{013}$ satisfy that $\xi_0 > 0$, $\xi_1 > 0$ and $\xi_3 > 0$.

**Conjecture 4**. Graphs from $\varepsilon_{013}$ satisfying $\xi_1 + 3\xi_3 \equiv 0 \text{ or } 3 \pmod{4}$ are graceful.

**Conjecture 5**. There exists a node of degree two in every graph of order p>5 from $\varepsilon_{013}$.

As a corollary it follows that;

**Conjecture 6.** Regular Euler graphs from $\varepsilon_{013}$ of order p>5 are nonexistent.

or

A regular Euler graph from $\varepsilon_{013}$ of order p>5 with three types of cycles also contains fourth type.

**Theorem 8.** There exists a node of degree two in every graph of order p>5 from $\varepsilon_{013}$ in the cases which stop.

Part Proof in each case:

**Case (a)**  Graphs with three intersecting cycles having a common path

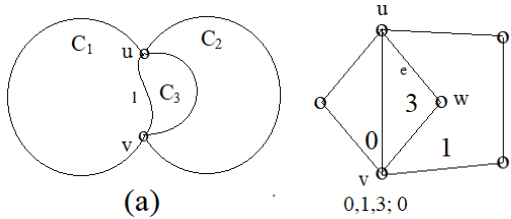

Fig.30 Case (a). Example of Euler graph from $\varepsilon_{013}$ for the case in Table-8.

Considering all possible cycle types for 0,1,3 cycles along with the three CCs results in the only case 0,1,3;0 as shown in the Table-8 below. Example of a graph of this type is shown in Fig.30b.

Table-8 Possibilities for the case (a).

| i | j | k | l | i+j-2l | i+k-2l | j+k-2l |
|---|---|---|---|---|---|---|
| 0 | 1 | 3 | 0 | 1 | 3 | 0 |

Consider Fig.31 with the cycles $C_i$, of type $T_i$, i=0,1,3. The u-w path is even and so there is a node v such that uv is an edge as shown. If there is a degree two node at this level we are done. If not, there is a v-x or v-v' path as shown in Fig.31. This path divides a cycle into two cycles. From the Table-7 the possible cycle combination and path parity can be found out and the cases are enumerated in Fig.31. The cases which stop at this level are shown in orange path indicating a cycle of type 2. Other cases continue. In this case, first case in the first row and fourth case in the second row.

Case-1 First consider v-x path with x on u-x-w path. This path divides $C_0$ into two cycles $C_{01}$, $C_{02}$ with the five possibilities:

0,0;d,

1,3;d, 3,1;d, are impossible as $C_{02}$, $C_1$ and $C_{01}$, $C_1$ intersect in odd path, respectively leading to a contradiction.

1,1,;e, 3,3;e are impossible as $C_{01}$ intersects with $C_3$ in odd path, a contradiction.

Case-2 Here, v-x path with x on $C_1$ divides it into four cases:

0,1;e,

1,0;e,

0,3;d,

3,0;d. v-x path is even. In the first two cases. $C_{11}$ ($C_{12}$) intersects $C_0$ in odd path, a contradiction, see Table-7,

Case-3 Here, v-x path with x on u-x-w path of $C_3$ and divides $C_3$ into two cycles with four possibilities:

0,1;d,

1,0;d,

0,3;e,

3,0;e. $C_{31}$ ($C_{32}$) intersect $C_0$ in odd paths, a contradiction.

Case-4 Here, v-x path with x on u-x-w path of $C_3$ and divides $C_3$ into two cycles with four possibilities:

0,1;d,

1,0;d,

0,3;e,

3,0;e. $C_{31}$ ($C_{32}$) intersect $C_0$ in odd paths, a contradiction.

MGB-CijkCP-013CPEven-03Aug15

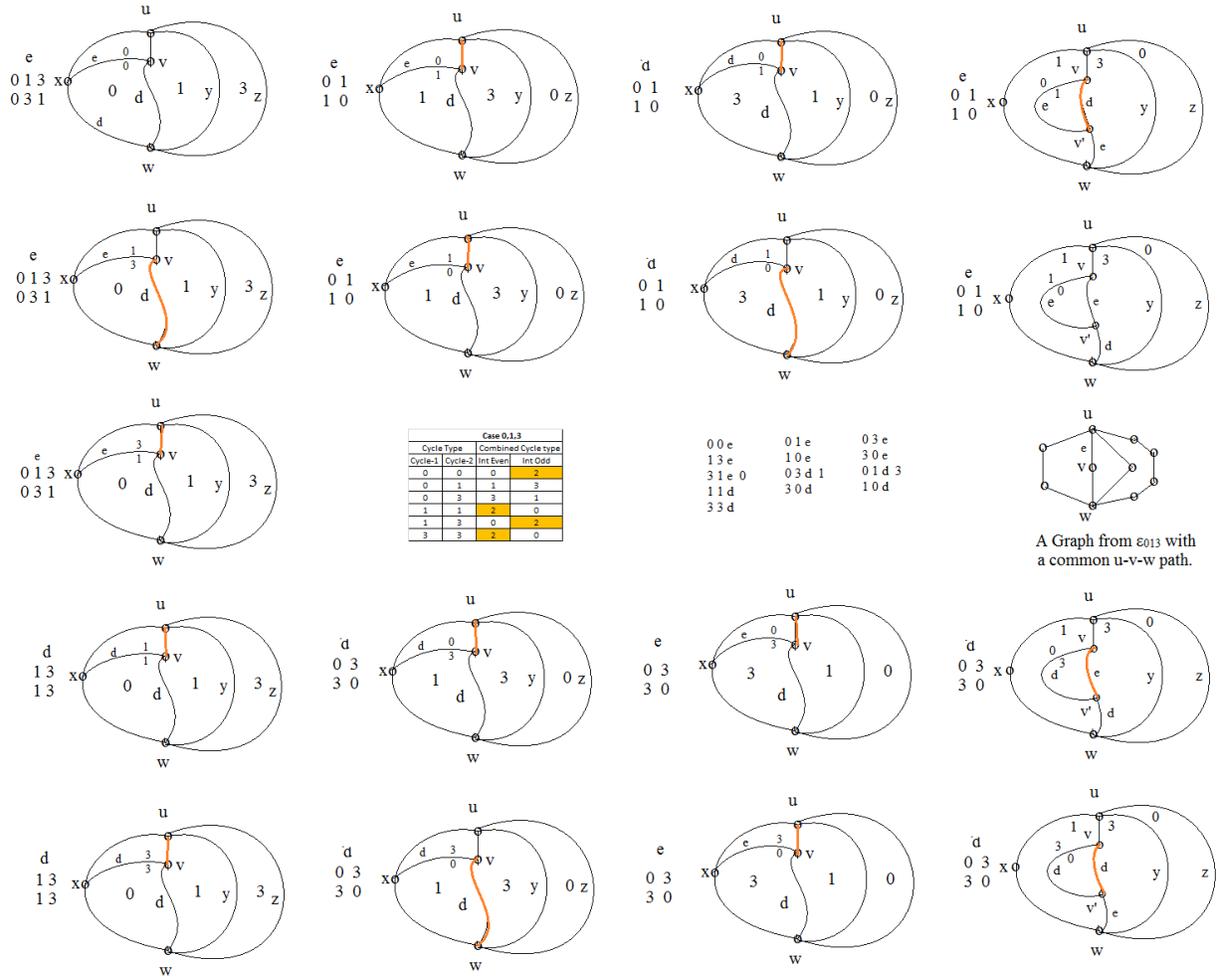

Fig.31 The 17 possibilities at the first level in the case 0,1,3;e with cycle divisions and stop cases.

## Case b

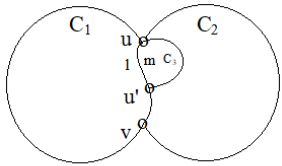

(b)

Table-9 Possible cases after eliminating ijk and jik same cases.

| i | j | k | l | m | i+j-2l | i+k-2m | j+k-2m |
|---|---|---|---|---|--------|--------|--------|
| 0 | 1 | 3 | 0 | 0 | 1 | 3 | 0 |
| 0 | 1 | 3 | 1 | 0 | 3 | 3 | 0 |
| 0 | 3 | 1 | 0 | 0 | 3 | 1 | 0 |
| 0 | 3 | 1 | 1 | 0 | 1 | 1 | 0 |
| 1 | 3 | 0 | 0 | 0 | 0 | 1 | 3 |
| 1 | 3 | 0 | 0 | 1 | 0 | 3 | 1 |

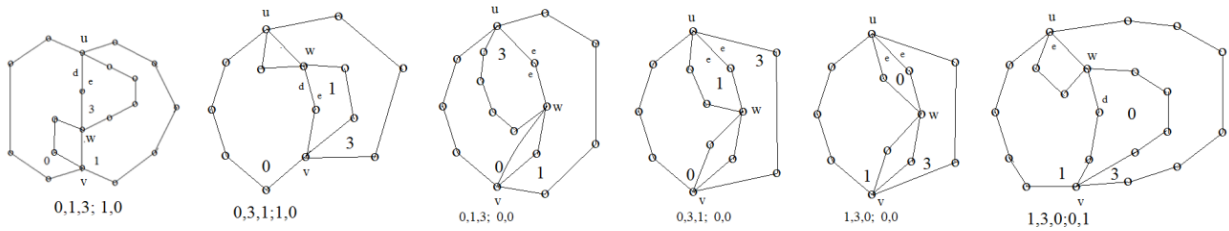

0,1,3; 1,0          0,3,1;1,0          0,1,3; 0,0          0,3,1; 0,0          1,3,0; 0,0          1,3,0;0,1

Fig.32 Examples of Euler graphs from $\varepsilon_{013}$ for the cases in the Table-9.

MGB-C3-case2-01300-19Oct16

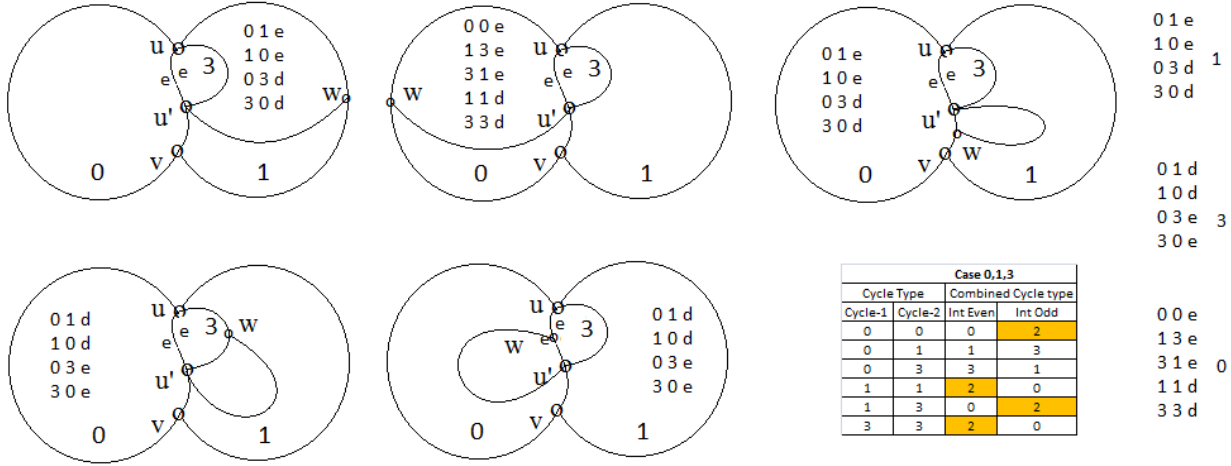

Fig.33 The five possibilities with u'-w path and cycle divisions for 0,1,3;e,e case.

MGB-C3-case2-03100-19Oct16

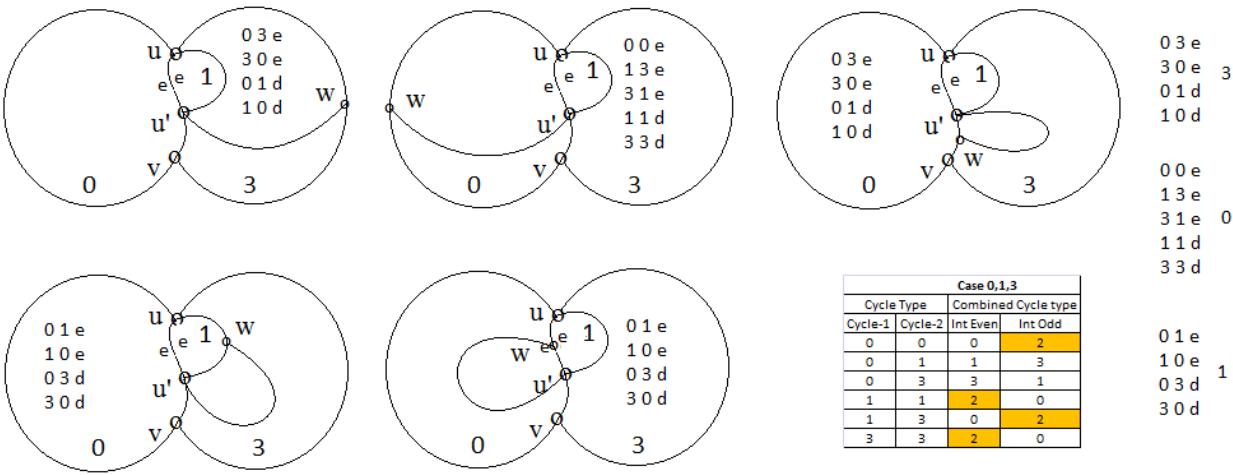

Fig.34 The five possibilities with u'-w path and cycle divisions for 0,3,1;e,e case.

The Fig.35 shows possibilities for the u'-x path in the case 0,1,3;e,e. Also shown are the cycle divisions containing the path using the rules from Table-7and stop cases.

MGB-C013CI-00-11Dec15.jpg

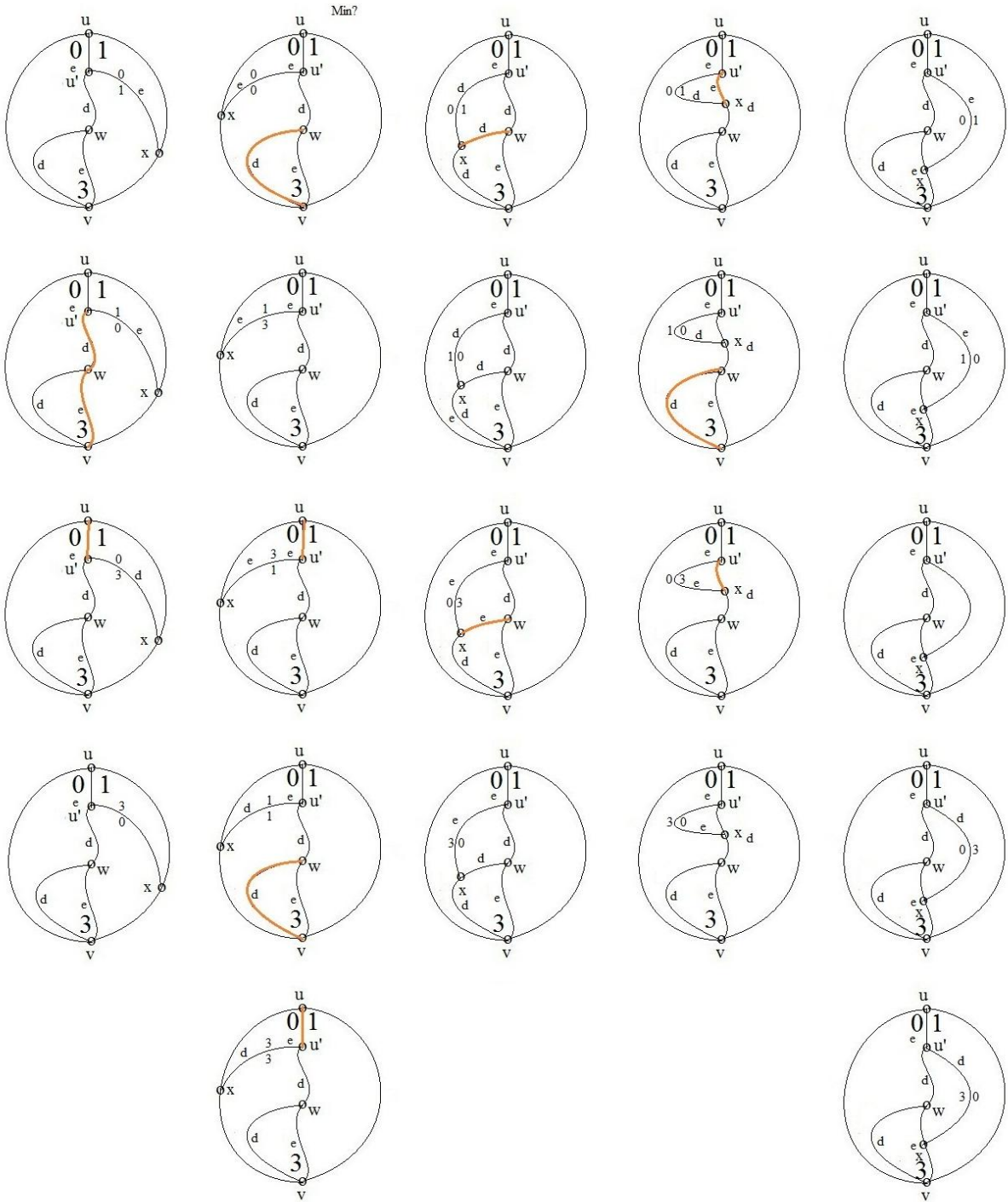

Fig.35 The 22 first level possibilities with u'-x path with parity and cycle division in the case 0,1,3;e,e.

The following Figs.36-40 show the possibilities for the remaining 5 cases with the cycle combination and path parities.

MGB-C013CI-10-11Dec15.jpg

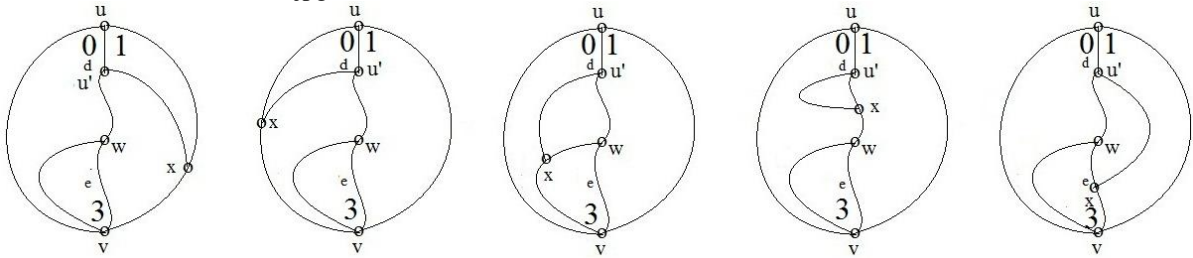

Fig.36 The u'-x path possibilities in 0,1,3;d,e case.

MGB-C031CI-00-11Dec15.jpg

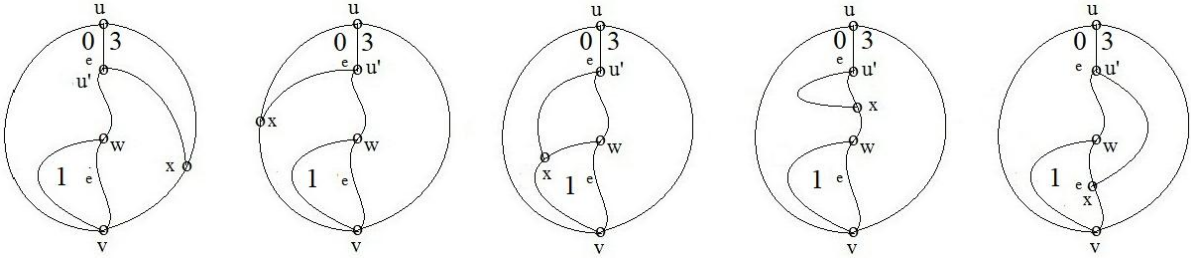

Fig.37 The u'-x path possibilities in 0,3,1;e,e case.

MGB-C031CI-10-11Dec15.jpg

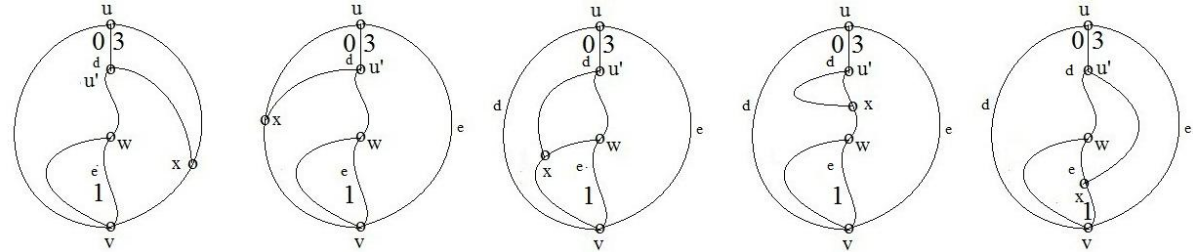

Fig.38 The u'-x path possibilities in 0,3,1;d,e case.

MGB-C130CI-00-11Dec15.jpg

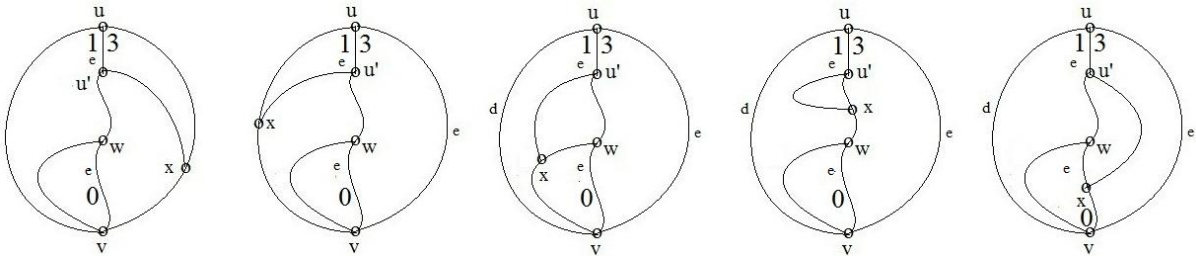

Fig.39 The u'-x path possibilities in 1,3,0;e,e case.

MGB-C130CI-10-11Dec15.jpg

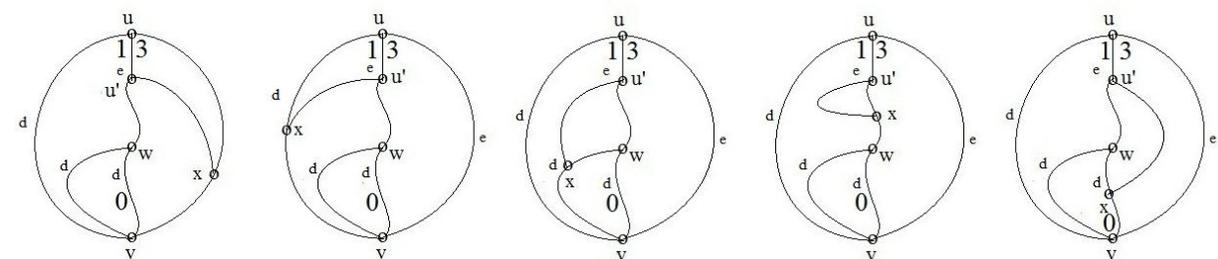

Fig.40 The u'-x path possibilities in 1,3,0;d,e case.

**Case c**

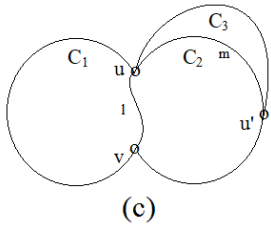

(c)

Table-10 Possible combinations after eliminating i,j,k and k,j,i same cases.

| i | j | k | l | m | i+j-2l | j+k-2m | i+j+k-2l-2m |
|---|---|---|---|---|--------|--------|-------------|
| 0 | 1 | 3 | 0 | 0 | 1 | 0 | 0 |
| 0 | 3 | 1 | 0 | 0 | 3 | 0 | 0 |
| 1 | 0 | 3 | 0 | 0 | 1 | 3 | 0 |
| 1 | 0 | 3 | 1 | 1 | 3 | 1 | 0 |

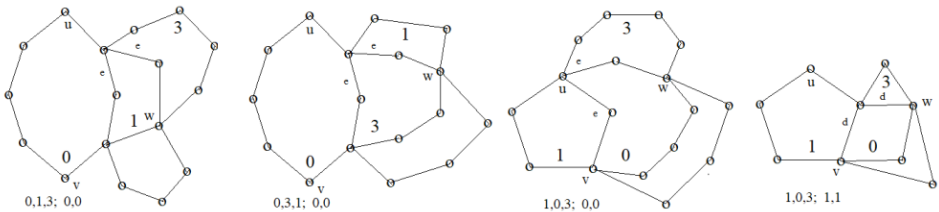

Fig.41 Examples of Euler graphs from $\varepsilon_{013}$ for each case in the Table-10.

MGB-C3-case1-01300-19Oct16

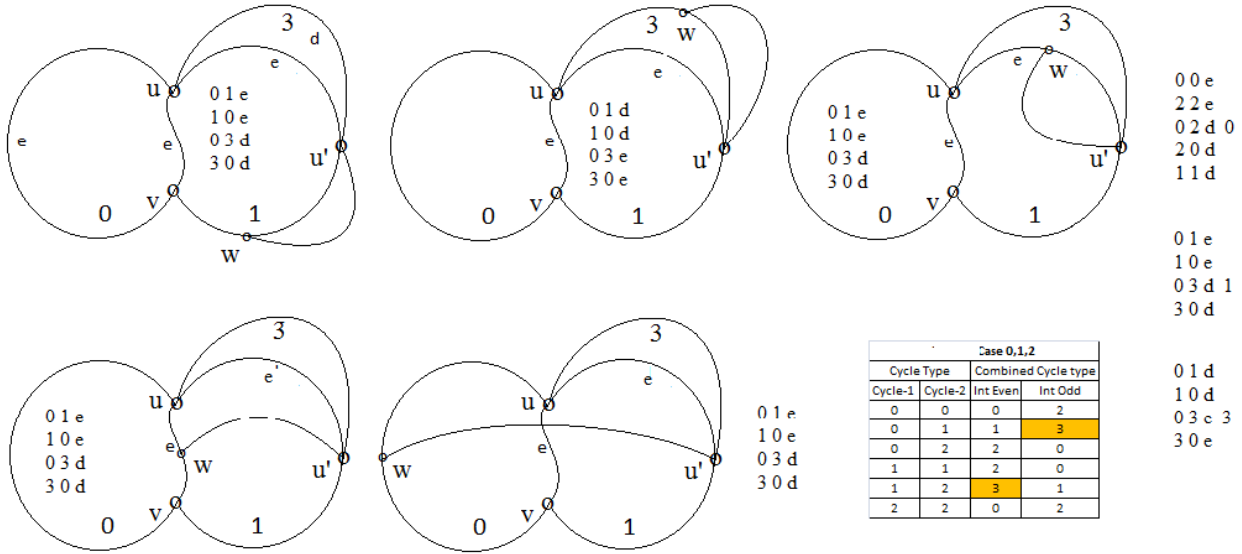

Fig.42 The five possibilities for u'-w path with division of cycles in 0,1,3;e,e case.

MGB-C3-case1-03100-19Oct16

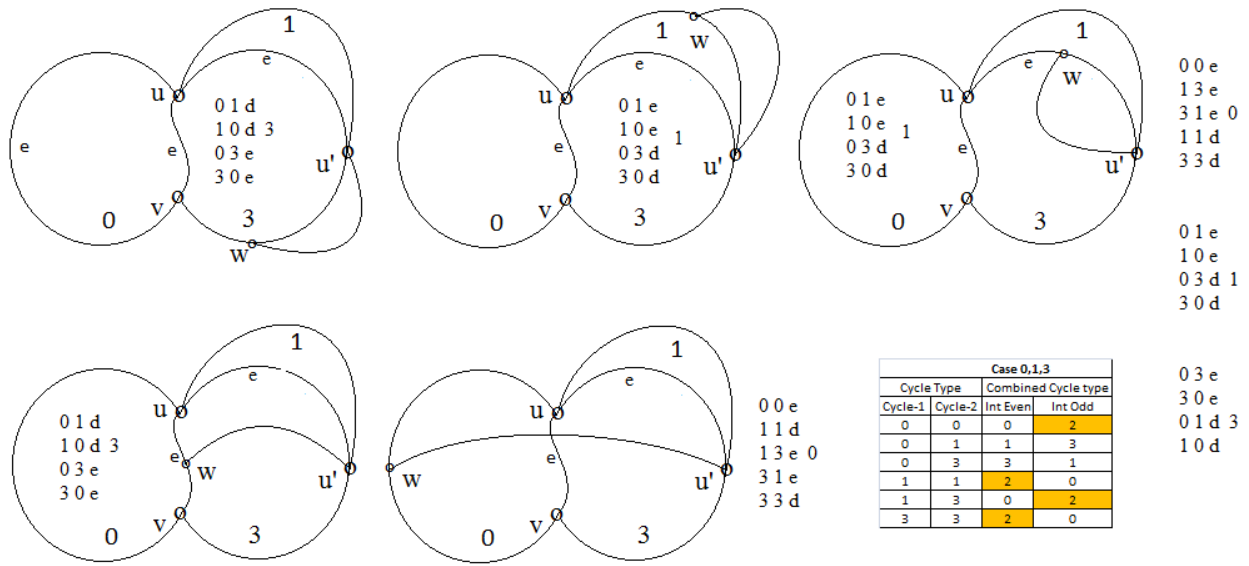

Fig.43 The five possibilities for u'-w path with division of cycles in 0,3,1;e,e case.

**Case (d)** Same as Case c

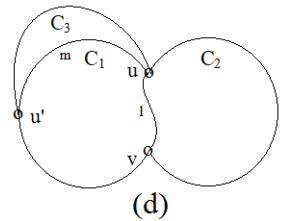

(d)

**Case (e)**

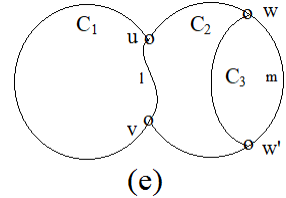

(e)

Table-11 Possible combinations after eliminating i,j,k and k,j,i same cases.

| I | j | k | l | m | i+j-2l | j+k-2m | i+j+k-2l-2m |
|---|---|---|---|---|--------|--------|-------------|
| 0 | 1 | 3 | 0 | 0 | 1 | 0 | 0 |
| 0 | 3 | 1 | 0 | 0 | 3 | 0 | 0 |
| 1 | 0 | 3 | 0 | 0 | 1 | 3 | 0 |
| 1 | 0 | 3 | 1 | 1 | 3 | 1 | 0 |

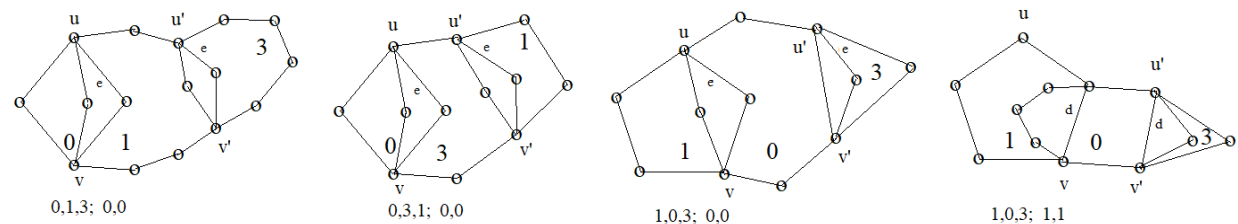

0,1,3;  0,0            0,3,1;  0,0            1,0,3;  0,0            1,0,3;  1,1

Fig.44a Examples of Euler graphs from $\varepsilon_{013}$ for each of the case in the Table-11.

MGB-Cijk-3CC-013-ee-14Jun20

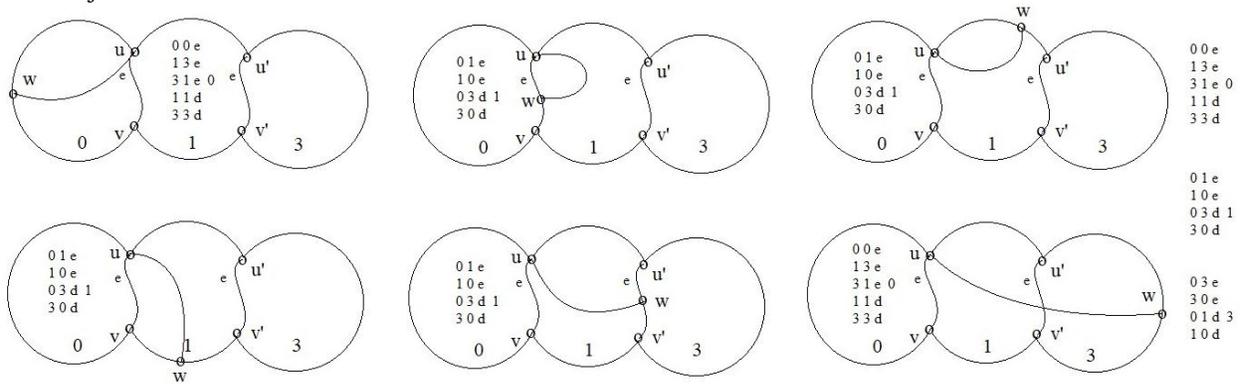

Fig.44b  The six first level possibilities for u-w path with cycle division in the case 0,1,3;e,e.

**Case (f)**     Graphs with three intersecting cycles, any two intersecting in different paths

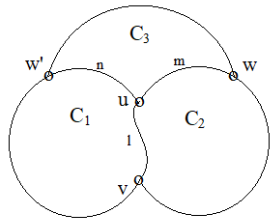

(f)

Table-12 Possible combinations after eliminating i,j,k; j,k,i and k,i,j same cases.

| I | j | k | l | m | n | i+j-2l | j+k-2m | i+k-2n | i+j+k-2l-2(m+n) |
|---|---|---|---|---|---|--------|--------|--------|------------------|
| 0 | 1 | 3 | 0 | 0 | 0 | 1 | 0 | 3 | 0 |
| 0 | 1 | 3 | 1 | 0 | 1 | 3 | 0 | 1 | 0 |
| 0 | 3 | 1 | 1 | 0 | 1 | 1 | 0 | 3 | 0 |

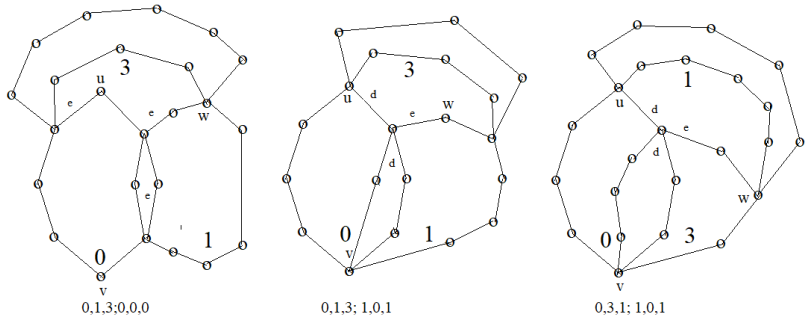

0,1,3;0,0,0          0,1,3; 1,0,1          0,3,1; 1,0,1

Fig.45 Examples of Euler graphs from $\varepsilon_{013}$ for the cases in the Table-12.

MGB-C013-dde-01d-30Jan2016

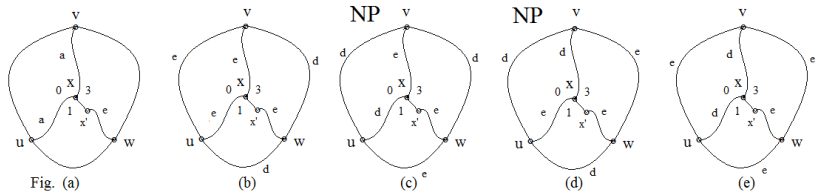

Fig.46 Possible parity combinations from general case (a) aae to eee, ede, dee,dde for the intersection paths.

Fig 46c is impossible as CC of $C_1$ and $C_2$ is of type $T_1$ and intersects with $C_3$ in odd path, a contradiction as the length of intersection must be even from Table-7.

Fig 46d is impossible as CC of $C_1$ and $C_2$ is of type $T_3$ and intersects with $C_3$ in even path, a contradiction as the length of intersection must be odd from Table-7.

**Case-I**: Case (b) of Fig.46. All the three paths are even; eee.

MGB-C013-eee-16Jul15

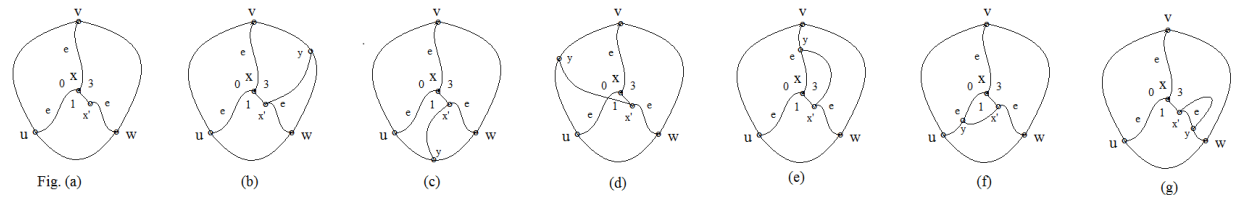

Fig.47 The six first level possibilities for x'-y path in the case eee.

Here, u-x and v-x paths may be of any parity and w-x path is even. Four cases arise as in Fig.47b to d. Since w-x path is even there is x' adjacent to x on w-x path. If x' is of degree two then we are done. Else, there is y such that x'-y path exists. Consider the first case as in Fig.47b. Six possible cases arise as shown in Figs.47b to g. In all these cases we shall show the existence of a cycle type $T_2$ or existence of two cycles with intersection not as in Table-7.

Case b the x'-y path divides $C_3$ into two cycles. From Table-7 it follows that $C_3$ gets divided into $C_{31}$:v-x,x'-y-u; $C_{32}$:w-x'-y-w with possibilities:

  0,1,d;

  1,0,d;

  0,3,e; is impossible as intersection of $C_2$ and $C_{31}$ is odd, a contradiction. It must be even from Table-7.

  3,0,e is impossible as intersection of $C_2$ and $C_{32}$ is odd, a contradiction. It must be even from Table-7.

Case c,d the x'-y path divides $C_1$ or CC of $C_0$ and $C_1$ of type $T_1$ into two cycles. From Table-7 it follows that $C_1$ gets divided into:

  0,1,e; 1,0,e; are impossible as intersection of $C_2$ and 3-cycle is odd.

  0,3,d; and 3,0,d.

Case e,f

Case g may be considered as x'-y path dividing $C_3$ into $C_{31}$:v-x,x'-y-w-v and $C_{32}$:x'-y-x':

  0,1,d; is  impossible as x-y'-y path is odd implying that  x'-y path is even, but intersection of $C_2$ and $C_{32}$ must be odd from Table-7.

  1,0,d; x'-y'-y path is odd and so x'-y is even, however, $C_{32}$ and $C_2$ must intersect in odd path, a contradiction.

  0,3,e; is  impossible as it follows that x'-y path is odd but intersection of $C_2$ and $C_{32}$ must be even. Lastly,

  3,0,e; is  impossible as it follows that x'-y'-y path is even and so x'-y path is odd implying that y-w path is even. Therefore, x,x'-y'-y-w path is odd, a contradiction.

MGB-C013-eee-16Jul15

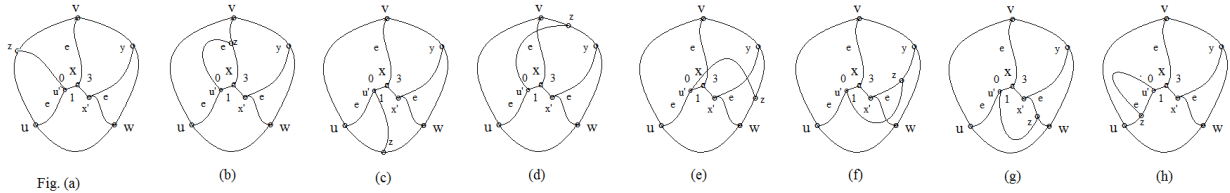

Fig.48 The eight second level possibilities for u'-z path in the case shown even x'-y path.

Case a,b the u'-z path divides $C_1$ into two cycles. From Table-7, it follows that $C_1$ gets divided into $C_{11}$:u-u'-z-u; $C_{12}$:x,u'-y-v-x with possibilities 0,0,e; is impossible as CC of $C_{12}$ and $C_2$ with odd intersection results in $T_3$ cycle with odd intersection with second $T_1$ cycle. 1,1,d; 1,3,e; 3,1,e; 3,3,d:

Case c,d,f the x'-y path divides $C_1$ or CC of $C_0$ and $C_1$ of type $T_1$ into two cycles. From Table-7 it follows that $C_1$ gets divided into:
0,1,e; 1,0,e; are impossible as intersection of $C_2$ and 3-cycle is odd.
0,3,d; and 3,0,d.

Case g may be considered as x'-y path dividing $C_3$ into:
0,1,d; is impossible as intersection of $C_2$ and 1-cycle is even,
1,0,d;
0,3,e; is impossible as intersection of $C_2$ and 3-cycle is odd,
3,0,e.

**Case-II:** Case (e) of Fig.46. The paths are odd, odd, even; dde

MGB-C013-dde-23Jan2016

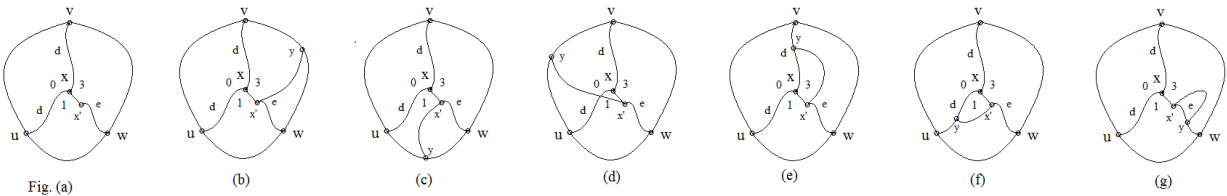

Fig.49 The six first level possibilities for x'-y path in the case dde.

Case b leads to four subcases for the division of $C_3$ into:
0,1;d, impossible as $C_1$ and $C_{31}$ intersect in odd path, a contradiction to the intersection must be even from Table-7.
1,0;d,
0,3;e, impossible as $C_2$ and $C_{32}$ intersect in odd path, a contradiction whereas the intersection must be even from Table-7.
3,0;e. impossible as $C_2$ and $C_{31}$ intersect in odd path, a contradiction whereas the intersection must be even from Table-7.

Case c leads to four subcases for the division of $C_1$ into:
0,1;e, impossible as $C_{21}$ and $C_3$ intersect in odd path, a contradiction whereas the intersection must be even from Table-7.
1,0;e, impossible as $C_{22}$ and $C_3$ intersect in odd path, a contradiction whereas the intersection must be even from Table-7.
0,3;d, impossible as $C_1$ and $C_{21}$ intersect in odd path, a contradiction whereas the intersection must be even from Table-7.
3,0;d.

MGB-C013-dde-23Jan2016

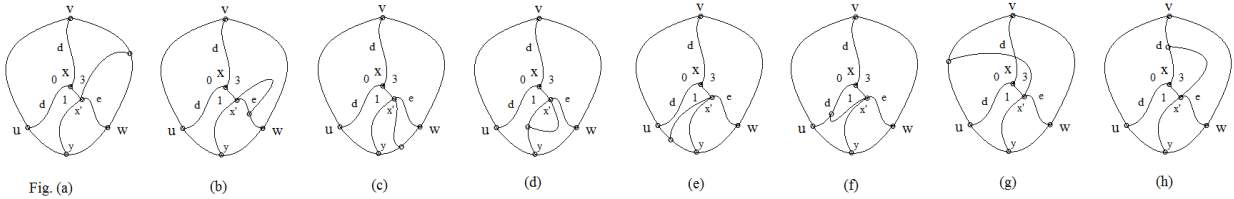

Fig.50 The eight second level possibilities for x'-y' path.

MGB-C013-dde-23Jan2016

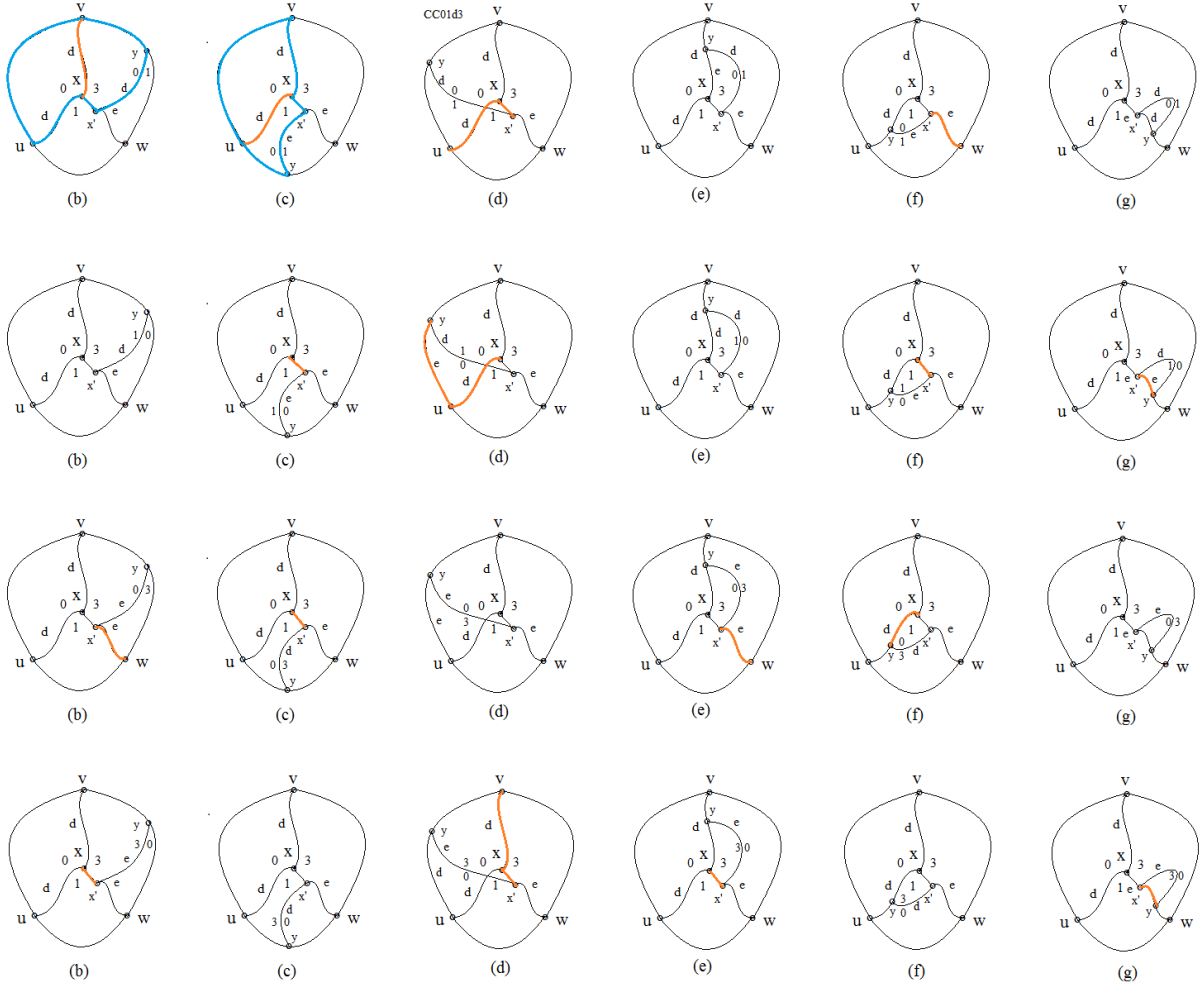

Fig.51 The 24 first level possibilities with cycle divisions for x'-y path with stop cases in the case 013-dde-01d.

MGB-C013-dde-01d-30Jan2016

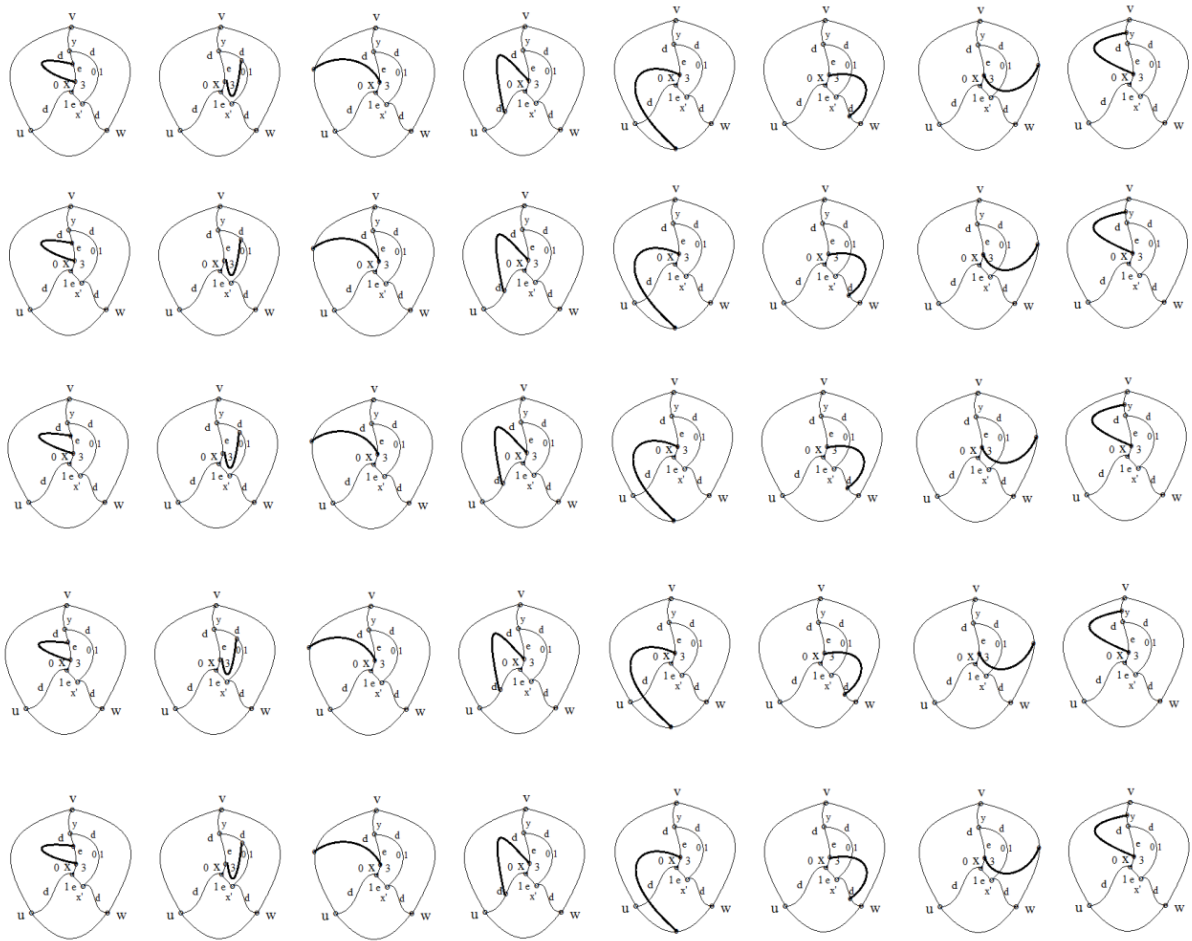

Fig.52 The 40 second level possibilities for x'-y path odd dividing cycle 3 in the case 013-dde-01d.

MGB-C013-dde-23Jan2016

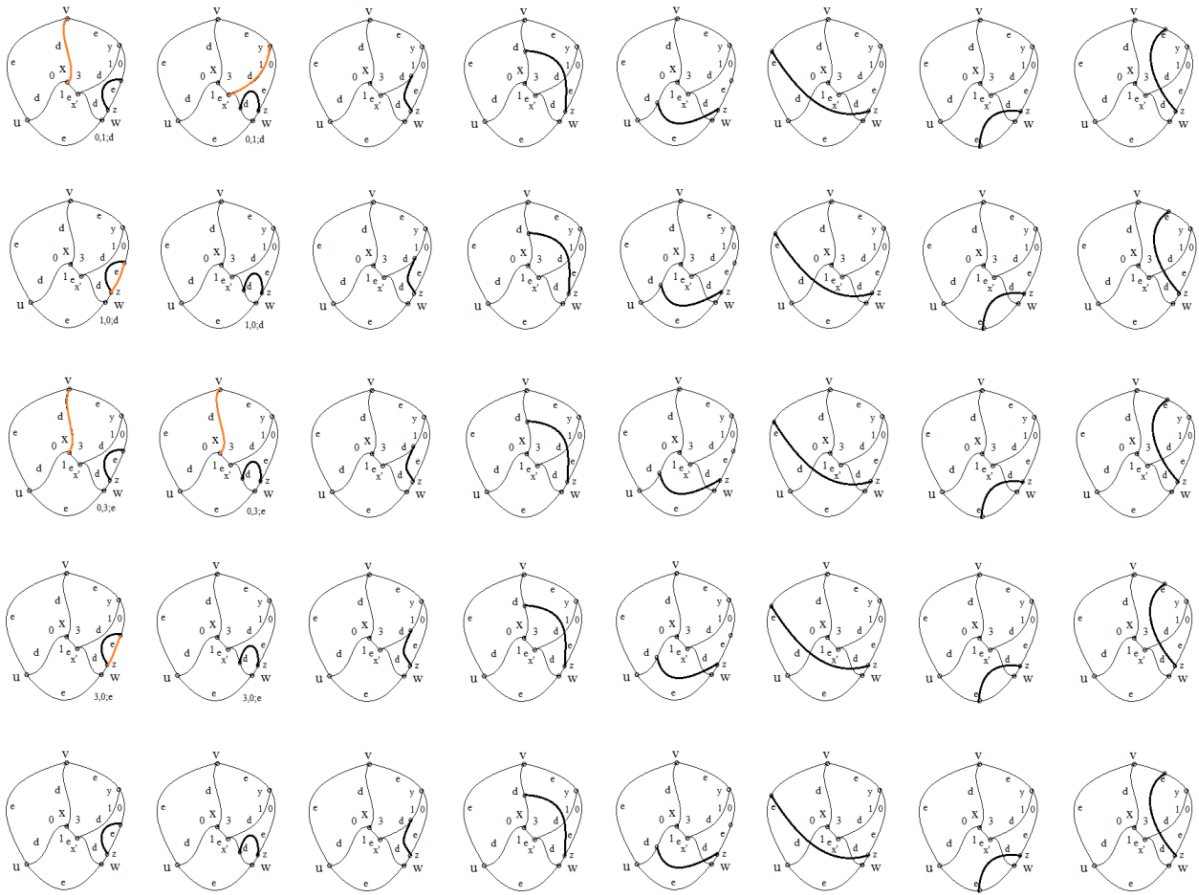

Fig.53 The 40 second level possibilities for the z-z' path, x'-y odd path dividing cycle 3 into 1,0 in the case 013-dde.

MGB-C013-eee-16Jul15

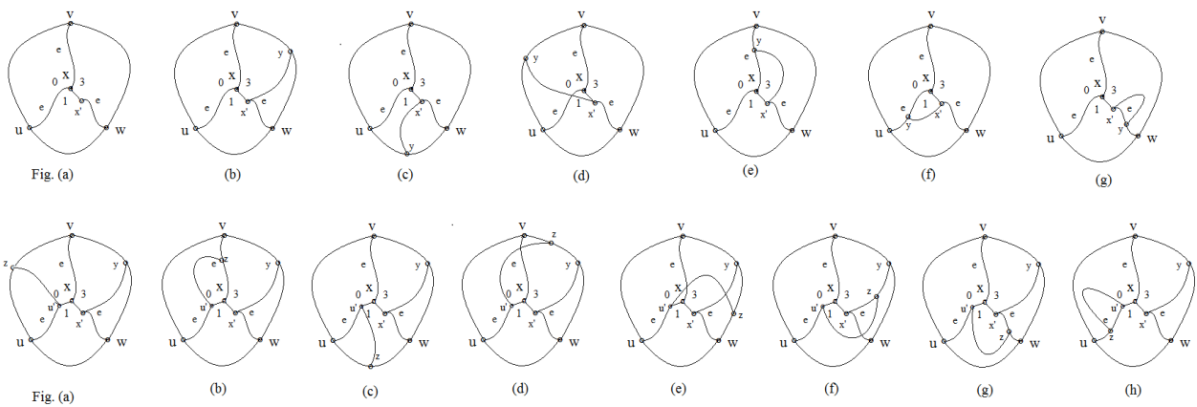

Fig.54 The seven first level possibilities for x'-y path and the eight second level possibilities for the case (b) in the first row in the case 013-eee.

MGB-C013-eee-16Jul15

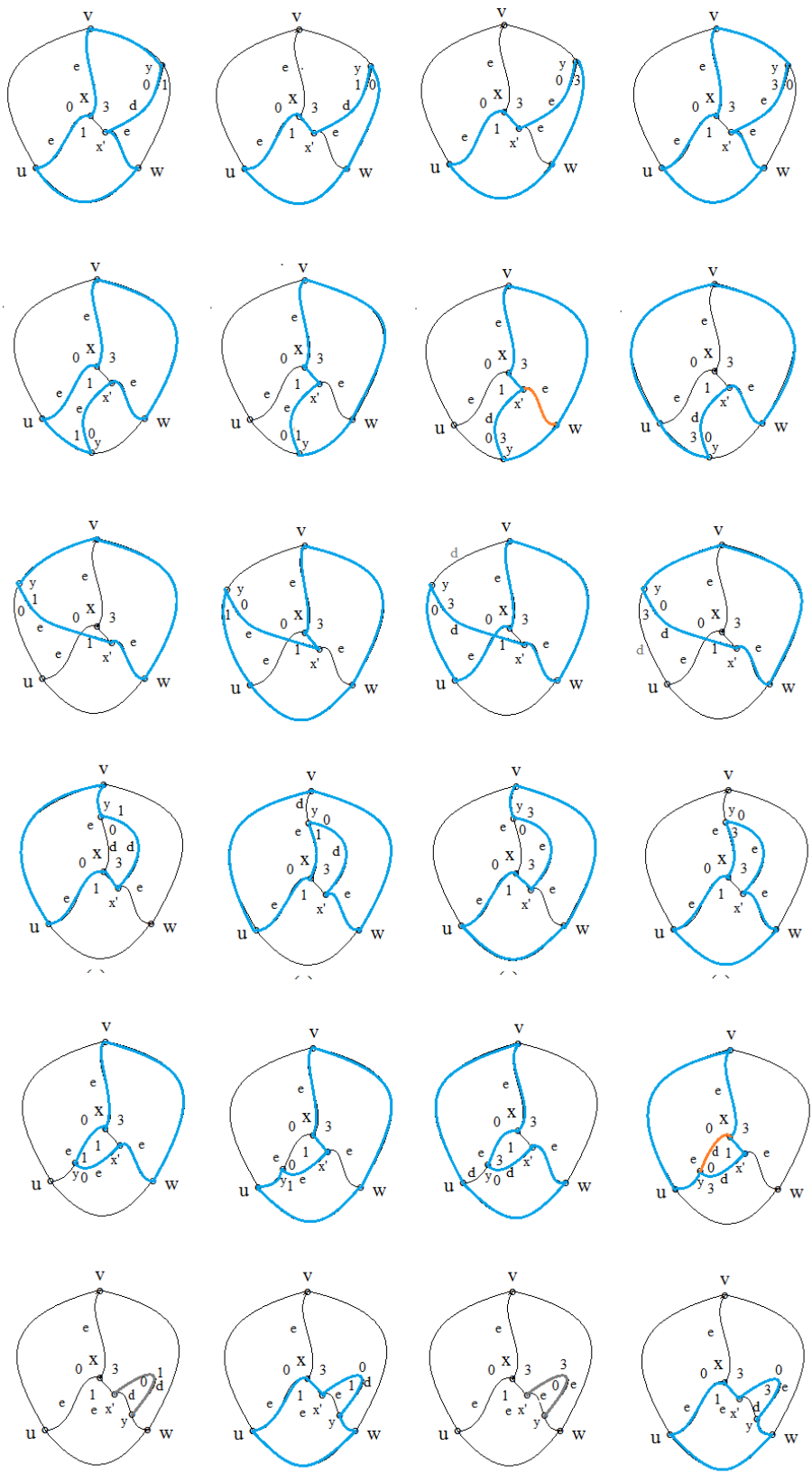

Fig.55 The 24 first level possibilities for the x'-y path with stop cases in case 013-eee.

Cases to continue to second level.

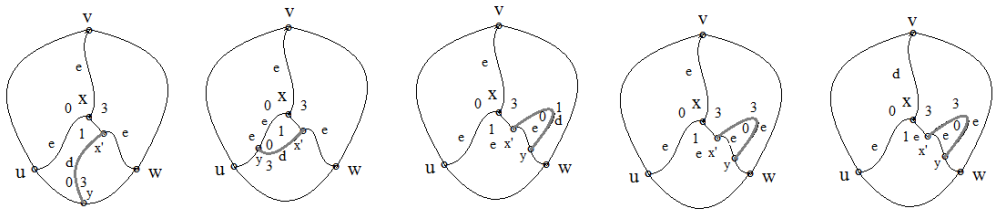

Fig.56 No solution cases

**Case (g)** Not studied.

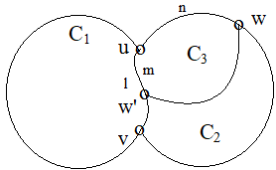

**Case (h)** Not studied.

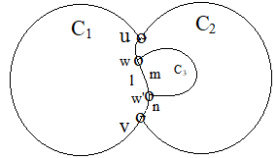

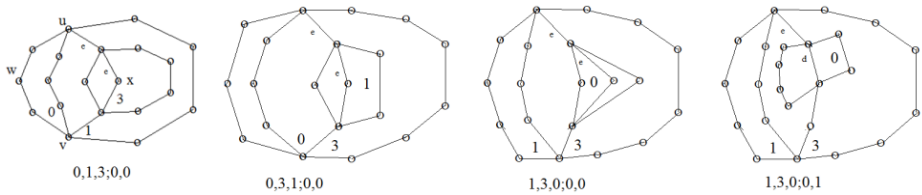

Fig.57 Examples of Euler graphs from $\varepsilon_{013}$.

**Corollary 8.1.** A regular Euler graph in $\varepsilon_{013}$ is nonexistent in the subcases which stop.

**Case 0,2&3. $\varepsilon_{023}$: Euler graphs with only three types of cycles $C_n$, n≡0,2&3(mod 4)**

**Observation 7.** Order of a graph from $\varepsilon_{023}$ satisfies p≥8. Such a graph has cycle types (0,2,3), so order is at least 6 and any 6-cycle has no odd chords. If an even chord exists then a 3-cycle is not possible with two nodes of 6-cycle in common else a 5-cycle is formed along with an even chord. Fig.58a is Euler (8,10)-graph has minimum order with a 3-cycle, a 4-cycle and a 6-cycle. 3-cycle leads to a 7-cycle which is type 3. Fig.58b,c are examples of Euler graphs from $\varepsilon_{023}$ of order 10 with size 13 and14. All three graphs are nongraceful by Rosa-Golomb criterion.

Graphs in $\varepsilon_{023}$ may be nonplanar. An example of nonplanar graph from $\varepsilon_{023}$ is shown in Fig.58d.

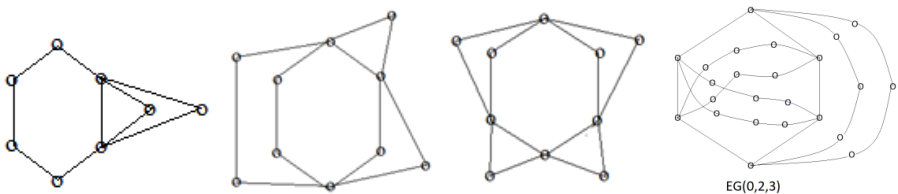

Fig.58 Graphs from $\varepsilon_{023}$.

| Table-14 | Case 0,2,3 | | |
|---|---|---|---|
| Cycle Type | | Combined Cycle type | |
| Cycle-1 | Cycle-2 | Int Even | Int Odd |
| 0 | 0 | 0 | 2 |
| 0 | 2 | 2 | 0 |
| 0 | 3 | 3 | 1 |
| 2 | 2 | 0 | 2 |
| 2 | 3 | 1 | 3 |
| 3 | 3 | 2 | 0 |

The CC rules in this case is given in the Table-14.

**Theorem 9.** The CC is of type 1 whenever two cycles of type 0 and 3 (2 and 3) of a graph G from $\varepsilon_{023}$ intersect in a path of length l>0 is odd (even).

**Theorem 10**. A graph in $\varepsilon_{023}$ with three cycles as shown in cases (c) and (e) where an ordered triple (i,j,k)=(0,2,3), (0,3,2) or (2,0,3) has a combined cycle $C_n$, n≡1(mod 4) whenever l+m is even.

**Theorem 11.** Size of a graph from $\varepsilon_{023}$ satisfies q≡$2\xi_2+3\xi_3$ (mod 4).

Proof follows from Equation (1) (Rao 2014 [10]). Further, if $2\xi_2+3\xi_3$≡1or2(mod 4) then the graphs are nongraceful; else $2\xi_2+3\xi_3$≡0or3(mod 4) and the graphs are candidates for gracefulness. Graphs in $\varepsilon_{023}$ satisfy that $\xi_0$>0, $\xi_2$>0 and $\xi_3$>0.

**Theorem 12.** Two necessary conditions follow: If $2\xi_2+3\xi_3$≡0(mod 4) then $\xi_3$ is even. If $2\xi_2+3\xi_3$≡3(mod 4) then $\xi_3$-1 is even implying $\xi_3$ is odd.

**Conjecture 7**. Graphs in $\varepsilon_{023}$ satisfying $2\xi_2+3\xi_3$≡0or3(mod 4) are graceful.

**Conjecture 8**. There exists a node of degree two in every graph of order p>7 from $\varepsilon_{023}$.

As a corollary it follows that:

**Conjecture 9**. Regular Euler graphs from $\varepsilon_{023}$ of order >7 are nonexistent.

or

A regular Euler graph from $\varepsilon_{023}$ of order p>7 with three types of cycles also contains fourth type.

**Theorem 13**. There exists a node of degree two in every graph of order p>6 from $\varepsilon_{023}$ in the cases which stop.

**Case (a)**   Not possible.

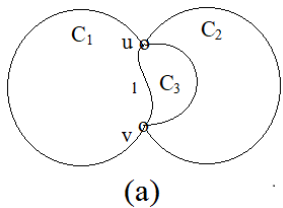

(a)

**Case (b)**

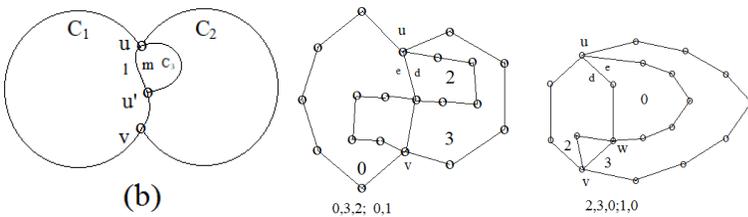

(b)

Fig.59 Case (b). Examples of Euler graphs from $\varepsilon_{023}$ for the cases in the Table-15.

Table-15 Possible cases in the Case (b).

| i | j | k | l | m | i+j-2l | i+k-2m | j+k-2m |
|---|---|---|---|---|--------|--------|--------|
| 0 | 3 | 2 | 0 | 1 | 3 | 0 | 3 |
| 2 | 3 | 0 | 1 | 0 | 3 | 2 | 3 |

The analysis of the cases is done as in the case 0,1,2.

MGB-C3-case2-02301-19Oct16

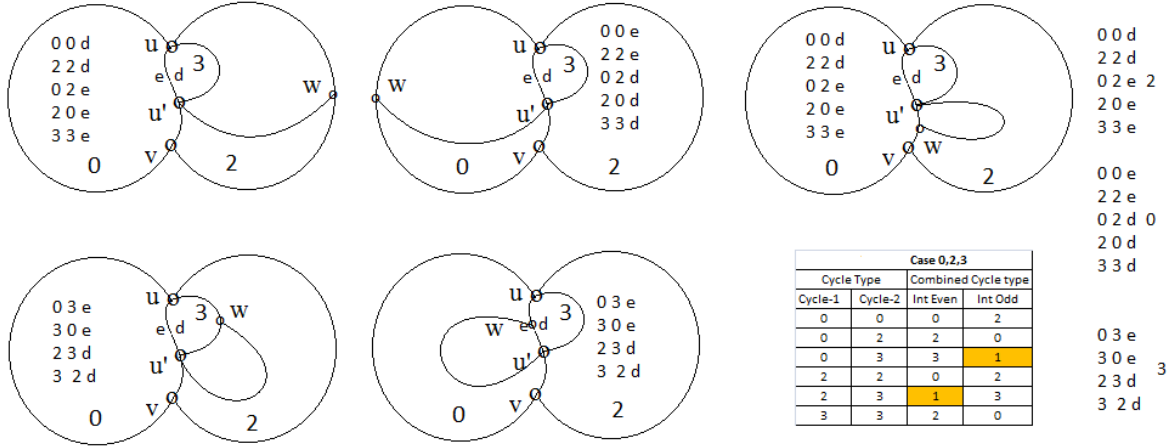

Fig.60 The five possibilities for u'-w path with division of cycles for the case 0,2,3;e,d.

MGB-C3-case2-03201-19Oct16

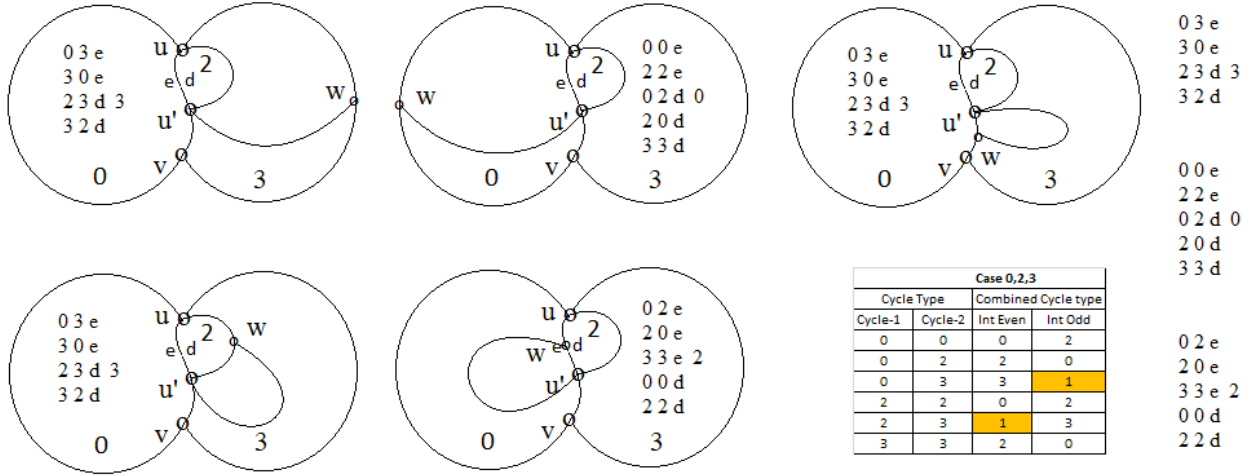

Fig.61 The five possibilities for u'-w path with division of cycles for the case 0,3,2;e,d.

**Case (c)**

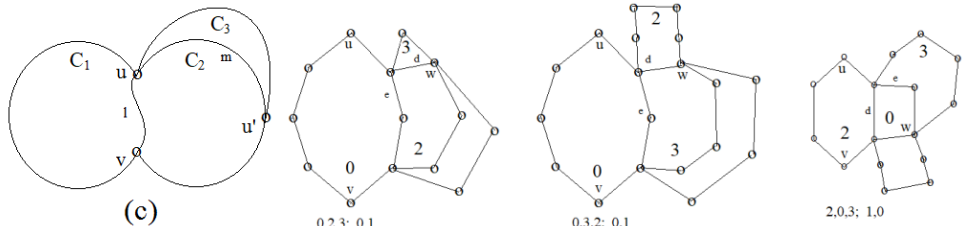

Fig.62 Case (c). Examples of Euler graphs from $\varepsilon_{023}$ for the cases in the Table-16.

Table-16 Possible cases in the case (c).

| i | j | k | l | m | i+j-2l | j+k-2m | i+j+k-2l-2m |
|---|---|---|---|---|--------|--------|-------------|
| 0 | 2 | 3 | 0 | 1 | 2 | 3 | 3 |
| 0 | 3 | 2 | 0 | 1 | 3 | 3 | 3 |
| 2 | 0 | 3 | 1 | 0 | 0 | 3 | 3 |

MGB-C3-case1-02301-19Oct16

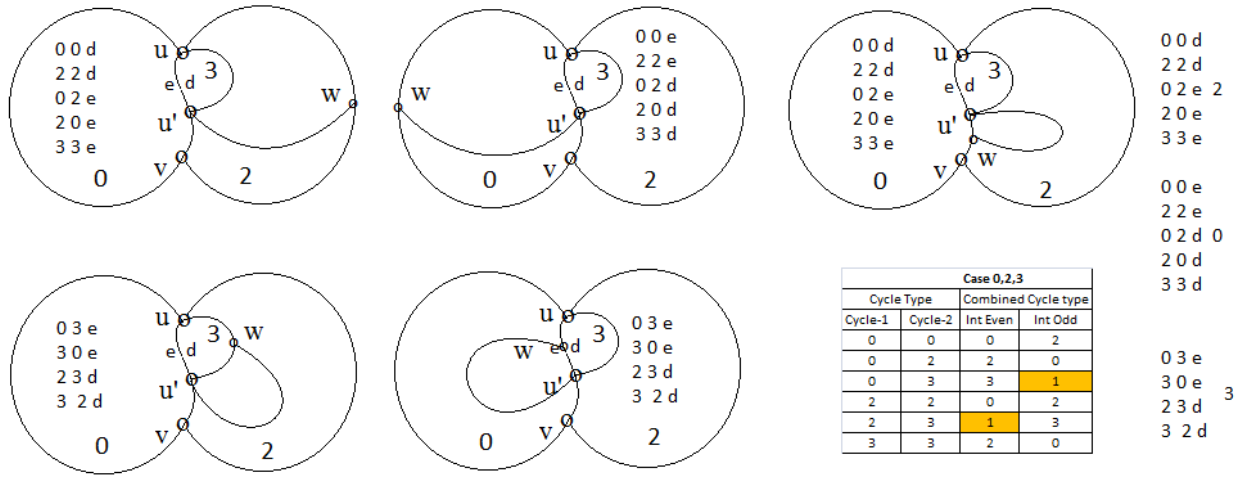

Fig.63 The five first level possibilities and cycle division for u'-w path for the case 0,2,3;e,d.

MGB-C3-case1-03201-19Oct16

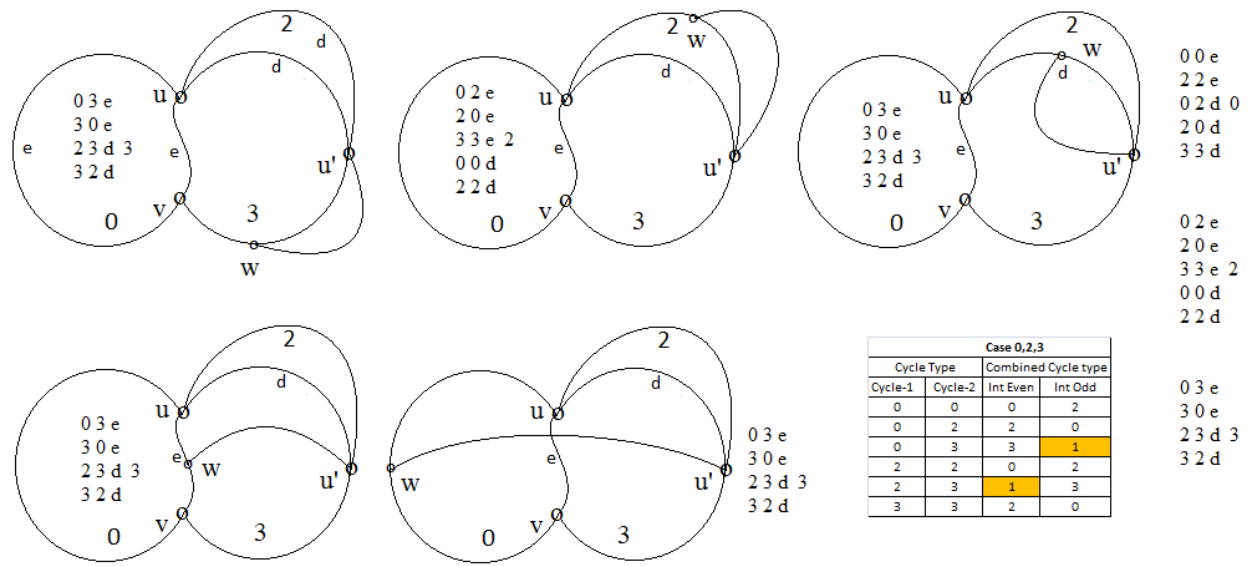

Fig.64 The five first level possibilities and cycle division for u'-w path in the case 0,3,2;e,d.

**Case (d)**   Same as Case c

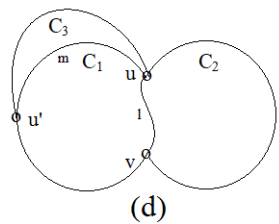

(d)

## Case (e)

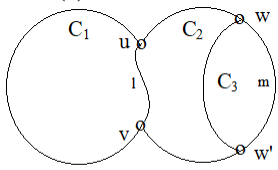

(e)

Table-17 Possible cases for the case (e)

| i | j | k | l | m | i+j-2l | j+k-2m | i+j+k-2l-2m |
|---|---|---|---|---|--------|--------|-------------|
| 0 | 2 | 3 | 0 | 1 | 2 | 3 | 3 |
| 0 | 3 | 2 | 0 | 1 | 3 | 3 | 3 |
| 2 | 0 | 3 | 1 | 0 | 0 | 3 | 3 |

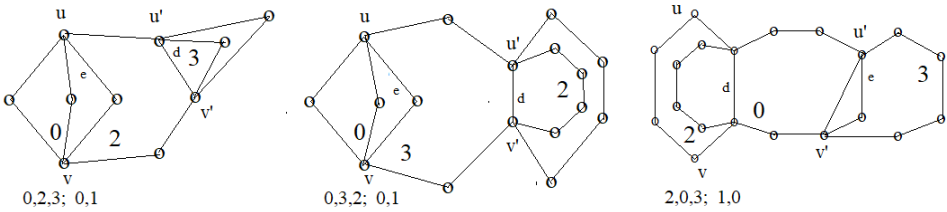

0,2,3; 0,1          0,3,2; 0,1          2,0,3; 1,0

Fig.65a Examples of Euler graphs from $\varepsilon_{023}$ for the cases in the Table-17.

MGB-Cijk-3CC-023-ed-20Jan16

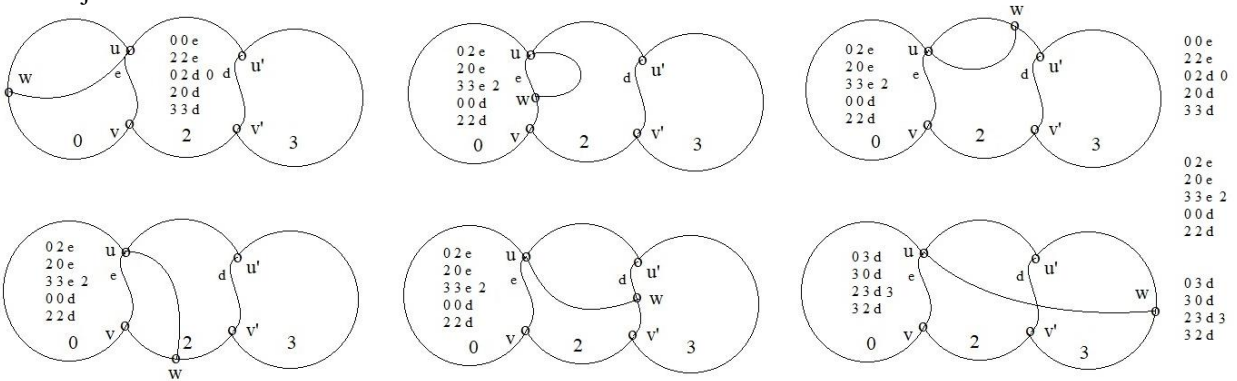

Fig.65b The six first level possibilities with cycle division in the case 0,2,3;e,d

**Case (f)** Graphs with three intersecting cycles, any two intersecting in different paths

Table-18 Possible cases for the case (f).

| i | j | k | l | m | n | i+j-2l | j+k-2m | i+k-2n | i+j+k-2l-2(m+n) |
|---|---|---|---|---|---|--------|--------|--------|-----------------|
| 0 | 2 | 3 | 0 | 1 | 0 | 2 | 3 | 3 | 3 |
| 0 | 3 | 2 | 0 | 1 | 0 | 3 | 3 | 2 | 3 |

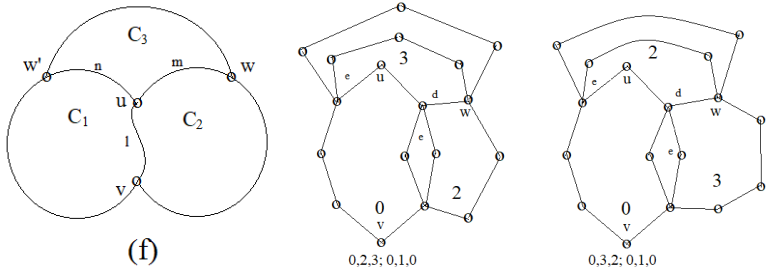

Fig.66 Examples of Euler graphs from $\varepsilon_{023}$ for the cases in the Table-18.

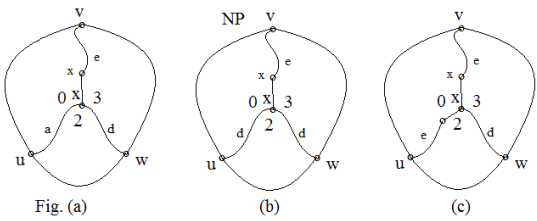

Fig.67 Possible intersection path parities. General case ead followed by the two subcases edd and eed.

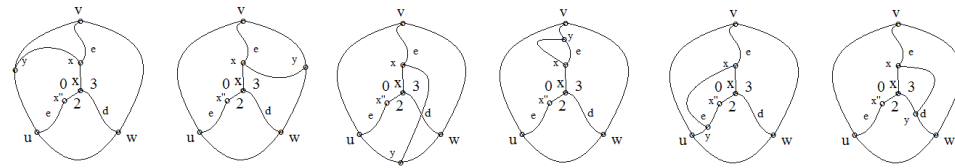

Fig.68 The six first level posibilities for x'-y path.

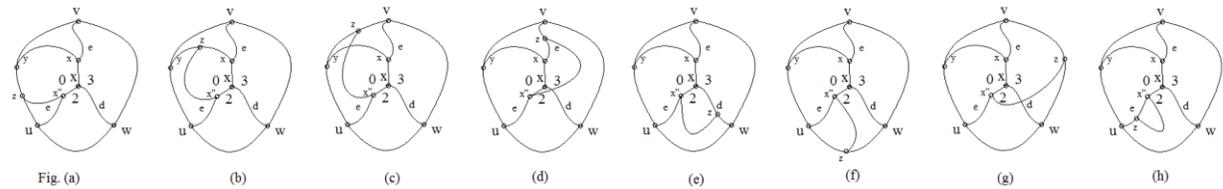

Fig.69 The eight posibilities for x''-z path for the first case of Fig.68

MGB-C023-eed-16Jul15

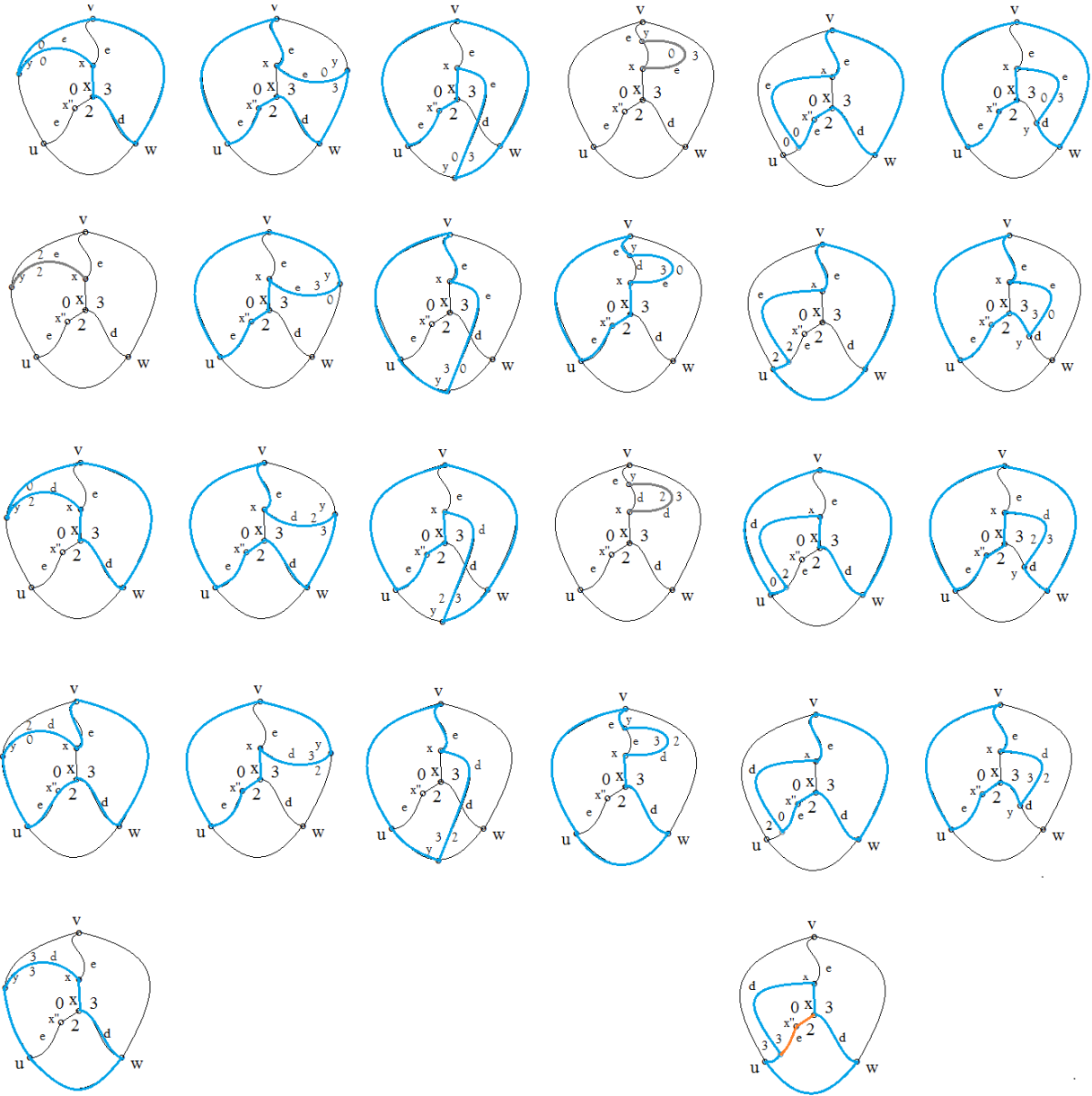

Fig.70 The 26 first level possibilities and stop cases with x'-y path.

MGB-C023-22-eed-16Jul15

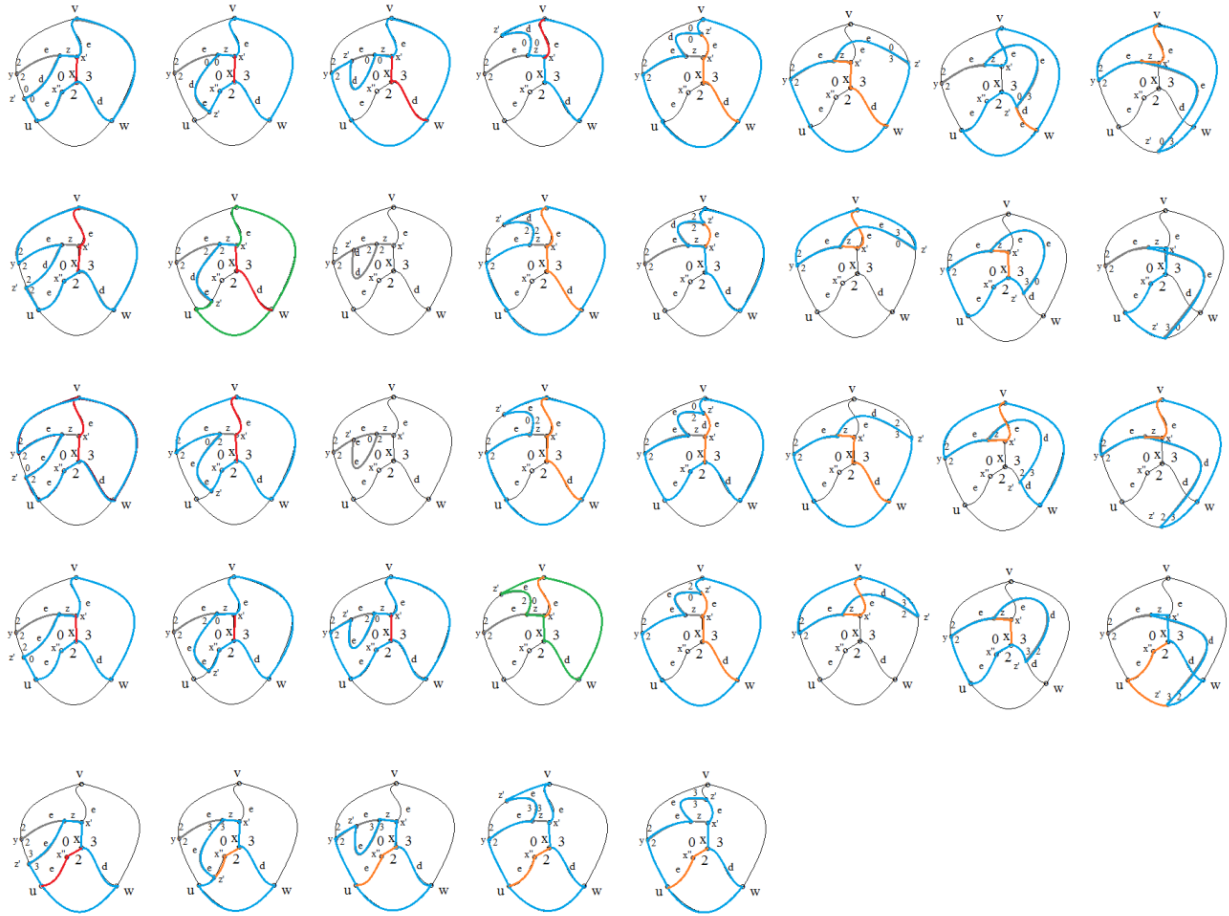

Fig.71 The 37 possibilities for z-z' path in the case 0,2,3;22-eed with stop cases.

Cases to continue to second level.

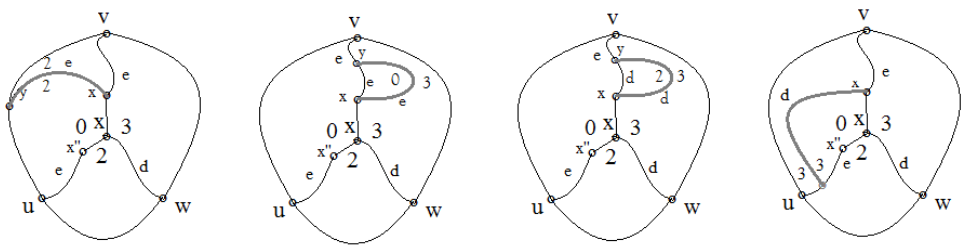

Fig.72 No solution cases

**Case (g)** Not studied.

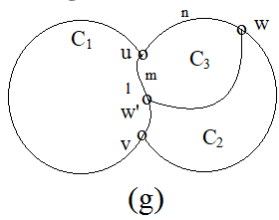

**Case (h)** Not studied.

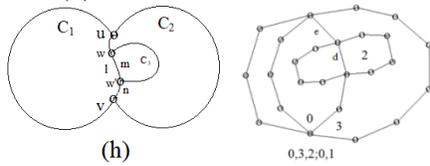

(h)

0,3,2;0,1

Fig.73 An example Euler graph from $\varepsilon_{023}$.

**Corollary 13.1.** A regular Euler graph in $\varepsilon_{023}$ is nonexistent in the cases which stop.

**Case 1,2&3. $\varepsilon_{123}$: Euler graphs with only three types of cycles $C_n$, $n \equiv 1,2\&3 \pmod 4$**

**Observation 8.** Order of graphs in $\varepsilon_{123}$ satisfies $p \geq 9$. Such a graph has at least one cycle of the type 1,2&3. That makes the graph of order at least six. Minimizing order implies that other cycles have maximum intersection with the 6-cycle. 5-cycle can have only a path of odd number of edges in common with the 6-cycle. This happens only when the 5-cycle has three edges in common with the 6-cycle, resulting in two nodes u,v of degree 3 with the node w not common with 6-cycle. u and v are of degree 3 and so there are nodes adjacent. Clearly they are not of 6-cycle. Two 3-cycles having uw and vw in common make the graph Euler and may be verified to be a minimum graph. Euler (9,12)-graph is shown in Fig.74b. Each block of Fig74c is a cycle graph of the type (1,2,3). It is not known whether the graphs in this case are nonplanar.

| Table-20   Case 1,2,3 | | | |
|---|---|---|---|
| Cycle Type | | Combined Cycle type | |
| Cycle-1 | Cycle-2 | Int Even | Int Odd |
| 1 | 1 | 2 | 0 |
| 1 | 2 | 3 | 1 |
| 1 | 3 | 0 | 2 |
| 2 | 2 | 0 | 2 |
| 2 | 3 | 1 | 3 |
| 3 | 3 | 2 | 0 |

CC rules in this case is given in the Table-20.

**Theorem 14.** The CC is of type 0 whenever two cycles of type 1,1 or 3,3 (1,3 or 2,2) of a graph G from $\varepsilon_{123}$ intersect in a path of length l>0 is odd (even).

**Theorem 15.** Size of a graph from $\varepsilon_{123}$ satisfies $q \equiv \xi_1 + 2\xi_2 + 3\xi_3 \pmod 4$.

Proof follows from Equation (1) (Rao 2014 [10]). Further, if $\xi_1 + 2\xi_2 + 3\xi_3 \equiv 1 \text{ or } 2 \pmod 4$ then the graphs are nongraceful by Rosa-Golumb criterion; else $\xi_1 + 2\xi_2 + 3\xi_3 \equiv 0 \text{ or } 3 \pmod 4$ and the graphs are candidates for gracefulness. Graphs in $\varepsilon_{123}$ satisfy that $\xi_1 > 0$, $\xi_2 > 0$ and $\xi_3 > 0$.

**Theorem 16.** Two necessary conditions follow: If $\xi_1 + 2\xi_2 + 3\xi_3 \equiv 0 \pmod 4$ then $\xi_1 + 3\xi_3$ is even. If $\xi_1 + 2\xi_2 + 3\xi_3 \equiv 3 \pmod 4$ then $\xi_1 + 3(\xi_3 - 1)$ is even.

**Conjecture 10**. Euler graphs in $\varepsilon_{123}$ satisfing $\xi_1 + 2\xi_2 + 3\xi_3 \equiv 0 \text{ or } 3 \pmod 4$ are graceful.

**Conjecture 11**. There exists a node of degree two in any graph of order p>8 from $\varepsilon_{123}$.

As a corollary it follows that:

**Conjecture 12**. Regular Euler graphs from $\varepsilon_{123}$ of order >8 are nonexistent.

or

A regular Euler graph from $\varepsilon_{123}$ of order p>8 with three types of cycles also contains fourth type.

Analysis follows as in the case 0,1,2.

**Theorem 17**. There exists a node of degree two in any graph of order p>8 from $\varepsilon_{123}$ in the cases which stop.

Part proof in the cases (a), (b), (c) and (f).

**Euler's Graph World - More Conjectures on Gracefulness Boundaries-III**

**Case (a)** Three cycles intersecting in a common path.

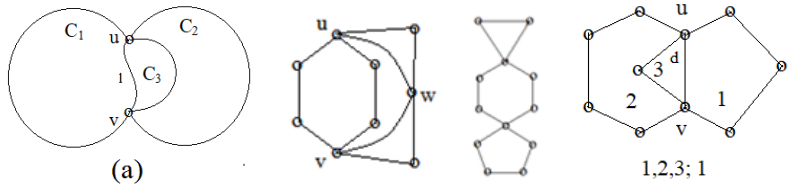

Fig.74 Case (a). Euler graphs from $\varepsilon_{123}$. Last graph is an example for the case in the Table-21.

Table-21 The only possibility in the case (a).

| i | j | k | l | i+j-2l | i+k-2l | j+k-2l |
|---|---|---|---|--------|--------|--------|
| 1 | 2 | 3 | 1 | 1 | 2 | 3 |

From Table-20, it follows that the three cycles of type i,j,k intersect in odd path. An example of such a graph is given in Fig.74d. Consider three cycles $C_1$, $C_2$, $C_3$ as in Fig.74a with common u-v path of odd length. This is done in two cases according as u-v-w path is of length one or more.

**Case-1** Common u-w path is an edge or u,w are adjacent.
MGB-CijkCP-123CPEdge-25May2020

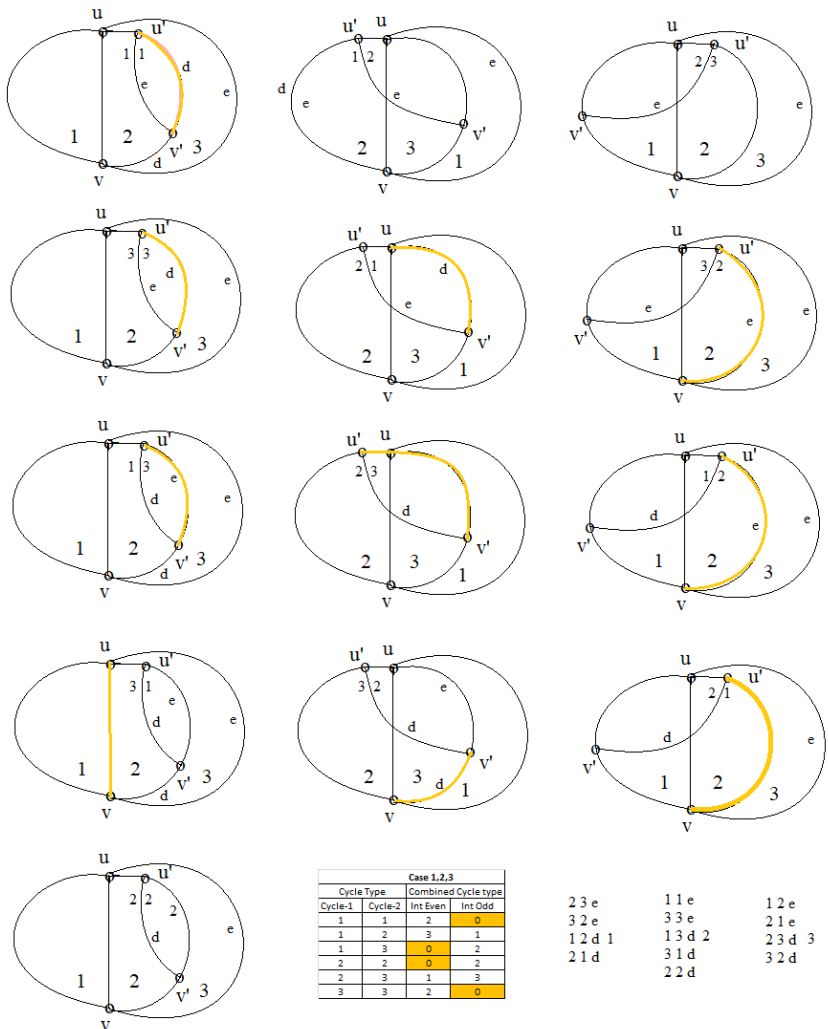

Fig.75 The 13 possibilities with common edge intersection and stop cases.

There is an edge uu'on the cycle 2. If u' is of degree two we are done. If not there are three possibilities for the u'-v' path as shown in Fig.75. From the rules of Table-20 we can enumerate the new cycles created by u'-v' path with parity for the intersection as shown in Fig.75 Some cases stop as shown in color.

**Case-2** Common u-w path is of odd length and three or more.

By assumption there is a node v such that uv is an edge on the u-w path. Note that v-w path is even. See Fig.76. Claim that v is of degree two. If not there is x so that v-x path exists. Three cases arise according as x is on u-x-w path as in row 1 of Figs.76a,b,c and v-v' path, v' on v-w path as in row 1 of Fig.76d.

MGB-CijkCP-123CPOdd-03Aug15

1,2,3 CP is Odd

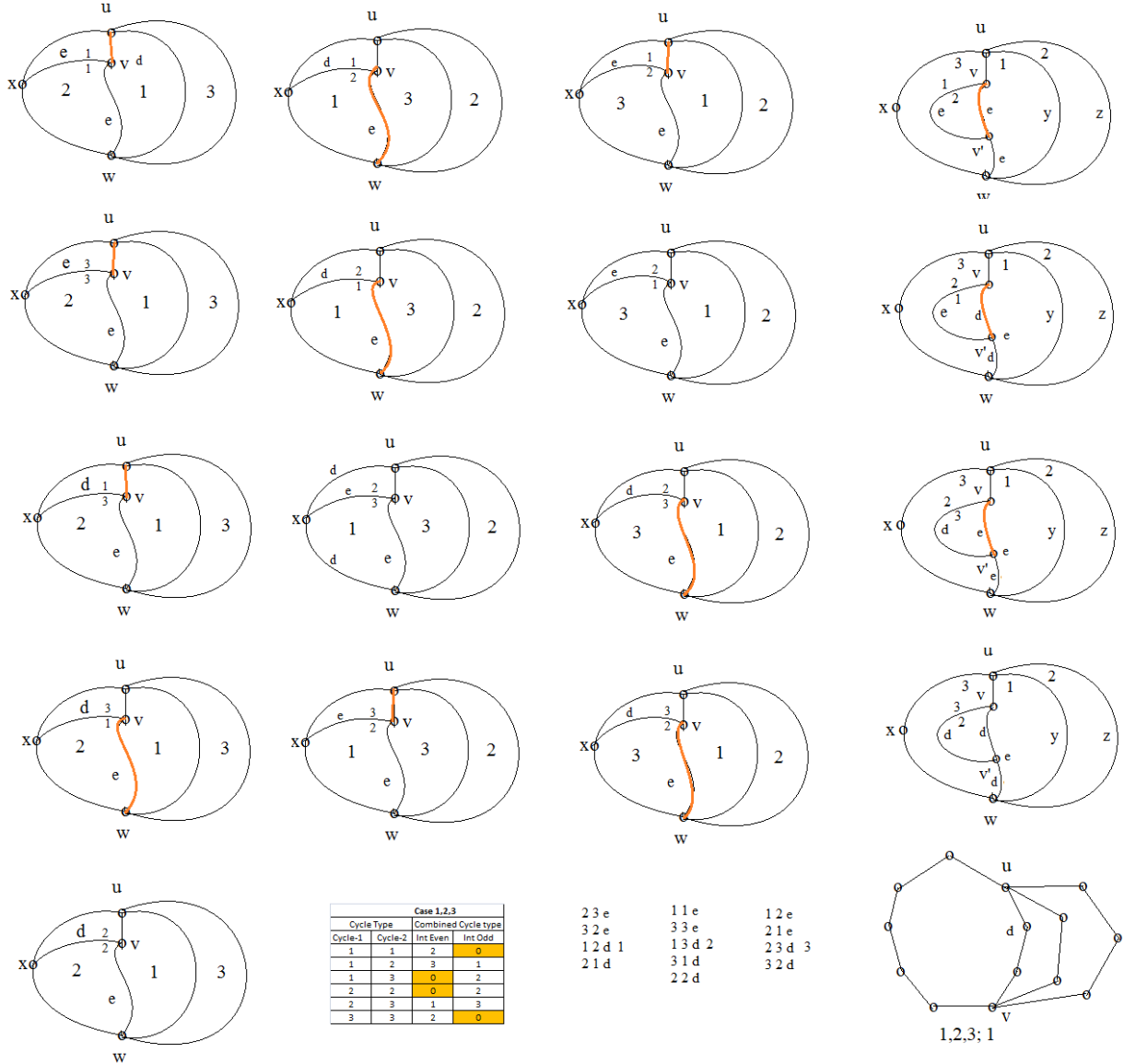

Fig.76 The 17 possibilities for v-x path and stop cases and an example of Euler graph from $\varepsilon_{123}$.

**Subcase-1** $C_2$ gets divided into two cycles: $C_{21}$:u-x-v,u; $C_{22}$:w-x-v-w. Five cases arise:

The case 1,1;e, with the cycles of type $T_1$ and $T_1$, respectively, is impossible as $C_{21}$ intersects with $C_1$ in odd path, a contradiction as it should be even, see Table-20.

The case 3,3;e, with the cycles of type $T_3$ and $T_3$, respectively, is impossible as $C_{12}$ intersects with $C_1$ in even path, a contradiction as it should be odd, see Table-20.

The case 1,3;d, with the cycles of type $T_1$ and $T_3$, respectively, is impossible as $C_{21}$ intersects with $C_1$ in odd path, a contradiction as it should be even, see Table-20.

The case 3,1;d, with the cycles of type $T_3$ and $T_1$, respectively, is impossible as $C_{22}$ intersects with $C_3$ in even path, a contradiction as it should be odd, see Table-20.

The case 2,2;d, with the cycles of type $T_2$ and $T_2$, respectively, is impossible as $C_{11}$ intersects with $C_3$ in odd path, a contradiction as it should be even, see Table-20.

**Subcase-2** $C_1$ gets divided into two cycles: $C_{11}$:u-x-v,u; $C_{12}$:w-x-v-w.  Four cases arise:

The case 1,2;d, with the cycles of type $T_1$ and $T_2$, respectively, is impossible as $C_{12}$ intersects with $C_2$ in even path, a contradiction as it should be odd, see Table-20.

The case 2,1;d, with the cycles of type $T_2$ and $T_1$, respectively, is impossible as $C_{12}$ intersects with $C_2$ in even path, a contradiction as it should be odd, see Table-20.

The case 2,3;e, with the cycles of type $T_2$ and $T_3$, respectively, is impossible as $C_{12}$ intersects with $C_2$ in even path, a contradiction as it should be odd, see Table-20.

The case 3,2;e, with the cycles of type $T_3$ and $T_2$, respectively, is impossible as $C_{11}$ intersects with $C_3$ in odd path, a contradiction as it should be even, see Table-20.

**Subcase-3** $C_3$ gets divided into two cycles: $C_{31}$:u-x-v,u; $C_{32}$:w-x-v-w. Four cases arise:

The case 1,2;e, with the cycles of type $T_1$ and $T_2$, respectively, is impossible as $C_{31}$ intersects with $C_1$ in odd path, a contradiction as it should be even, see Table-20.

The case 2,1;e, with the cycles of type $T_2$ and $T_1$, respectively, is impossible as $C_{31}$ intersects with $C_2$ in even path, a contradiction as it should be odd, see Table-20.

The case 2,3;d, with the cycles of type $T_2$ and $T_3$, respectively, is impossible as $C_{32}$ intersects with $C_1$ in even path, a contradiction as it should be odd, see Table-20.

The case 3,2;d, with the cycles of type $T_3$ and $T_2$, respectively, is impossible as $C_{32}$ intersects with $C_2$ in odd path, a contradiction as it should be even, see Table-20.

**Subcase-4** C3 gets divided into two cycles: $C_{31}$:u-x-w-v'-v,u; C32:v-v'-v. Four cases arise:

The case 1,2;e, with the cycles of type $T_1$ and $T_2$, respectively, is impossible as $C_{32}$ intersects with $C_2$ in even path, a contradiction as it should be odd, see Table-20.

The case 2,1;e, with the cycles of type $T_2$ and $T_1$, respectively, is impossible as $C_{32}$ intersects with $C_1$ in odd path, a contradiction as it should be even, see Table-20.

The case 2,3;d, with the cycles of type $T_2$ and $T_3$, respectively, is impossible as $C_{32}$ intersects with $C_1$ in even path, a contradiction as it should be odd, see Table-20.

The case 3,2;d, with the cycles of type $T_3$ and $T_2$, respectively, is impossible as $C_{32}$ intersects with $C_2$ in odd path, a contradiction as it should be even, see Table-20.

**Case (b)**

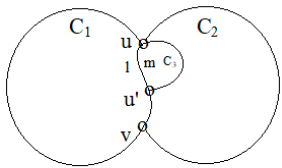

(b)

Table-22 Possible cases for the case (b).

| i | j | k | l | m | i+j-2l | i+k-2m | j+k-2m |
|---|---|---|---|---|---|---|---|
| 1 | 2 | 3 | 0 | 1 | 3 | 2 | 3 |
| 1 | 2 | 3 | 1 | 1 | 1 | 2 | 3 |
| 1 | 3 | 2 | 1 | 0 | 2 | 3 | 1 |
| 1 | 3 | 2 | 1 | 1 | 2 | 1 | 3 |
| 2 | 3 | 1 | 0 | 1 | 1 | 1 | 2 |
| 2 | 3 | 1 | 1 | 1 | 3 | 1 | 2 |

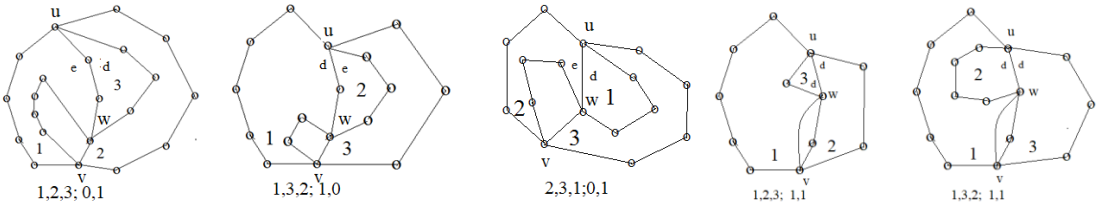

Fig.77 Examples of Euler graphs from $\varepsilon_{123}$ for the first five cases in the Table-22.

MGB-C3-case2-12301-19Oct16

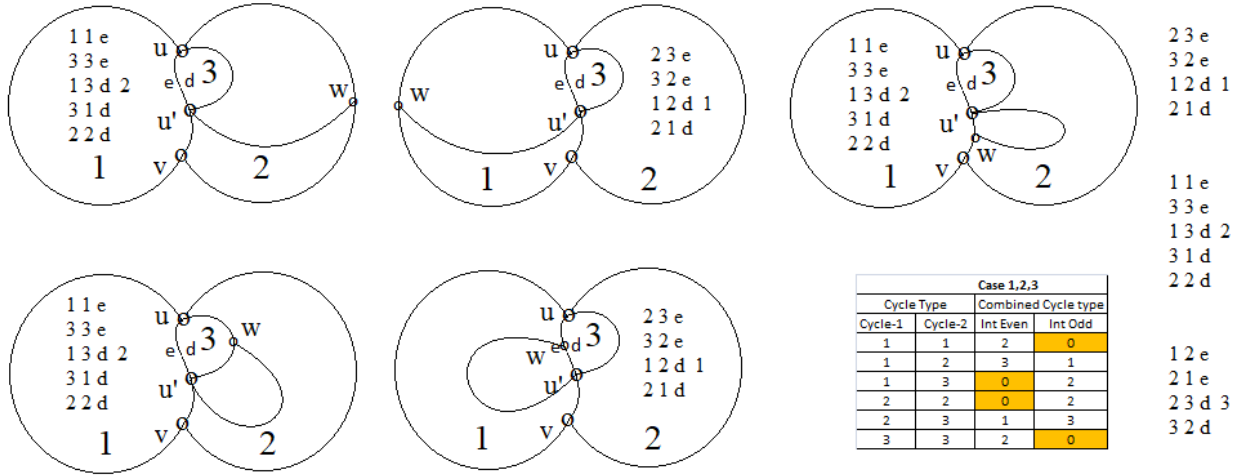

Fig.78 The five possible u'-w paths in the case 1,2,3;e,d and division of a cycle in which the path lies.

**Case (c)**

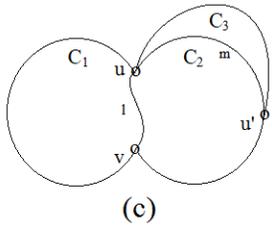

(c)

Table-23 Possible case in the case (c).

| i | j | k | l | m | i+j-2l | j+k-2m | i+j+k-2l-2m |
|---|---|---|---|---|--------|--------|-------------|
| 1 | 2 | 3 | 0 | 0 | 3 | 1 | 2 |
| 1 | 2 | 3 | 1 | 1 | 1 | 3 | 2 |
| 1 | 3 | 2 | 1 | 1 | 2 | 3 | 2 |
| 2 | 1 | 3 | 1 | 1 | 1 | 2 | 2 |

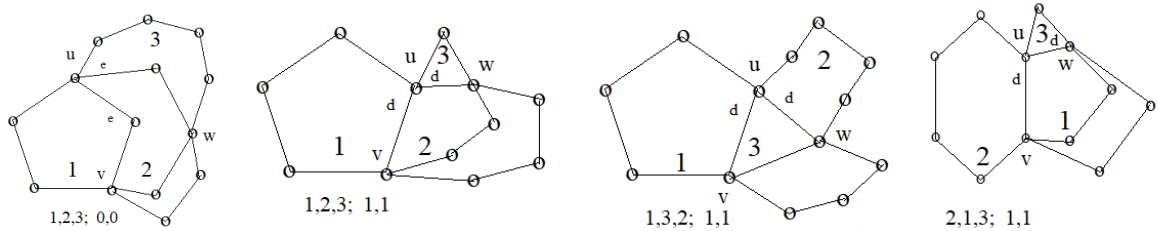

Fig.79 Examples of Euler graphs from $\varepsilon_{123}$ for the cases in the Table-23.

MGB-C3-case1-12300-19Oct16

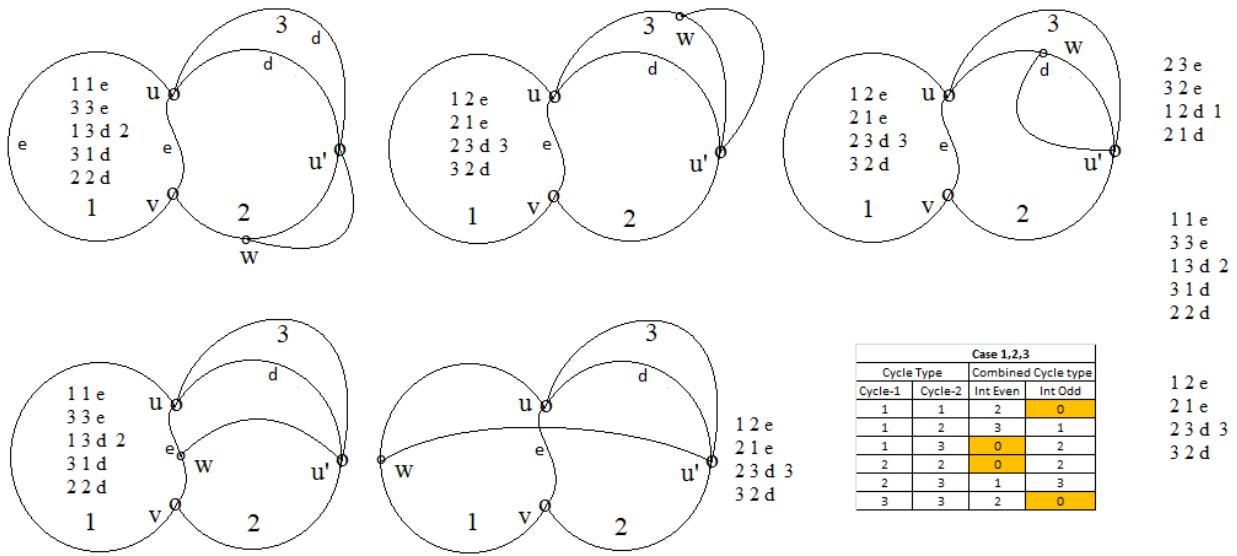

Fig.80 The five possibilities for u'-w path in the case 1,2,3;e,e and division of cycle in which the path lies.

**Case (d)**

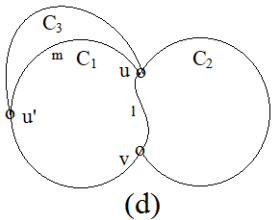

(d)

Same as Case c.

**Case (e)**

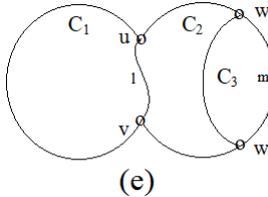

(e)

Table-24 Possible case in the case (c).

| i | j | k | l | m | i+j-2l | j+k-2m | i+j+k-2l-2m |
|---|---|---|---|---|--------|--------|-------------|
| 1 | 2 | 3 | 0 | 0 | 3 | 1 | 2 |
| 1 | 2 | 3 | 1 | 1 | 1 | 3 | 2 |
| 1 | 3 | 2 | 1 | 1 | 2 | 3 | 2 |
| 2 | 1 | 3 | 1 | 1 | 1 | 2 | 2 |

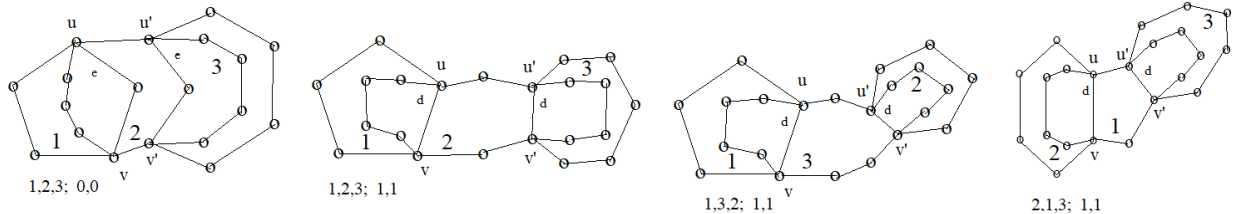

Fig.81a Examples of Euler graphs from $\varepsilon_{123}$ for the cases in the Table-24.

MGB-Cijk-3CC-123-ee-20Jan16

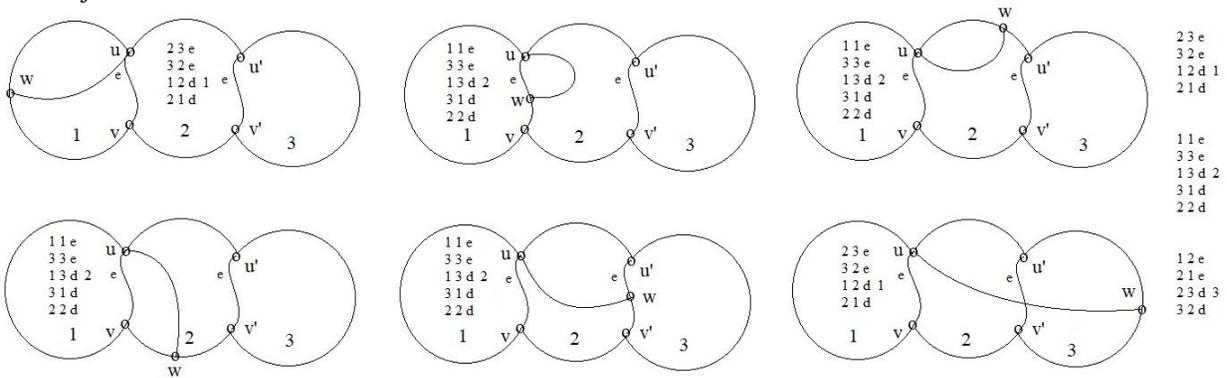

Fig.81b The six possibilities for u-w path with cycle division in the case 1,2,3;e,e.

**Case (f)**  Graphs with three intersecting cycles, any two intersecting in different path

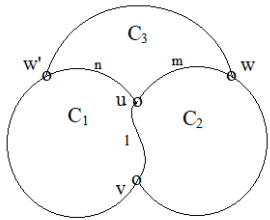

(f)

Table-25 Possible case in the case (c).

| i | j | k | l | m | n | i+j-2l | j+k-2m | i+k-2n | i+j+k-2l-2(m+n) |
|---|---|---|---|---|---|--------|--------|--------|------------------|
| 1 | 2 | 3 | 0 | 1 | 1 | 3 | 3 | 2 | 2 |
| 1 | 2 | 3 | 1 | 0 | 1 | 1 | 1 | 2 | 2 |
| 1 | 3 | 2 | 1 | 0 | 1 | 2 | 1 | 1 | 2 |

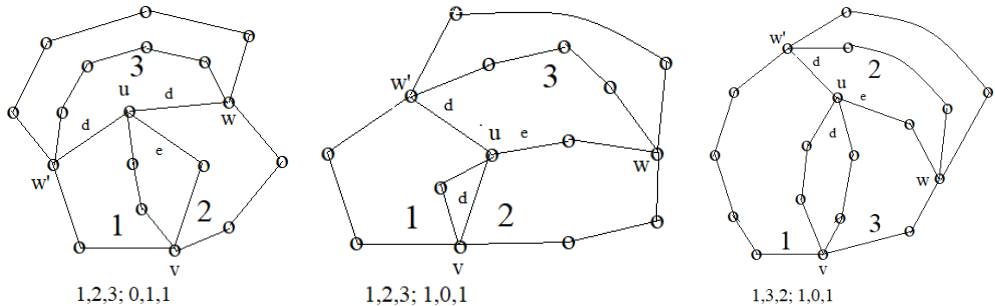

1,2,3; 0,1,1        1,2,3; 1,0,1        1,3,2; 1,0,1

Fig.82 Examples of Euler graphs from $\varepsilon_{123}$ for the cases in the Table-25.

MGB-C123-12-1-ded-15Sep15

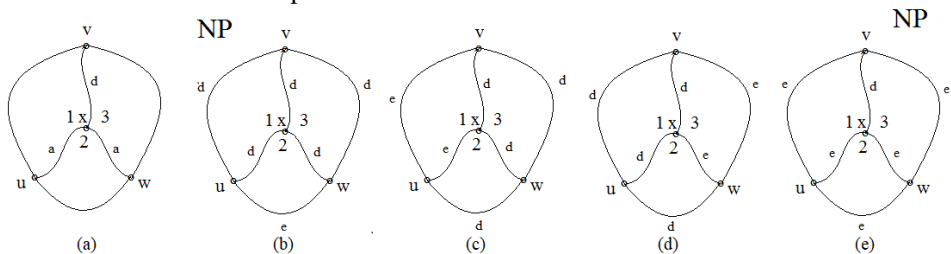

Fig.83 General Case (a) daa leading to four cases ddd, ded, dde and dee with (b) and (e) not possible.

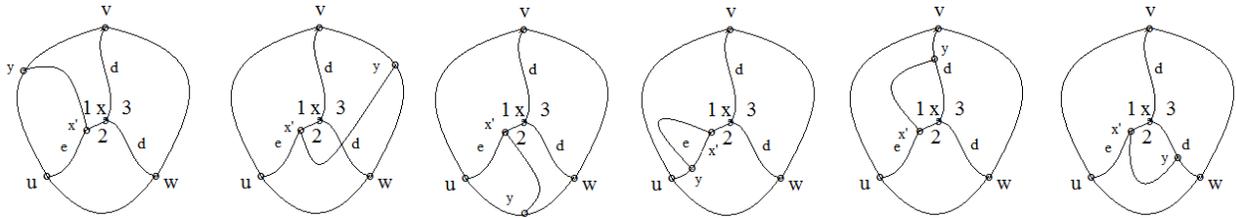

Fig.84 Case ded with possible x'-y paths.

The x'-y path divides the cycle in which it falls. Using intersection rules each of these cases lead to subcases as shown in Fig.85 with stop cases. We proceed with the cases which continue. For example, first case in the first row of Fig.85.

MGB-C123-ded-15Sep15

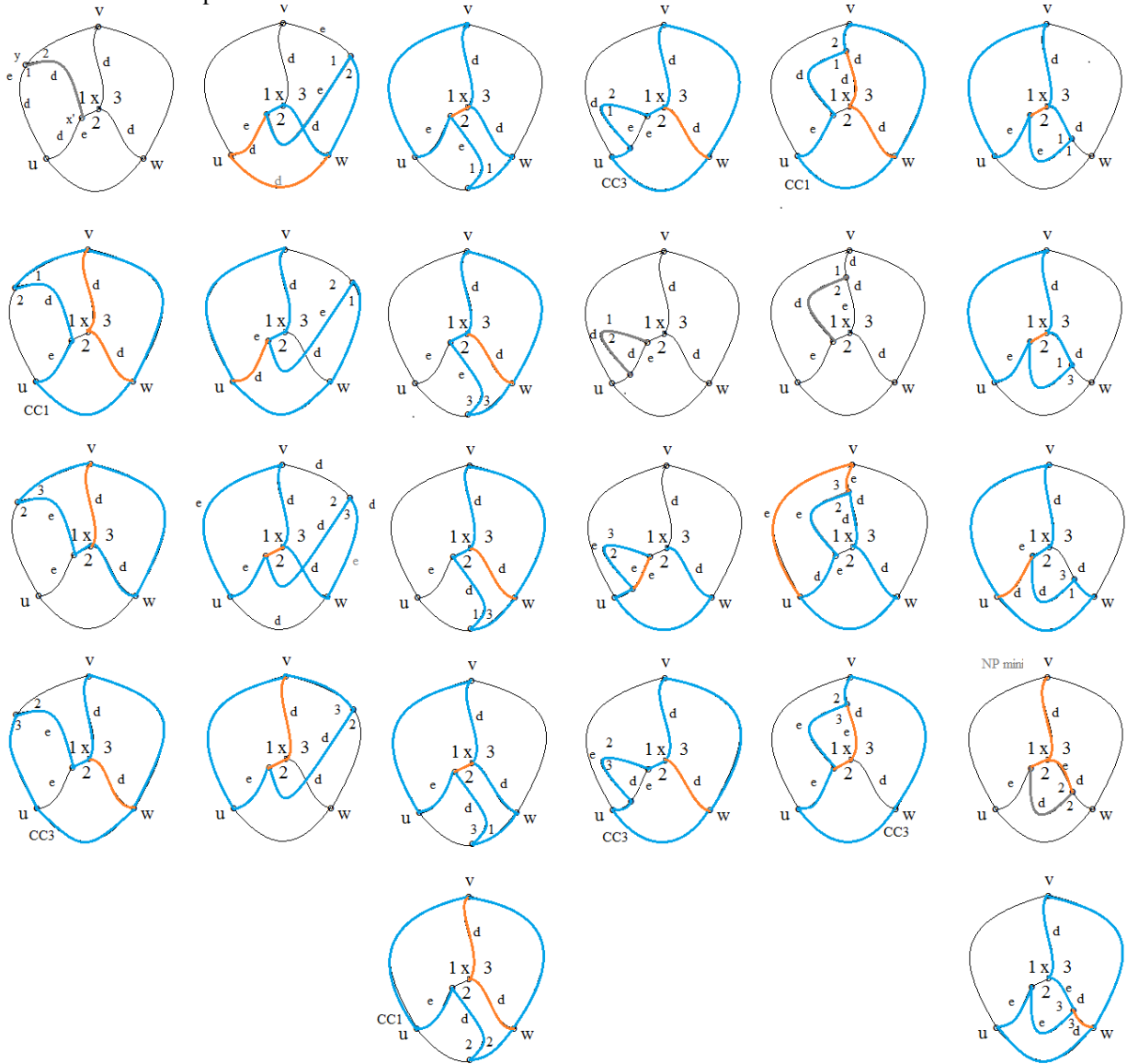

Fig.85 The 26 possible first level cases with cycle division in the case 1,2,3-ded and stop cases.

**Euler's Graph World - More Conjectures on Gracefulness Boundaries-III**

The first case 2,1;d in row 1 of Fig.85 may be continued if there is a node with degree >1. If such a node is degree two then we are through else continue with a path. Possible cases and the cycle division with stop cases are shown in Fig.86.

Case 1,2;d fifth case in row 2 of Fig.84 continues. Consider a node adjacent to x. If this node is degree two we are done else we continue. Fig.87 shows possible cases for the path from a node on x-v path with stop cases.

MGB-C123-12-1-ded-15Sep15

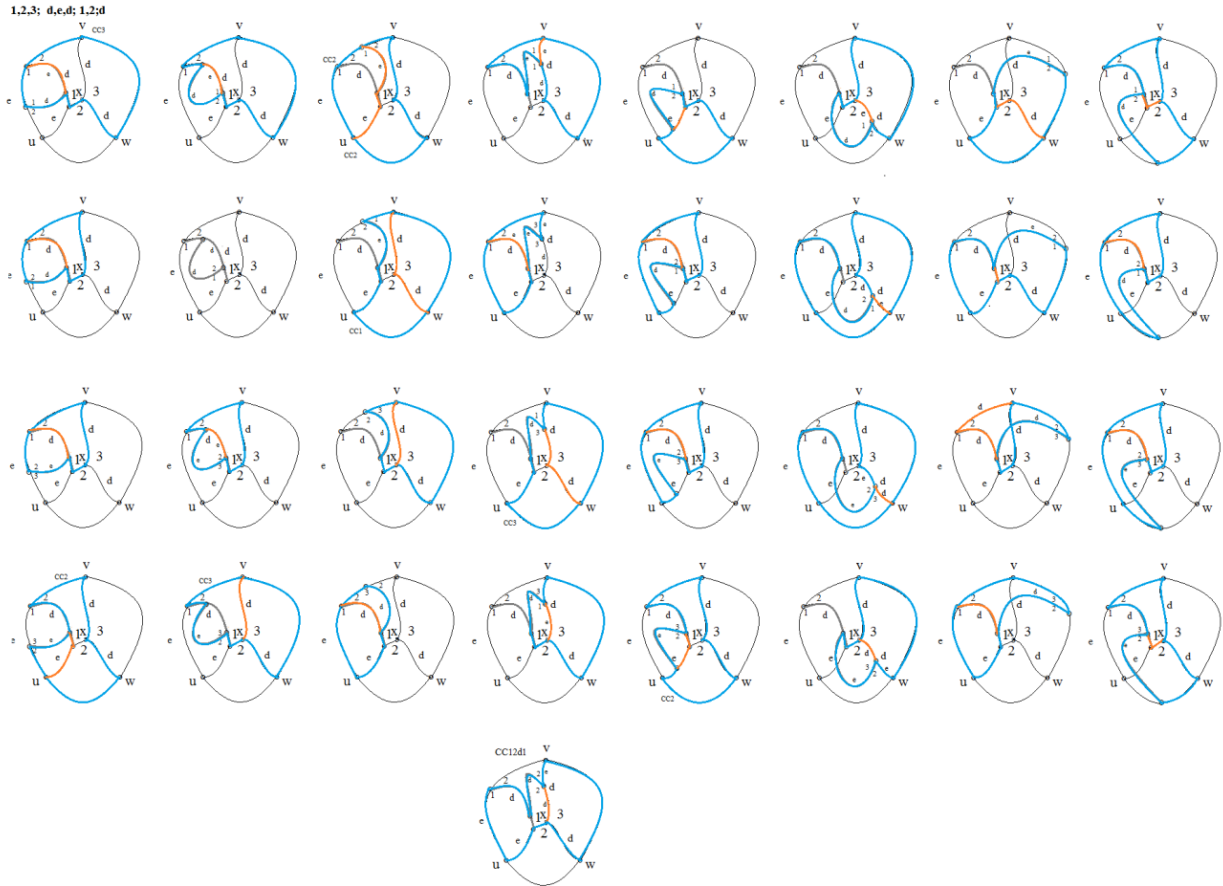

Fig.87 The 33 second level possibilities of the case 1,2,3-ded-12d with cycle division and stop cases.

MGB-C123-12-2-ded-15Sep15

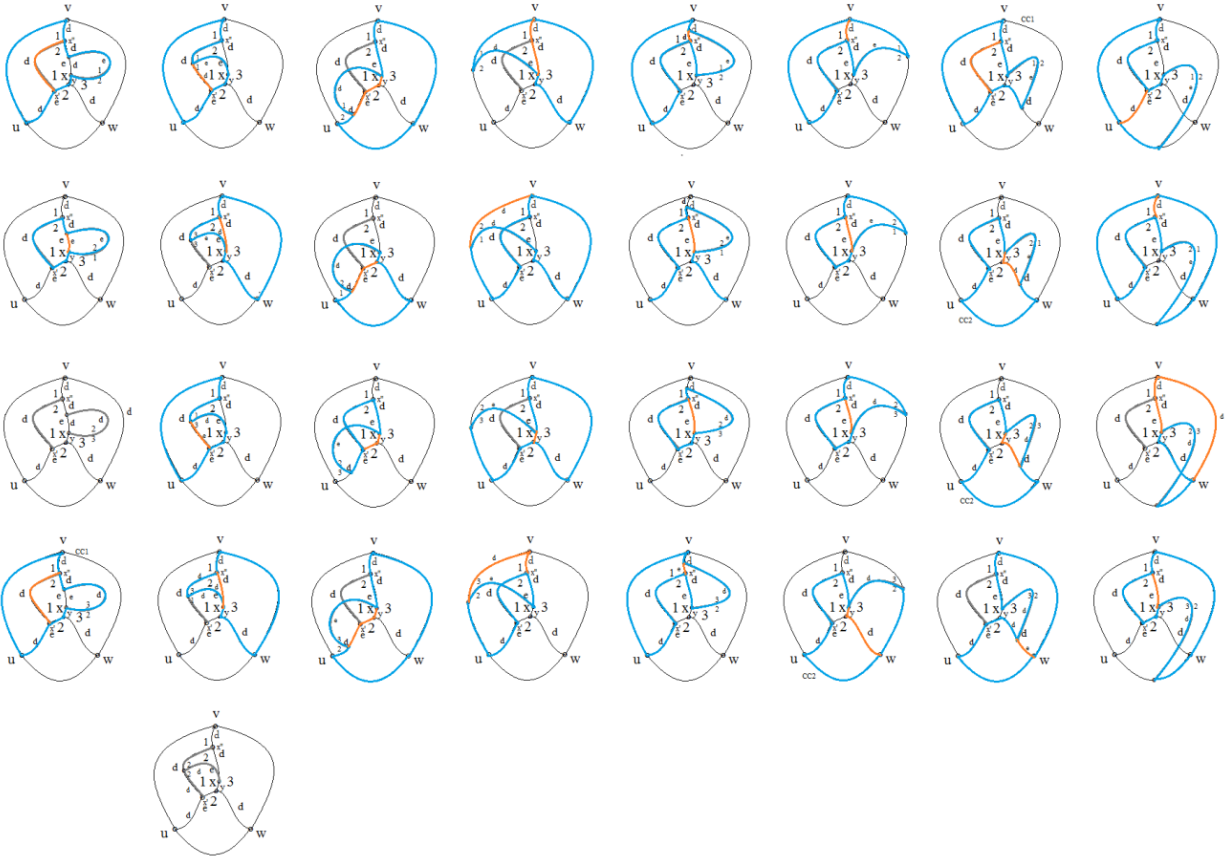

Fig.88 The 33 possibilities of the case 123-ded-12d with cycle division and stop cases.

Same steps are followed in the case dde in Figs.89 to 94.

MGB-C123-dde-15Sep15

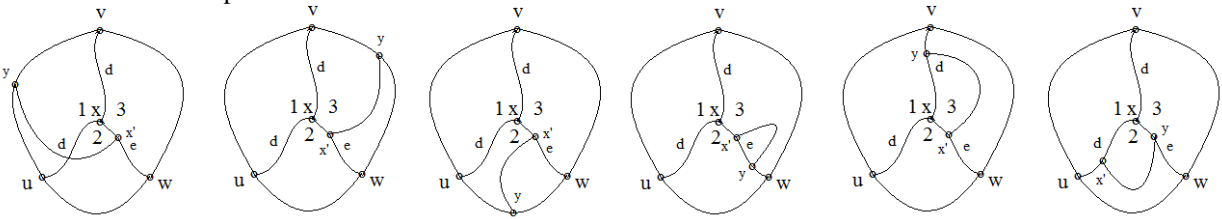

Fig.89 The six possibilities of the x'-y path in the case 1,2,3-dde.

MGB-C123-dde-15Sep15

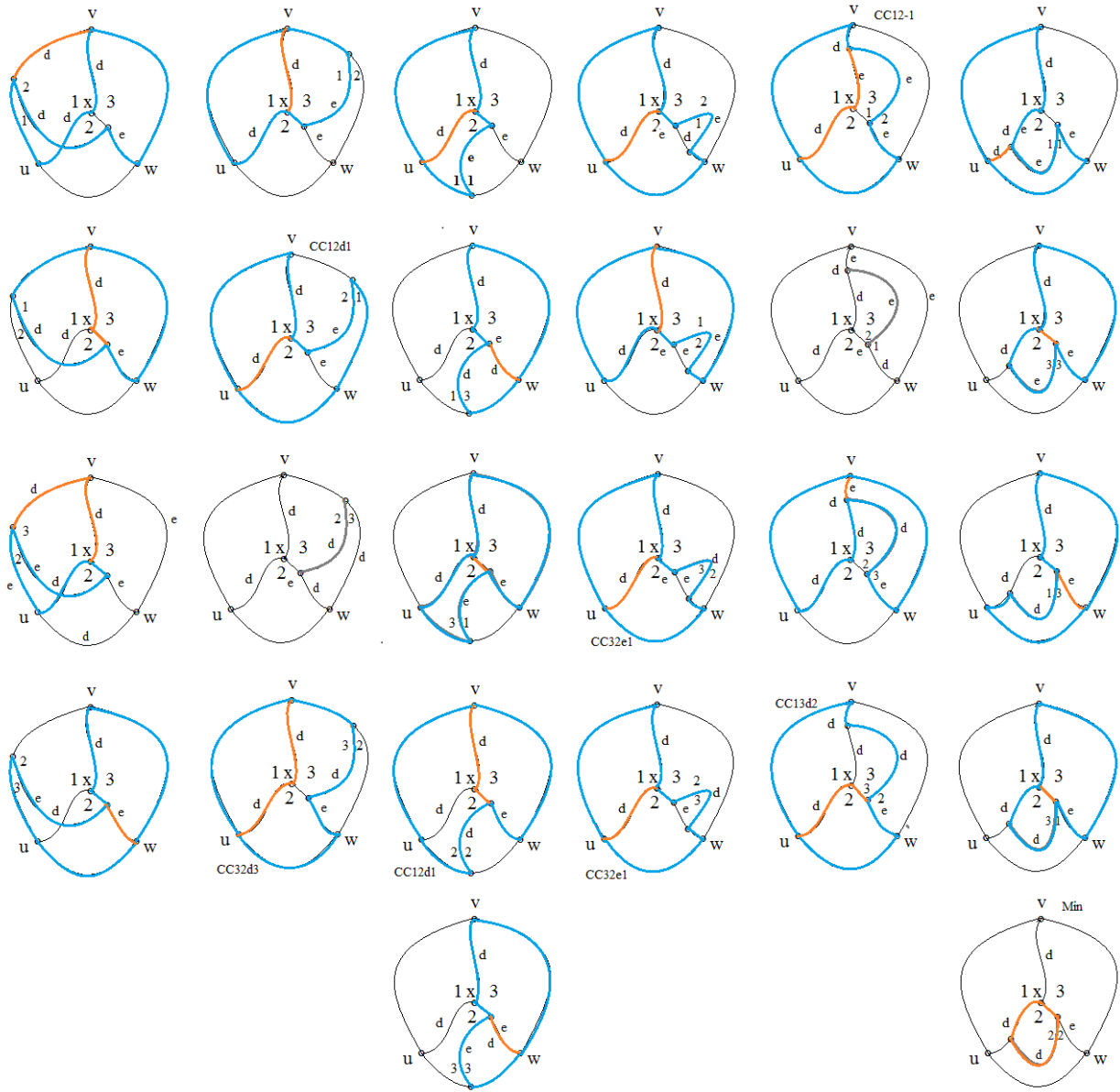

Fig. 90 The 26 first level possibilities with cycle division and stop cases for the case 1,2,3;dde

MGB-C123-21e-dde-15Sep15

1,2,3; d,d,e;   2,1;e

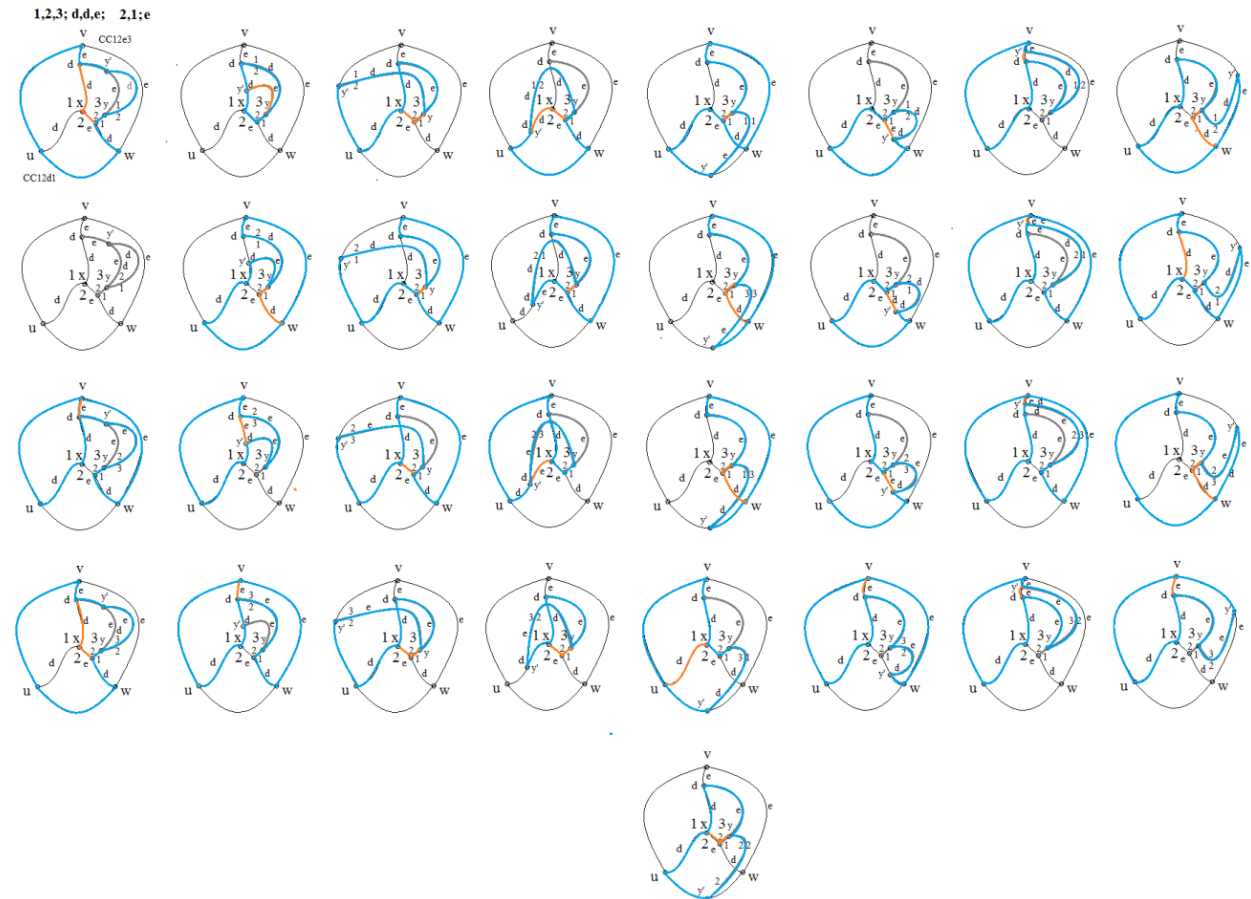

Fig.91 The 33 second level possibilities for 1,2,3-dde-21e case with division of cycles and stop cases.

MGB-C123-23d-Comp-dgt1-dde-15Sep15

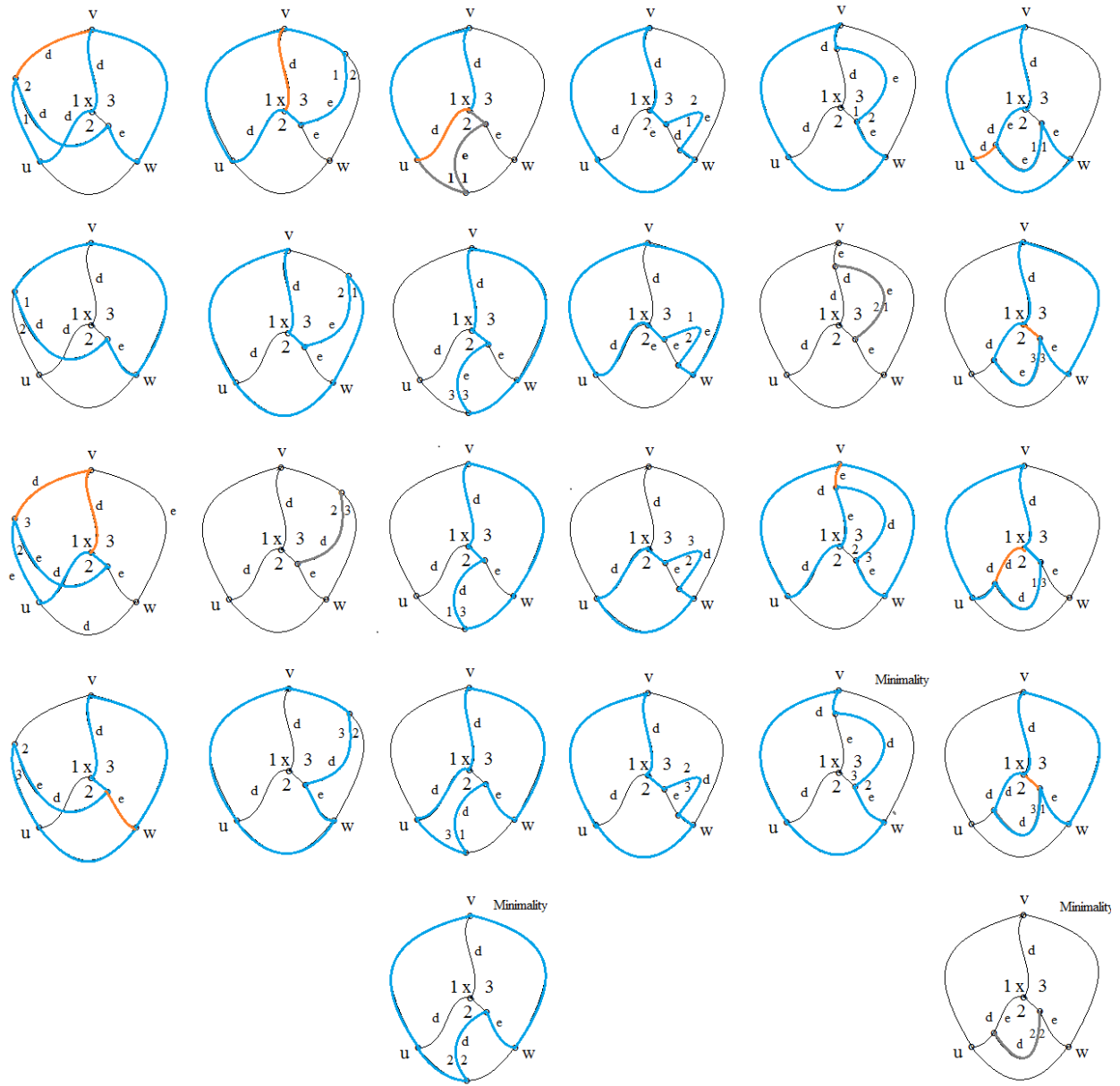

Fig.92 The 26 first level possibilities of case 123-dde-23d with division of cycles and stop cases.

MGB-C123-23d-Comp-dgt1-dde-15Sep15

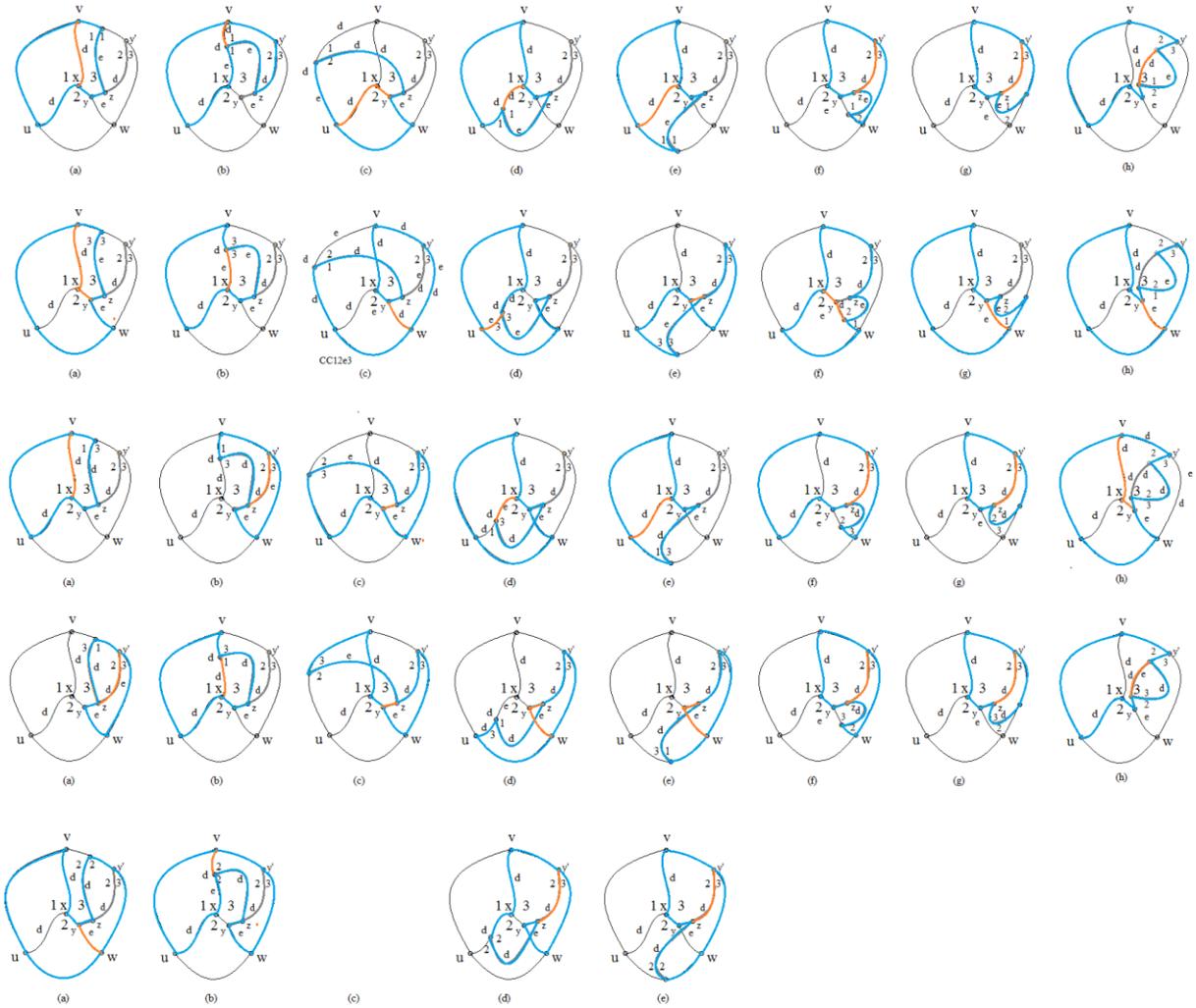

Fig.93 The 36 second level possibilities of case C123-dde-23d with division of cycles. Complete- All cases stop.

Cases to continue to second level.

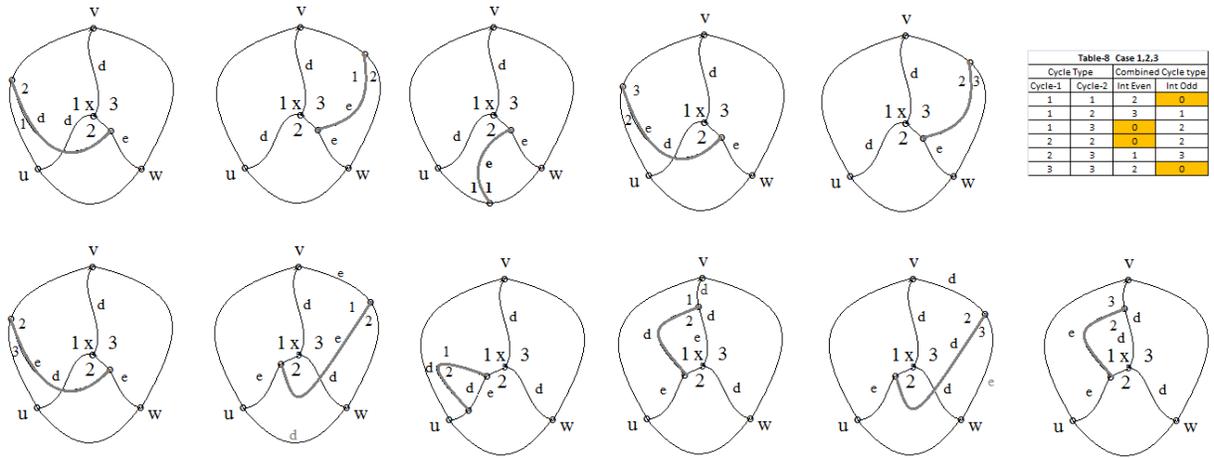

Fig.94 No solution cases

**Case (g)** Not studied.

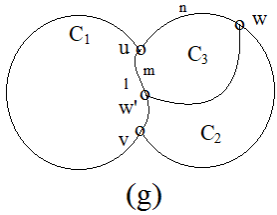

(g)

**Case (h)** Not studied.

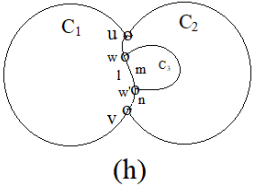

(h)

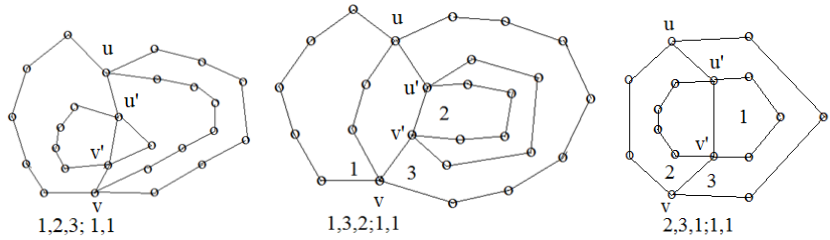

Fig.95 Examples of Euler graphs from $\varepsilon_{123}$.

**Corollary 17.1.** A regular Euler graph in $\varepsilon_{123}$ is nonexistent in the stop cases.

**Euler Graphs with All Four Types of Cycles**

Order of graphs satisfies p≥7.

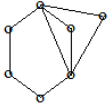

Fig.96 A Euler (7,9)-graph.

**Theorem 18.** Size of Euler graph from $\varepsilon_{0123}$ satisfies q≡$\xi_1$+2$\xi_2$+3$\xi_3$ (mod 4).

Equation (1) (Rao 2014 [10]) reduces to: q≡1$\xi_1$+2$\xi_2$+3$\xi_3$ (mod 4). If $\xi_1$+2$\xi_2$+3$\xi_3$≡0or3(mod 4) is satisfied then the graphs are candidates for gracefulness. If $\xi_1$+2$\xi_2$+3$\xi_3$≡1or2(mod 4) then the graphs are nongraceful by Rosa-Golumb criterion. Graphs in $\varepsilon_{0123}$ satisfy that $\xi_0$>0, $\xi_1$>0, $\xi_2$>0 and $\xi_3$>0.

**Theorem 19.** Two necessary conditions follow: If $\xi_1$+2$\xi_2$+3$\xi_3$≡0(mod 4) then $\xi_1$+3$\xi_3$ is even. If $\xi_1$+2$\xi_2$+3$\xi_3$≡3(mod 4) then $\xi_1$+3($\xi_3$-1) is even or $\xi_1$+3$\xi_3$ is odd.

**Conjecture 13**. Euler graphs in $\varepsilon_{0123}$ satisfying $\xi_1$+2$\xi_2$+3$\xi_3$≡0or3(mod 4) are graceful.

**Summary**

In support of the general conjecture, Conjecture-1 (Rao (2014), [9]) that *graphforest of a graceful pendant free graph is graceful,* we have:

**Conjecture 14**. Eulerforest of a Euler graph in Conjectures 1, 4. 7, 10 and 13 is graceful.

**Problems**

Study Euler graphs
- For extremal graphs
- Under regularity
- Under planarity

Enumerate Euler Graph types under (mod 4).

## Acknowledgements

The author, an Oil & Gas Professional with works in Graph Theory, expresses deep felt gratitude to
The *Sonangol Pesquisa e Produção*, Luanda, Angola
for excellent facilities, support and encouragement (2011-'13) with special thanks to
Mr. *Joao Noguiera*, Managing Director, Development Directorate.

Deep felt gratitude to *Alexander Rosa*, Emeritus Professor,
Dept. of Mathematics & Statistics, McMaster University, West Hamilton, Ontario, Canada and
*Solomon W Golumb*, an American Mathematician, Engineer and Professor of Electrical Engineering
at the University of Southern California (A recipient of the USC Presidential Medallion,
the IEEE Shannon Award of the Information Theory Society and three honorary doctorate degrees,
National Medal of Science presented by President Barack Obama)
for their initial inspiring works on Graceful Labelings;
to Emeritus Professor *Lanka Radhakrishna*, Department of Mathematics, Shivaji University, Kolhapur;
to Emeritus Professor *G.A. Patwardhan*, Combinatorics, Department of Mathematics, IIT Bombay;
and to Executive Director, Mr. *J.L. Narasimham*, ONGC, Mumbai